\documentclass[11pt,journal,comsoc]{article}
\usepackage[latin9]{inputenc}
\usepackage{geometry}
\usepackage{url}
\usepackage{amsmath}
\usepackage{amssymb}
\usepackage{stmaryrd}

\makeatletter

\providecommand{\tabularnewline}{\\}

\date{}

\usepackage{amsthm}\usepackage{bm}\usepackage{mathrsfs}\usepackage{enumitem}

\usepackage{comment}

\usepackage{tikz}
\usepackage{pgfplots}

\newtheorem{assumption}{Assumption}
\newtheorem{lemma}{Lemma}
\newtheorem{definition}{Definition}
\newtheorem{theorem}{Theorem}
\newtheorem{corollary}{Corollary}
\newtheorem{remark}{Remark}

\makeatother

\begin{document}
\title{Phase transitions in nonparametric regressions}
\author{Ying Zhu\\
First draft on arXiv: December 7, 2021. This draft: October 4, 2023\thanks{Assistant Professor of Economics at UC San Diego. yiz012@ucsd.edu.
I thank the co-editor Elie Tamer, the anonymous AE and two anonymous
referees for their valuable comments. I would also like to thank Xiaohong
Chen, Fang Han, Hidehiko Ichimura, Oliver Linton and Yixiao Sun for
substantial comments. I am also grateful to Don Andrews, Matias Cattaneo,
Julie Cullen, Bo Honoré, Michael Jansson, Jason Klusowski, Shakeeb
Khan, Esfandiar Maasoumi, Ulrich Müller, Mikkel Plagborg-Moller, Dimitris
Politis, Chris Sims, and Zhijie Xiao. I am thankful to the Society
of Hellman Fellows at University of California and the Cowles Foundation
at Yale University, and appreciate everyone who attended my talks
at a number of seminars and conferences.}}
\maketitle
\begin{abstract}
When the unknown regression function of a single variable is known
to have derivatives up to the $(\gamma+1)$th order bounded in absolute
values by a common constant everywhere or a.e. (i.e., $(\gamma+1)$th
degree of smoothness), the minimax optimal rate of the mean integrated
squared error (MISE) is stated as $\left(\frac{1}{n}\right)^{\frac{2\gamma+2}{2\gamma+3}}$
in the literature. This paper shows that: (i) if $n\leq\left(\gamma+1\right)^{2\gamma+3}$,
the minimax optimal MISE rate is $\frac{\log n}{n\log(\log n)}$ and
the optimal degree of smoothness to exploit is roughly $\max\left\{ \left\lfloor \frac{\log n}{2\log\left(\log n\right)}\right\rfloor ,\,1\right\} $;
(ii) if $n>\left(\gamma+1\right)^{2\gamma+3}$, the minimax optimal
MISE rate is $\left(\frac{1}{n}\right)^{\frac{2\gamma+2}{2\gamma+3}}$
and the optimal degree of smoothness to exploit is $\gamma+1$.

The fundamental contribution of this paper is a set of metric entropy
bounds we develop for smooth function classes. Some of our bounds
are original, and some of them improve and/or generalize the ones
in the literature (e.g., Kolmogorov and Tikhomirov, 1959). Our metric
entropy bounds allow us to show phase transitions in the minimax optimal
MISE rates associated with some commonly seen smoothness classes as
well as non-standard smoothness classes, and can also be of independent
interest outside the nonparametric regression problems.\\
\\
\end{abstract}

\section{Introduction\label{sec:Introduction}}

Estimation of an unknown univariate function $f$ from the nonparametric
regression model 
\begin{equation}
Y_{i}=f(X_{i})+\varepsilon_{i},\,i=1,...,n\label{eq:model}
\end{equation}
has been a central research topic in econometrics, machine learning,
numerical analysis and statistics. Many semiparametric estimators
involve nonparametric regressions as an intermediate step, and some
of the classical examples in economics can be found in several \textit{Handbook
of Econometrics} chapters such as Powell (1994), Chen (2007), and
Ichimura and Todd (2007). The typical assumption about $f$ in (\ref{eq:model})
is that it has derivatives up to a given $(\gamma+1)$th order bounded
in absolute values by a common constant everywhere or almost everywhere
(a.e.). Given an estimator of $f$, an important object of interest
is the convergence rate of the mean integrated squared error (MISE)
of this estimator and the minimax optimality property of the MISE
rate, which tells one how fast the population mean squared distance
between the estimator and $f$ shrinks to zero uniformly when $f$
ranges over a smoothness class, as the sample size $n$ increases.
In particular, MISE is a global mean squared error criterion by integrating
over all possible values of the covariate and noise with respect to
some underlying distribution (see Pagan and Ullah, 1999).

When the smoothness degree $\gamma+1$ is finite and the sample size
$n\rightarrow\infty$, existing \textit{asymptotic} results show that
the minimax optimal MISE rate is $\left(\frac{1}{n}\right)^{\frac{2\gamma+2}{2\gamma+3}}$
which decreases as $\gamma+1$ increases. In the literature of statistical
learning theory, $\left(\frac{1}{n}\right)^{\frac{2\gamma+2}{2\gamma+3}}$
is also stated as the \textit{nonasymptotic} minimax optimal rate
and derived without assuming $n\rightarrow\infty$. See Tsybakov (2009)
for a comprehensive review of the literature that show the classical
asymptotic minimax optimal rate $\left(\frac{1}{n}\right)^{\frac{2\gamma+2}{2\gamma+3}}$
under the assumption of $\gamma+1$ being finite and $n$ tending
to infinity. See Wainwright (2019) for the nonasymptotic derivation.
The rate $\left(\frac{1}{n}\right)^{\frac{2\gamma+2}{2\gamma+3}}$
gives arise to the so called ``blessing of smoothness'' folklore
(i.e., the more smoothness one can exploit, the better) in the theoretical
community. 

In terms of the minimax optimal MISE rates associated with the \textit{standard}
$(\gamma+1)$th degree smoothness classes\footnote{In this paper, a \textit{standard} $(\gamma+1)$th degree smoothness
class refers to one that consists of functions with derivatives belonging
to a ball of some constant radius independent of the derivative order
and bounded away from zero and from above, with respect to either
an $l_{\infty}$ (max) norm or a Hilbert norm.}, we show the following results: (i) if $n\leq\left(\gamma+1\right)^{2\gamma+3}$
(the ``small $n$'' regime), the minimax optimal MISE rate is $\frac{\log n}{n\log(\log n)}$
and the optimal degree of smoothness to exploit is roughly $\max\left\{ \left\lfloor \frac{\log n}{2\log\left(\log n\right)}\right\rfloor ,\,1\right\} $;\footnote{``$\log$'' in this paper is natural logarithm. The more precise
characterization of the optimal degree of smoothness to exploit under
$n\leq\left(\gamma+1\right)^{2\gamma+3}$ is detailed in Section \ref{sec:Minimax-standard}.} (ii) if $n>\left(\gamma+1\right)^{2\gamma+3}$ (the ``large $n$''
regime, which clearly includes the case of $\gamma$ being finite
and $n$ tending to infinity), the minimax optimal MISE rate is $\left(\frac{1}{n}\right)^{\frac{2\gamma+2}{2\gamma+3}}$
and the optimal degree of smoothness to exploit is $\gamma+1$. To
our knowledge, this paper is the first in the literature to show the
minimax optimal rate in the ``small $n$'' regime and the sample
size (i.e., $\left(\gamma+1\right)^{2\gamma+3}$) at which the minimax
optimal rate transitions from $\frac{\log n}{n\log(\log n)}$ to $\left(\frac{1}{n}\right)^{\frac{2\gamma+2}{2\gamma+3}}$. 

The definition of our minimax optimal rates is based on the existing
literature. A rate is said to be minimax optimal in our problem if
we can show: (1) the MISE for \textit{any} estimators (by taking the
infimum over all estimators) in the worst case scenario (by taking
the supremum over a $(\gamma+1)$th degree smoothness class) is bounded
from below, and such a bound is called a minimax lower bound; (2)
there exists an estimator such that its MISE in the worst case scenario
has an upper bound that matches the lower bound in terms of the rate,
the so-called achievability result; to be precise, apart from some
universal constant independent of $n$ (and in our interest, also
independent of $\gamma$), the upper bound matches the lower bound
and the matching part is the minimax optimal rate. The difference
in the ``constants'' between our results and the existing literature
is that the constants in our bounds are truly \textit{universal constants
that are independent of $n$ and $\gamma$ and bounded away from zero
and from above}, while those in the existing literature cannot be
independent of $\gamma$. This difference allows us to provide an
asymptotic interpretation for the rate $\frac{\log n}{n\log(\log n)}$
in the ``small $n$'' regime (in Section \ref{sec:Minimax-standard})
relative to the classical rate $\left(\frac{1}{n}\right)^{\frac{2\gamma+2}{2\gamma+3}}$:
\begin{itemize}
\item as $n\rightarrow\infty$ and $\left(\gamma+1\right)^{2\gamma+3}\rightarrow\infty$,
while 
\begin{equation}
n=o\left(\left(\gamma+1\right)^{2\gamma+3}\right),\label{eq:asymptotic condition}
\end{equation}
then 
\begin{equation}
\left(\frac{1}{n}\right)^{\frac{2\gamma+2}{2\gamma+3}}=o\left(\frac{\log n}{n\log(\log n)}\right);\label{eq:asymptotic interpretation}
\end{equation}
\item as $n\rightarrow\infty$ and $\left(\gamma+1\right)^{2\gamma+3}\rightarrow\infty$,
while 
\begin{equation}
\left(\gamma+1\right)^{2\gamma+3}=o(n),\label{eq:asymptotic condition_1}
\end{equation}
then 
\[
\frac{\log n}{n\log(\log n)}=o\left(\left(\frac{1}{n}\right)^{\frac{2\gamma+2}{2\gamma+3}}\right);
\]
that is, the rate $\left(\frac{1}{n}\right)^{\frac{2\gamma+2}{2\gamma+3}}$
kicks in.\footnote{It is worth pointing out that the results in Section \ref{sec:Minimax-standard}
hold both non-asymptotically and asymptotically. One may interpret
the results in Section \ref{sec:Minimax-standard} in the context
of a triangular data generating process where the smoothness degree
is indexed by the sample size. } 
\end{itemize}
Particularly, we show in this paper that, if the maximum smoothness
degree of $f$ is $\gamma+1$, estimators which minimize the sum of
squared residuals and are constrained to exploit the optimal degree
of smoothness achieve the abovementioned minimax optimal MISE rates.
These estimators will be referred to as the constrained nonparametric
least squares estimator (CNLS) in the following. CNLS estimators constrained
to be in a Sobolev class associated with a Reproducing Kernel Hilbert
Space (RKHS) radius have nice closed form expressions via kernel functions
and are easy to implement in the regularized form, often referred
to as the kernel ridge regression (KRR) estimators (among the most
popular nonparametric estimators) in machine learning. There is a
rich theory based on RKHS for the asymptotic properties of such estimators
under the regime where $\gamma+1$ is finite and $n\rightarrow\infty$
(see, e.g., Schölkopf and Smola, 2002; Berlinet and Thomas-Agnan,
2004). These estimators are closely related to smoothing splines methods
and Gaussian process regressions (see, Wahba, 1990; Rasmussen and
Williams, 2006).

It is worth mentioning the connection between our theoretical results
and some numerical findings from the literature. From a practical
view-point, Marron and Wand (1992) find the exact MISE of kernel density
estimators based on Gaussian kernels can increase with the order of
kernel being exploited (the assumed degree of smoothness) when the
sample size is moderate. Marron (1994) further shows in simulations
that the second order kernel produces a smaller MISE than the fourth
order kernel when the sample size is between $70$ and $10000$, and
the fourth order kernel is dominantly better than the second order
kernel when $n>10000$. Despite that our focus in this paper is on
the minimax optimal rates and our achievability results concern global
nonparametric procedures such as KRR (instead of kernel density estimators
based on Gaussian kernels), Marron and Wand (1992) and our paper share
one general message: the classical rate $\left(\frac{1}{n}\right)^{\frac{2\gamma+2}{2\gamma+3}}$
is an underestimate of the MISE when $n$ is not large enough. Another
message conveyed by our paper is that the optimal degree of smoothness
to exploit should depend on the sample size. 

In econometrics, there is increasing empirical evidence against exploiting
higher degrees of smoothness. Graham, et. al (2010) comment, ``As
is usual in semiparametric estimation, higher order kernels are required
for bias reduction, although the use of such kernels in practice may
be ill-advised.'' Based on several empirical studies with sample
sizes ranging from a couple of thousands to at most thirties of thousands,
Gelman and Imbens (2019) recommend researchers to avoid using high
order polynomials but use local linear or local quadratic polynomials
to estimate the two conditional mean functions of a pretreatment variable
in regression discontinuity designs (RDD) analyses. 

The implication of our results extends to other applications. Establishing
the convergence rate of a semiparametric procedure often requires
establishing an MISE rate or its sample analogue concerning a first-step
nonparametric regression. Below are a couple of examples:
\begin{itemize}
\item When applying a 2SLS-type procedure to estimate a triangular system
where the first-stage equations linking the endogenous regressors
with instruments take the form of (\ref{eq:model}), the MSE of the
2SLS estimator for the parameters of interest in the second-stage
(main) equation depends on the MISE of the first-stage estimators.
\item When applying the partialling-out type strategy to estimate the parameters
of interest in a partially linear model, the first step uses a nonparametric
regression to obtain the partial residuals, and the second step uses
a least squares procedure or a regularized least squares procedure
based on the estimated residuals from the first step. The MSE of the
second step estimator depends on the MISE of the first-step estimator.\footnote{For the partially linear models in Item 2, following derivations similar
to those in Zhu (2017) would reveal the dependence. For the triangular
systems in Item 1, following derivations similar to those in Zhu (2018)
would reveal the dependence; in particular, the modifications involve
replacing terms like $\left(Z_{ij}\hat{\pi}_{j}-Z_{ij}\pi_{j}^{*}\right)^{2}$
in Zhu (2018) with $\left(\hat{f}\left(Z_{ij}\right)-f\left(Z_{ij}\right)\right)^{2}$,
where $j$ is the index of the endogenous variables and $Z_{ij}$
is an instrumental variable for the $j$th endogenous variable.}
\end{itemize}

\section{Overview of our results}

\subsection{Preliminaries }

\textbf{Notation}. Let $\left\lfloor x\right\rfloor $ denote the
largest integer smaller than or equal to $x$. For two functions $f(n,\gamma)$
and $g(n,\gamma)$, let us write $f(n,\gamma)\succsim g(n,\gamma)$
if $f(n,\gamma)\geq cg(n,\gamma)$ for a universal constant $c\in(0,\,\infty)$;
similarly, we write $f(n,\gamma)\precsim g(n,\gamma)$ if $f(n,\gamma)\leq c^{'}g(n,\gamma)$
for a universal constant $c^{'}\in(0,\,\infty)$; and $f(n,\gamma)\asymp g(n,\gamma)$
if $f(n,\gamma)\succsim g(n,\gamma)$ and $f(n,\gamma)\precsim g(n,\gamma)$.
Throughout this paper, we use various $c$ and $C$ letters to denote
positive universal constants that are: finite and bounded away from
zero (denoted by $\asymp1$) and independent of $n$ and $\gamma$
and the dimension $d$ of the covariates (when $d-$dimensional covariates
are of interest); these constants may vary from place to place. 

For a $J-$dimensional vector $\theta$, the $l_{q}-$norm $\left|\theta\right|{}_{q}:=\left(\sum_{j=1}^{J}|\theta_{j}|^{q}\right)^{1/q}$
if $1\leq q<\infty$ and $\left|\theta\right|{}_{q}:=\max_{j\in\left\{ 1,...,J\right\} }\left|\theta_{j}\right|$
if $q=\infty$. Let $\mathbb{B}_{q}^{J}\left(R\right):=\left\{ \theta\in\mathbb{R}^{J}:\left|\theta\right|_{q}\leq R\right\} $.
For functions on $\left[a,\,b\right]$, the unweighted $L^{2}-$norm
$\left|f-g\right|_{2}:=\sqrt{\int_{a}^{b}\left[f\left(x\right)-g\left(x\right)\right]^{2}dx}$,
and the weighted $L^{2}\left(\mathbb{P}\right)-$norm $\left|f-g\right|_{2,\mathbb{P}}:=\sqrt{\int_{a}^{b}\left[f\left(x\right)-g\left(x\right)\right]^{2}\mathbb{P}(dx)}$. 

For functions on $\left[a,\,b\right]^{d}$, the supremum norm $\left|f-g\right|_{\infty}:=\sup_{x\in\left[a,\,b\right]^{d}}\left|f\left(x\right)-g\left(x\right)\right|$. 

Finally, the $\mathcal{L}^{2}(\mathbb{P}_{n})-$norm of the vector
$f:=\left\{ f(x_{i})\right\} _{i=1}^{n}$, denoted by $\left|f\right|_{n}$,
is $\left[\frac{1}{n}\sum_{i=1}^{n}\left(f(x_{i})\right)^{2}\right]^{\frac{1}{2}}$.

\subsubsection{Classes of smooth functions\label{subsec:Classes-of-smooth}}

\begin{definition}[Generalized Hölder subclass] For a non-negative
integer $\gamma$, we let the Hölder\textit{ class} $\mathcal{U}_{\gamma+1}\left(\left(R_{k}\right)_{k=0}^{\gamma+1},\,\left[-1,\,1\right]\right)$
be the class of functions such that any function $f\in\mathcal{U}_{\gamma+1}\left(\left(R_{k}\right)_{k=0}^{\gamma+1},\,\left[-1,\,1\right]\right)$
satisfies: (1) $f$ is continuous on $\left[-1,\,1\right]$ and all
derivatives of $f$ exist; (2) $\left|f^{k}\left(x\right)\right|\leq R_{k}$
for all $k=0,...,\gamma$ and $x\in\left[-1,\,1\right]$, where $f^{0}\left(x\right)=f\left(x\right)$;
(3) $\left|f^{\gamma}(x)-f^{\gamma}(x^{'})\right|\leq R_{\gamma+1}\left|x-x^{'}\right|$
for all $x,\,x^{'}\in\left[-1,\,1\right]$. 

\end{definition}

\begin{remark} Note that in our definition of $\mathcal{U}_{\gamma+1}$,
the absolute values of the derivatives of any member are allowed to
depend on the derivative order. This is motivated by nonparametric
procedures for recovering solutions of ordinary differential equations
in the literature of functional data analysis. We refer interested
readers to Zhu and Mirzaei (2021) for the details.

\end{remark}

Any function $f$ in a Hölder class on $\left[-1,\,1\right]$ can
be written as

\begin{equation}
f(x)=\underset{poly}{\underbrace{\sum_{k=0}^{\gamma}\frac{x^{k}}{k!}f^{(k)}(0)}}+\frac{x^{\gamma}}{\gamma!}f^{(\gamma)}(z)-\frac{x^{\gamma}}{\gamma!}f^{(\gamma)}(0)\label{eq:expansion}
\end{equation}
where $z$ is some intermediate value between $x$ and $0$, and we
follow the convention $0^{0}=0!=1$. 

\begin{definition}[Generalized polynomial subclass] The \textit{generalized
polynomial subclass consisting of functions in the form of ``$poly$''
in (\ref{eq:expansion}) can be expressed as}
\[
\mathcal{U}_{\gamma+1,1}=\left\{ f(x)=\sum_{k=0}^{\gamma}\theta_{k}x^{k}:\textrm{\ensuremath{\left(\theta_{k}\right)_{k=0}^{\gamma}}\ensuremath{\ensuremath{\in}}\ensuremath{\mathcal{P}_{\gamma}}},\,x\in\left[-1,\,1\right]\right\} 
\]
with the $\left(\gamma+1\right)-$dimensional polyhedron 
\[
\ensuremath{\mathcal{P}_{\gamma}=\left\{ \ensuremath{\left(\theta_{k}\right)_{k=0}^{\gamma}}\ensuremath{\ensuremath{\in}}\mathbb{R}^{\gamma+1}:\theta_{k}\in\left[\frac{-R_{k}}{k!},\,\frac{R_{k}}{k!}\right]\right\} }.
\]

\end{definition} 

\begin{definition}[Generalized Hölder subclass] The \textit{generalized
Hölder subclass} $\mathcal{U}_{\gamma+1,2}$ is the class of functions
such that any function $f\in\mathcal{U}_{\gamma+1,2}$ satisfies:
$f\in\mathcal{U}_{\gamma+1}$ such that $f^{(k)}\left(0\right)=0$
for all $k=0,...,\gamma$. 

\end{definition} 

Consequently, we have the following relationships: 
\begin{eqnarray}
 & \mathcal{U}_{\gamma+1,1}\subseteq\mathcal{U}_{\gamma+1}, & \mathcal{U}_{\gamma+1,2}\subseteq\mathcal{U}_{\gamma+1},\label{eq:inner}\\
\mathcal{U}_{\gamma+1}\subseteq & \mathcal{U}_{\gamma+1,1}+\mathcal{U}_{\gamma+1,2} & :=\left\{ f_{1}+f_{2}:f_{1}\in\mathcal{U}_{\gamma+1,1},\,f_{2}\in\mathcal{U}_{\gamma+1,2}\right\} .\label{eq:outter}
\end{eqnarray}

As discussed in Wainwright (2019, Chapter 12), any function $f$ in
the Sobolev space on $\left[0,\,1\right]$ has the expansion 
\[
f\left(x\right)=\sum_{k=0}^{\gamma}f^{(k)}(0)\frac{x^{k}}{k!}+\int_{0}^{1}f^{(\gamma+1)}(t)\frac{\left(x-t\right)_{+}^{\gamma}}{\gamma!}dt\quad\textrm{with }(a)_{+}=a\vee0,
\]
and one (RK)HS norm associated with a Sobolev space takes the form
\begin{equation}
\left|f\right|_{\mathcal{H},\gamma+1}=\sqrt{\sum_{k=0}^{\gamma}\left(f^{(k)}(0)\right)^{2}+\int_{0}^{1}\left[f^{(\gamma+1)}\left(t\right)\right]^{2}dt}.\label{eq:sobolev norm}
\end{equation}
Therefore, the Sobolev space can be decomposed into a polynomial subspace
and a Sobolev subspace imposed with the restrictions that $f^{(k)}(0)=0$
for all $k=0,...,\gamma$ and $f^{(\gamma+1)}$ belongs to the space
$\mathcal{L}^{2}\left[0,\,1\right]$ (see Wahba, 1990, Chapter 1;
Wainwright, 2019, Chapter 12, Examples 12.17 and 12.29). A Sobolev
subclass is a special case of the ellipsoid subclass; in particular,
a Sobolev subclass can be expressed in the following way (Wainwright
2019, Chapters 5, 12 and 15). 

\begin{definition}[Generalized ellipsoid subclass] The \textit{generalized
ellipsoid subclass} of smooth functions 

\begin{equation}
\mathcal{H}_{\gamma+1}=\left\{ f=\sum_{m=1}^{\infty}\theta_{m}\phi_{m}:\textrm{for }\left(\theta_{m}\right)_{m=1}^{\infty}\ensuremath{\in}\ell^{2}\left(\mathbb{N}\right)\textrm{ such that }\sum_{m=1}^{\infty}\frac{\theta_{m}^{2}}{\mu_{m}}\leq R_{\gamma+1}^{2}\right\} \label{eq:Ellipsoid-1}
\end{equation}
where $\ell^{2}\left(\mathbb{N}\right):=\left\{ \left(\theta_{m}\right)_{m=1}^{\infty}\vert\sum_{m=1}^{\infty}\theta_{m}^{2}<\infty\right\} $,
$\left(\mu_{m}\right)_{m=1}^{\infty}$ and $\left(\phi_{m}\right)_{m=1}^{\infty}$
are the eigenvalues and eigenfunctions (that forms an orthonormal
basis of $\mathcal{L}^{2}\left[0,\,1\right]$), respectively, of an
RKHS associated with a continuous and semidefinite kernel function.
\end{definition} 

When considering $\mathcal{H}_{\gamma+1}$ in (\ref{eq:Ellipsoid-1}),
we will assume that $\mu_{m}=\left(cm\right)^{-2\left(\gamma+1\right)}$
for a positive constant $c\asymp1$ independent of $\gamma$ and $R_{\gamma+1}$.
The decay rate of the eigenvalues follows the standard assumption
for $(\gamma+1)-$degree smooth functions in the literature (see,
e.g., Steinwart and Christmann, 2008; Wainwright, 2019) and $R_{\gamma+1}\asymp1$
in (\ref{eq:Ellipsoid-1}) gives the \textit{standard} ellipsoid subclass.
Moreover, (\ref{eq:Ellipsoid-1}) is equipped with the inner product
$\left\langle h,\,g\right\rangle _{\mathcal{H}}=\sum_{m=1}^{\infty}\frac{\left\langle h,\,\phi_{m}\right\rangle \left\langle g,\,\phi_{m}\right\rangle }{\mu_{m}}$
where $\left\langle \cdot,\cdot\right\rangle $ is the inner product
in $\mathcal{L}^{2}\left[0,\,1\right]$. 

\subsubsection{Metric entropy}

\begin{definition}[Covering and packing numbers]

Given a set $\Lambda$, a set $\left\{ \mathrm{\eta}^{1},\,\mathrm{\eta}^{2},...,\mathrm{\eta}^{N}\right\} \subset\Lambda$
is a $\delta-$cover of $\Lambda$ in the metric $\rho$ if for each
$\eta\in\Lambda$, there exists some $i\in\left\{ 1,...,N\right\} $
such that $\rho(\eta,\,\eta^{i})\leq\delta$. The $\delta-$covering
number of $\Lambda$, denoted by $N_{\rho}(\delta,\,\Lambda)$, is
the cardinality of the smallest $\delta-$cover. A set $\left\{ \mathrm{\eta}^{1},\,\mathrm{\eta}^{2},...,\mathrm{\eta}^{M}\right\} \subset\Lambda$
is a $\delta-$packing of $\Lambda$ in the metric $\rho$ if for
any distinct $i,j\in\left\{ 1,...,M\right\} $, $\rho(\eta^{i},\,\eta^{j})>\delta$.
The $\delta-$packing number of $\Lambda$, denoted by $M_{\rho}(\delta,\,\Lambda)$,
is the cardinality of the largest $\delta-$packing. Throughout this
paper, we use $N_{q}\left(\delta,\,\mathcal{\mathcal{F}}\right)$
and $M_{q}\left(\delta,\,\mathcal{\mathcal{F}}\right)$ to denote
the $\delta-$covering number and the $\delta-$packing number, respectively,
of a function class $\mathcal{F}$ with respect to the function norm
$\left|\cdot\right|_{q}$ where $q\in\left\{ 2,\,\infty\right\} $;
moreover, $N_{2,\mathbb{P}}\left(\delta,\,\mathcal{\mathcal{F}}\right)$
and $M_{2,\mathbb{P}}\left(\delta,\,\mathcal{\mathcal{F}}\right)$
denote the $\delta-$covering number and the $\delta-$packing number,
respectively, of a function class $\mathcal{F}$ with respect to the
weighted $L^{2}\left(\mathbb{P}\right)-$norm $\left|\cdot\right|_{2,\mathbb{P}}$.

\end{definition}

The following is a standard textbook result that summarizes the relationships
between covering and packing numbers: 
\begin{equation}
M_{\rho}(2\delta,\,\Lambda)\leq N_{\rho}(\delta,\,\Lambda)\leq M_{\rho}(\delta,\,\Lambda).\label{eq:sandwich}
\end{equation}
Given this sandwich result, a lower bound on the packing number gives
a lower bound on the covering number, and vice versa; similarly, an
upper bound on the covering number gives an upper bound on the packing
number, and vice versa.

Metric entropy is an important concept in approximation theory and
discrete geometry. In mathematical statistics and machine learning
theory, metric entropy is a fundamental building block. Combined with
the Fano's inequality from information theory (see, Cover and Thomas,
2005), it allows one to derive the minimax lower bounds for the MISE
rates; combined with the notion of ``local complexity'' in empirical
processes theory, it allows one to derive upper bounds on the MISE
rates.\footnote{For the \textit{standard} smoothness classes, one may use methods
based on ``ranks'' (for the polynomial subclass) and ``eigenvalues''
(for the nonparametric subclass) to derive upper bounds on the MISE.
However, ``ranks'' are not very useful for deriving the minimax
lower bounds in general. Even for upper bounds in the case of a Hölder
class, once we allow $R_{k}$ to depend on the order of derivative,
$k$, the ``rank''--based argument is hard to generalize as it
does not account for the impact of $R_{k}$.}

\subsection{Our contributions\label{subsec:Our-contributions}}

In the existing literature, the minimax lower bound and the upper
bound (for achievability) for MISE is typically stated as $\mathtt{constant}\cdot\left(\frac{1}{n}\right)^{\frac{2\gamma+2}{2\gamma+3}}$
and $\mathtt{constant'}\cdot\left(\frac{1}{n}\right)^{\frac{2\gamma+2}{2\gamma+3}}$,
respectively, where the specific forms of $\mathtt{constant}$ and
$\mathtt{constant}'$ are not derived. Then the matching part ``$\left(\frac{1}{n}\right)^{\frac{2\gamma+2}{2\gamma+3}}$''
is claimed as the minimax optimal rate (e.g., Wainwright, 2019). Deriving
constants for global criteria such as MISE and demonstrating achievability
(upper bounds) in the context of global nonparametric procedures (the
focus of our paper) is known to be difficult, which is why the existing
literature does not make an attempt in deriving the constants. 

Relative to the existing literature, our results take one step further
by revealing more explicit dependence on $n$, $\gamma$ and $\left\{ R_{k}\right\} _{k=0}^{\gamma+1}$.
The difference in the ``constants'' between our results and the
existing literature is that the constants in our bounds are truly
\textit{universal} constants that are independent of $n$ and $\gamma$
and bounded away from zero and from above, while those in the existing
literature cannot be independent of $\gamma$, as explained in the
next section. This difference allows us to provide an asymptotic interpretation
for the rate $\frac{\log n}{n\log(\log n)}$ in the ``small $n$''
regime (in Section \ref{sec:Minimax-standard}) relative to the classical
rate $\left(\frac{1}{n}\right)^{\frac{2\gamma+2}{2\gamma+3}}$, as
discussed in Section \ref{sec:Introduction}.\footnote{Because of the complexity of our problems, we make no attempt to derive
the explicit universal constants that are \textit{independent} of
$n$ and $\gamma$. The following process is a standard one used in
this literature for proving minimax optimal rates in MISE and explains
why such a derivation is not carried out: First, results on metric
entropy are derived. Then, these results are applied with other technical
lemmas to establish minimax optimal rates for the MISE. Throughout,
the universal constants become complex rather quickly and it becomes
physically impossible to track them from lines to lines. Based on
the existing mathematical tools from approximation theory, information
theory, minimax lower bounds theory, and empirical processes theory
(for upper bounds), deriving meaningful universal constants (letting
alone sharp universal constants) in the particular context of this
paper is infeasible and remains an open problem until a new mathematical
foundation can be introduced.} 

The fundamental contribution of this paper lies in a set of metric
entropy bounds. Some of them are original, and some of them improve
and/or extend the ones in the literature to allow general (possibly
$k-$dependent) $R_{k}$s. Besides their applications in the empirical
processes theory and the theory of minimax lower bounds, our bounds
for the covering and packing numbers are of independent interest for
researchers working in approximation theory and discrete geometry.

\subsection{Related literature and the novelty of our results}

Here we start with a short discussion of why the phase transition
phenomenon has been overlooked in the literature while leaving the
details to the subsequent sections. For $\mathcal{U}_{\gamma+1}$
under the assumption that $R_{k}=\overline{C}\asymp1$ (which we will
refer to as the \textit{standard} $\mathcal{U}_{\gamma+1}$), Kolmogorov
and Tikhomirov (1959) show that 
\begin{eqnarray}
\log(\delta-\textrm{covering number}) & \precsim & \left(\gamma+1\right)\log\frac{1}{\delta}+\delta^{\frac{-1}{\gamma+1}},\label{eq:Kolmogorov_covering}\\
\log(\delta-\textrm{packing number}) & \succsim & \delta^{\frac{-1}{\gamma+1}},\label{eq:Kolmogorov_packing}
\end{eqnarray}
for $\delta-$approximation accuracy. By (\ref{eq:sandwich}), (\ref{eq:Kolmogorov_packing})
implies 
\begin{equation}
\log(\delta-\textrm{covering number})\succsim\delta^{\frac{-1}{\gamma+1}}.\label{eq:Kolmogorov_cover_lower}
\end{equation}

The minimax optimal MISE rates (as functions of $n$ and $\gamma$)
are related to choices of $\delta$ (see, e.g., Yang and Barron, 1999;
Wainwright, 2019). Virtually every paper including recent textbooks
on nonasymptotic statistics (such as Wainwright, 2019, eq. 5.17 on
p. 129) takes $\log(\delta-\textrm{covering/packing number})\asymp\delta^{\frac{-1}{\gamma+1}}$.
It would take some diligence for one to recognize that: 
\begin{itemize}
\item when $\gamma+1$ and $\delta$ are large enough, taking $\log(\delta-\textrm{covering/packing number})\asymp\delta^{\frac{-1}{\gamma+1}}$
is problematic even for the \textit{standard} $\mathcal{U}_{\gamma+1}$,
because the term $\left(\gamma+1\right)\log\frac{1}{\delta}$ would
dominate $\delta^{\frac{-1}{\gamma+1}}$ in (\ref{eq:Kolmogorov_covering}),
and therefore, there is a significant gap between (\ref{eq:Kolmogorov_covering})
and (\ref{eq:Kolmogorov_cover_lower}); 
\item when $\gamma+1$ and $\delta$ are small enough, taking $\log(\delta-\textrm{covering/packing number})\asymp\delta^{\frac{-1}{\gamma+1}}$
is good enough for the \textit{standard} $\mathcal{U}_{\gamma+1}$,
because the term $\left(\gamma+1\right)\log\frac{1}{\delta}$ would
be dominated by $\delta^{\frac{-1}{\gamma+1}}$ in (\ref{eq:Kolmogorov_covering}),
and (\ref{eq:Kolmogorov_cover_lower}) is sharp for the \textit{standard}
$\mathcal{U}_{\gamma+1}$. 
\end{itemize}
Here is a list of our discoveries:
\begin{itemize}
\item the derivation of the lower bound $\delta^{\frac{-1}{\gamma+1}}$
under $R_{k}=\overline{C}\asymp1$ in (\ref{eq:Kolmogorov_packing})
as well as the following literature for \textit{Hölder and Sobolev
classes} ignore the \textit{polynomial subclass} $\mathcal{U}_{\gamma+1,1}$;
see e.g., Wainwright (2019), Example 5.11, which inherits the lower
bound derivation in Kolmogorov and Tikhomirov (1959). The lower bound
$\delta^{\frac{-1}{\gamma+1}}$ is a lower bound for the \textit{standard
Hölder or Sobolev/ellipsoid subclass} only, and is not sharp even
for the \textit{standard} smoothness classes if $\gamma+1$ and $\delta$
are large enough; 
\item the upper bound based on the arguments in Kolmogorov and Tikhomirov
(1959) for the \textit{polynomial subclass} $\mathcal{U}_{\gamma+1,1}$
is far from being tight when $R_{k}$s become large enough; 
\item the upper bound based on the arguments in Mityagin (1961) and the
following literature (Wainwright, 2019, Example 5.12) for the \textit{Sobolev/ellipsoid
subclass} does not give the sharp dependence on $\gamma$ and $R_{\gamma+1}$;
\item the minimax optimal rate for MISE or its sample analogue in the existing
nonasymptotic literature is stated as $\left(\frac{1}{n}\right)^{\frac{2\gamma+2}{2\gamma+3}}$,
which is derived based on $\delta^{\frac{-1}{\gamma+1}}$, the metric
entropy for the \textit{standard Hölder or Sobolev subclass only}.
For example, see the last sentence in Example 13.15 of Wainwright
(2019, p.436). Also see Yang and Barron (1999, p.1591, 3rd paragraph).\footnote{For a \textit{non-standard} Hölder class, Zhu and Mirzaei (2021) apply
the argument in Kolmogorov and Tikhomirov (1959) to derive an upper
bound on the covering number. For the lower bound on the covering/packing
number, Zhu and Mirzaei (2021) simply take the classical result $\delta^{\frac{-1}{\gamma+1}}$
from Kolmogorov and Tikhomirov (1959). The consequence is, the lower
and upper error bounds have different rates, neither of which is sharp.
See Theorem \ref{thm:MISE_(k-1)!} in Appendix \ref{sec:MISE-nonstandard}
for the sharp rates.}
\end{itemize}
All these issues are addressed in this paper. Particularly, in deriving
the metric entropy bounds for the polynomial subclass $\mathcal{U}_{\gamma+1,1}$,
we develop our own arguments; in deriving the metric entropy bounds
for the ellipsoid subclass, we base our arguments on the existing
literature but use an improved truncation strategy that gives our
resulting bounds the sharp dependence on $\gamma$ and $R_{\gamma+1}$.
In contrast to those in the existing literature, our entropy bounds
enable us to reveal the phase transition phenomenon and obtain the
sharp MISE rates in the respective ``small $n$'' and ``large $n$''
regimes, associated with the standard smoothness classes as well as
non-standard smoothness classes motivated in Zhu and Mirzaei (2021).

\subsection{Organization }

The results and proofs of the paper are organized as follows.
\begin{itemize}
\item Section \ref{sec:Minimax-standard}: Minimax optimal rates in commonly
seen cases. Appendix \ref{sec:MISE-nonstandard}: Minimax optimal
rates in non-standard cases. The proofs for Section \ref{sec:Minimax-standard}
and Appendix \ref{sec:MISE-nonstandard} are given in Appendix \ref{sec:Proofs-for-MISE}.
\item Section \ref{sec:Covering-and-packing}: Covering and packing numbers.
These results are used to prove the minimax optimal MISE rates in
Section \ref{sec:Minimax-standard} and Appendix \ref{sec:MISE-nonstandard}.
The proofs for Section \ref{sec:Covering-and-packing} are given in
Appendix \ref{sec:proofs_entropy}. 
\item Section \ref{sec:Discussions}: Discussions. Recommendation from the
applied literature (Section \ref{subsec:Practical-implications});
Open questions (Section \ref{subsec:Open-questions}).
\item Appendix \ref{sec:multi-dim-extensions}: Some insights about multivariate
smooth functions. The proofs for Appendix \ref{sec:multi-dim-extensions}
are given in Appendix \ref{sec:Proofs-for-multi-dim-extension}.
\item Appendix \ref{sec:Supporting-lemmas}: Additional supporting lemmas
for Appendix \ref{sec:Proofs-for-MISE}.
\end{itemize}

\section{Minimax optimal rates in commonly seen cases\label{sec:Minimax-standard}}

In this section, we revisit the minimax optimal MISE rates associated
with some commonly seen smoothness classes in the literature. 

\begin{definition}[standard smoothness classes]\label{Standard smoothness classes}
Let $\overline{C}\asymp1$ be a universal constant independent of
$n$ and $\gamma$. The \textit{standard} Hölder class corresponds
to $\mathcal{U}_{\gamma+1}$ with $R_{k}=\overline{C}$ for all $k=0,...,\gamma+1$.
We define the \textit{standard} Sobolev class as follows:
\begin{align}
\mathcal{S}_{\gamma+1}:= & \{f:\,\left[0,\,1\right]\rightarrow\mathbb{R}\vert\,f\textrm{ is \ensuremath{\gamma+1} times differentiable a.e.,}\nonumber \\
 & f^{(\gamma)}\textrm{ is absolutely continuous and}\left|f\right|_{\mathcal{H},\gamma+1}\leq\overline{C}\}\label{eq:sobolev_std}
\end{align}
with $\left|f\right|_{\mathcal{H},\gamma+1}$ defined in (\ref{eq:sobolev norm}).

\end{definition}

When we write $\mathcal{U}_{\gamma+1}$ or $\mathcal{S}_{\gamma+1}$
in this section, $\mathcal{U}_{\gamma+1}$ refers to the \textit{standard}
Hölder class and $\mathcal{S}_{\gamma+1}$ refers to the \textit{standard}
Sobolev class. The regression model (\ref{eq:model}) is subject to
the following conditions.

\begin{assumption}\label{Assumption 1} $\left\{ \varepsilon_{i}\right\} _{i=1}^{n}$
are independent $\mathcal{N}\left(0,\,\sigma^{2}\right)$ where $\sigma\asymp1$;
$\left\{ \varepsilon_{i}\right\} _{i=1}^{n}$ are independent of $\left\{ X_{i}\right\} _{i=1}^{n}$;
$\left\{ X_{i}\right\} _{i=1}^{n}$ are independent draws from a distribution
with density $p(x)$ on a bounded interval associated with our smoothness
classes. 

\end{assumption} 

The proofs for the minimax optimal rates in this section as well as
Appendix \ref{sec:MISE-nonstandard} are based on our results in Section
\ref{sec:Covering-and-packing} as well as techniques from empirical
processes, machine learning theory, and information theory, collected
in Appendix \ref{sec:Supporting-lemmas}. 

\subsubsection*{Allowing for heteroscedasticity and non-Gaussian noise}

Assumption \ref{Assumption 1} is the most standard in the literature
for minimax lower bounds; see, e.g., Yang and Barron (1999), Raskutti,
et. al (2011), and Wainwright (2019). The main reason is, a minimax
lower bound derived under a more restrictive model would be a lower
bound under a less restrictive model that includes the former as a
special case. For example, a homoscedastic model is a special case
of a heteroscedastic one, and the independence of $\varepsilon_{i}$
and $X_{i}$ along with $\mathbb{E}\left(\varepsilon_{i}\right)=0$
implies $\mathbb{E}\left(\varepsilon_{i}\vert X_{i}\right)=0$. Meanwhile,
upper bounds for known estimators under less restrictive model assumptions
are typically available or can be established before developing the
minimax lower bounds. See the introduction of Chapter 15 in Wainwright
(2019) for the insights provided by minimax lower bounds. 

Until very recently, Zhao and Yang (2022) find a way to allow the
noise terms to be heteroscedastic and have non-Gaussian distributions.
Specifically, Zhao and Yang (2022) consider 
\begin{equation}
Y_{i}=f(X_{i})+\underset{\varepsilon_{i}}{\underbrace{\sigma(X_{i})w_{i}}},\,i=1,...,n,\label{eq:heter model}
\end{equation}
where, conditional on $X_{i}$, $Y_{i}$ has a density from a location-scale
family subjective to some rather mild conditions. This setup allows
Zhao and Yang (2022) to derive an upper bound for the Kullback-Leibler
divergence between distributions that are not Gaussian (but can be
Gaussian). The bound in Zhao and Yang (2022) can be used to relax
homoscedasticity without changing the main arguments of our proofs.
For more details related to proofs under (\ref{eq:heter model}),
we refer the interested readers to Appendix \ref{subsec:proof-allowing-for-heter-nonGaussian}
and the remark following Lemma C.1 in Appendix \ref{sec:Supporting-lemmas}
of our paper. It is worth pointing out that, when $\sigma(\cdot)$
in (\ref{eq:heter model}) ranges over the standard smoothness classes
with positive members bounded away from zero and bounded from above,
the minimax lower bounds have the same rate as those under Assumption
\ref{Assumption 1}. 

As for the upper bounds (that demonstrate achievability), our proofs
are developed under a weaker assumption: conditional on $\left\{ X_{i}\right\} _{i=1}^{n}$,
$\left\{ \varepsilon_{i}\right\} _{i=1}^{n}$ are independent sub-Gaussian
random variables with zero mean and sub-Gaussian parameters at most
$\overline{\sigma}\asymp1$. This assumption is very standard in the
literature of empirical process theory (e.g., van der Vaart and Wellner,
1996; van de Geer, 2000; Zhu, 2017; Wainwright, 2019). However, to
establish minimax optimality, it is more sensible to consider a model
with the same set of assumptions on $\left\{ \varepsilon_{i}\right\} _{i=1}^{n}$
for the minimax lower bounds and upper bounds. Therefore, to focus
on the key points, we will state the results under Assumption \ref{Assumption 1}
in the main paper and Appendix \ref{sec:MISE-nonstandard} to be consistent
with the way the well-known results are presented in the literature
for minimax optimality. 

\subsection{Mean integrated squared error rates\label{subsec:Mean-integrated-squared}}

\begin{theorem}[lower bounds] \textit{\label{thm:MISE-lower-standard}Suppose
Assumption \ref{Assumption 1} holds with density $p(x)$ bounded
away from zero; that is, $p(x)\geq c>0$ for some universal constant
$c$.}

\textit{(i) If 
\begin{equation}
\frac{n}{\sigma^{2}}>\left(\gamma+1\right)^{2\gamma+3},\label{eq:n-lower-Theorem3.1}
\end{equation}
then we have 
\begin{eqnarray*}
\inf_{\tilde{f}}\sup_{f\in\mathcal{S}_{\gamma+1}}\mathbb{E}\left(\left|\tilde{f}-f\right|_{2,\mathbb{P}}^{2}\right) & \geq & \underline{c}_{0}\left(\frac{\sigma^{2}}{n}\right)^{\frac{2\left(\gamma+1\right)}{2\left(\gamma+1\right)+1}},\\
\inf_{\tilde{f}}\sup_{f\in\mathcal{U}_{\gamma+1}}\mathbb{E}\left(\left|\tilde{f}-f\right|_{2,\mathbb{P}}^{2}\right) & \geq & \underline{c}_{0}\left(\frac{\sigma^{2}}{n}\right)^{\frac{2\left(\gamma+1\right)}{2\left(\gamma+1\right)+1}},
\end{eqnarray*}
for some universal constant $\underline{c}_{0}\in(0,\,1]$ independent
of $n$ and $\gamma$, and bounded away from zero. Note that $\frac{\sigma^{2}\left(\gamma+1\right)}{n}<\left(\frac{\sigma^{2}}{n}\right)^{\frac{2\left(\gamma+1\right)}{2\left(\gamma+1\right)+1}}$
in this case.}

\textit{(ii) If 
\begin{equation}
\frac{n}{\sigma^{2}}\leq\left(\gamma+1\right)^{2\gamma+3},\label{eq:small n regime}
\end{equation}
we let $\gamma^{*}\in\left\{ 0,...,\gamma\right\} $ be the smallest
integer such that $\frac{n}{\sigma^{2}}\leq\left(\gamma^{*}+1\right)^{2\gamma^{*}+3}$.
Then we have 
\begin{eqnarray*}
\inf_{\tilde{f}}\sup_{f\in\mathcal{S}_{\gamma+1}}\mathbb{E}\left(\left|\tilde{f}-f\right|_{2,\mathbb{P}}^{2}\right) & \geq\underline{c} & \frac{\sigma^{2}\left(\gamma^{*}+1\right)}{n},\\
\inf_{\tilde{f}}\sup_{f\in\mathcal{U}_{\gamma+1}}\mathbb{E}\left(\left|\tilde{f}-f\right|_{2,\mathbb{P}}^{2}\right) & \geq\underline{c} & \frac{\sigma^{2}\left(\gamma^{*}+1\right)}{n},
\end{eqnarray*}
for some universal constant $\underline{c}\in(0,\,1]$ independent
of $n$ and $\gamma$ such that $\underline{c}\leq\frac{\underline{c}_{0}}{2}$,
and bounded away from zero. Note that $\frac{\sigma^{2}\left(\gamma^{*}+1\right)}{n}\geq\left(\frac{\sigma^{2}}{n}\right)^{\frac{2\left(\gamma+1\right)}{2\left(\gamma+1\right)+1}}$
in this case.}

\end{theorem}

The proof for Theorem \ref{thm:MISE-lower-standard} is given in Section
\ref{subsec:Proof-for-MISE-lower-standard}, which relies on the constructions
behind the bound $\underline{B}_{2}$ in Lemma \ref{lm:entropy-generalized-polynomial},
Lemma \ref{lm:entropy-generalized-H=0000F6lder} and Lemma \ref{lm:entropy-ellipsoid-subclass}
in Section \ref{sec:Covering-and-packing}. 

\subsubsection*{The intuition behind}

Recall from Section \ref{subsec:Classes-of-smooth} that a smoothness
class contains two subclasses: a polynomial subclass and an infinite-dimensional
nonparametric subclass such that any function $f$ in this subclass
has $f^{(k)}\left(0\right)=0$ for all $k=0,...,\gamma$. The proof
for Theorem \ref{thm:MISE-lower-standard} applies a version of the
Fano's inequality, which converts the problem into a multiple classification
problem that tries to distinguish among $M$ members in the function
classes of interest. The set of $M$ members need be sufficiently
large and is naturally connected with a packing set of the function
class.

Based on the constructions behind the bound $\underline{B}_{2}$ in
Lemma \ref{lm:entropy-generalized-polynomial}, we can show, the error
rate that prevents from distinguishing any pair of members in the
packing set of the standard $\gamma$th degree polynomial subclass
(in $\mathcal{U}_{\gamma+1}$ and $\mathcal{S}_{\gamma+1}$, respectively)
is $\frac{\sigma^{2}\left(\gamma+1\right)}{n}$. Based on Lemmas \ref{lm:entropy-generalized-H=0000F6lder}
and \ref{lm:entropy-ellipsoid-subclass}, we can show, the error rate
that prevents from distinguishing any pair of members in the packing
set of the standard $(\gamma+1)$th degree Hölder subclass and Sobolev
subclass is $\left(\frac{\sigma^{2}}{n}\right)^{\frac{2\left(\gamma+1\right)}{2\left(\gamma+1\right)+1}}$.
Note that 
\begin{itemize}
\item when $\frac{n}{\sigma^{2}}>\left(\gamma+1\right)^{2\gamma+3}$, one
has $\frac{\sigma^{2}\left(\gamma+1\right)}{n}<\left(\frac{\sigma^{2}}{n}\right)^{\frac{2\left(\gamma+1\right)}{2\left(\gamma+1\right)+1}}$,
in which case the \textit{rate-minimizing} estimator (the ``min''
part of ``minimax'') would exploit the smoothness degree $\gamma$
(Theorem \ref{thm:MISE-lower-standard}(i));
\item when $\frac{n}{\sigma^{2}}\leq\left(\gamma+1\right)^{2\gamma+3}$,
one has $\frac{\sigma^{2}\left(\gamma+1\right)}{n}\geq\left(\frac{\sigma^{2}}{n}\right)^{\frac{2\left(\gamma+1\right)}{2\left(\gamma+1\right)+1}}$,
in which case the \textit{rate-minimizing} estimator (the ``min''
part of ``minimax'') would exploit the smoothness degree $\gamma^{*}+1$
(as detailed in Theorem \ref{thm:MISE-lower-standard}(ii)) such that
$\frac{\sigma^{2}\left(\gamma^{*}+1\right)}{n}\approx\left(\frac{\sigma^{2}}{n}\right)^{\frac{2\left(\gamma^{*}+1\right)}{2\left(\gamma^{*}+1\right)+1}}$.
In particular, $\frac{\gamma^{*}+1}{n}\asymp\frac{\log n}{n\log(\log n)}$.
\end{itemize}
In deriving the minimax lower bounds, the existing nonasymptotic literature
uses the metric entropy $\delta^{\frac{-1}{\gamma+1}}$ (e.g., Wainwright,
2019, Examples 5.11 and 15.23). The entropy $\delta^{\frac{-1}{\gamma+1}}$
ignores the polynomial subclass. For more details, we refer the interested
readers to Section \ref{sec:Covering-and-packing}. Using $\delta^{\frac{-1}{\gamma+1}}$
to derive the minimax lower bound results in a rate $\left(\frac{\sigma^{2}}{n}\right)^{\frac{2\left(\gamma+1\right)}{2\left(\gamma+1\right)+1}}$
and overlooks the (sharper) rate in Theorem \ref{thm:MISE-lower-standard}(ii),
the ``small $n$'' regime.

The following two theorems show the achievability of the rates in
Theorem \ref{thm:MISE-lower-standard}. In particular, we show that
estimators constrained to exploit the optimal degree of smoothness
in each regime achieve the respective rate in Theorem \ref{thm:MISE-lower-standard}.
We consider the constrained nonparametric least squares estimator
(CNLS)
\begin{equation}
\hat{f}\in\arg\min_{\check{f}\in\mathcal{F}}\frac{1}{2n}\sum_{i=1}^{n}\left(y_{i}-\check{f}\left(x_{i}\right)\right)^{2}.\label{eq:least squares-1}
\end{equation}
The CNLS is commonly seen in the minimax optimality literature (e.g.,
Raskutti, et. al 2011) and the Kernel Ridge Regression (KRR) literature.
Let $k=\gamma$ in the ``large $n$'' regime and $k=\gamma^{*}$
in the ``small $n$'' regime. We consider either $\mathcal{F}=\mathcal{U}_{k+1}$
or 
\begin{align*}
\mathcal{F}= & \mathcal{S}_{k+1}\\
:= & \{f:\,\left[0,\,1\right]\rightarrow\mathbb{R}\vert\,f\textrm{ is \ensuremath{k+1} times differentiable a.e.,}\\
 & f^{(k)}\textrm{ is absolutely continuous and}\left|f\right|_{\mathcal{H},k+1}\leq\overline{C}\}
\end{align*}
where $\left|\cdot\right|_{\mathcal{H},k+1}$ is the norm defined
in (\ref{eq:sobolev norm}). Both cases can be of interest but the
latter is more widely implemented in practice, as we explain below. 

\subsubsection*{Kernel Ridge Regression (KRR) in machine learning}

Constraining the estimators to be in $\mathcal{S}_{k+1}$ allows one
to implement (\ref{eq:least squares-1}) via kernel functions. The
following discussion is based on Wainwright (2019). Let the matrix
$\mathbb{K}_{k+1}\in\mathbb{R}^{n\times n}$ consist of entries $\frac{1}{n}\mathcal{K}_{k+1}\left(x_{i},\,x_{i^{'}}\right)$,
taking the form
\begin{eqnarray}
\mathcal{K}_{k+1}\left(x_{i},\,x_{i^{'}}\right) & = & 1+\left(x_{i}\wedge x_{i^{'}}\right)\quad\textrm{for }k=0,\label{eq:kernel_1}\\
\mathcal{K}_{k+1}\left(x_{i},\,x_{i^{'}}\right) & = & \sum_{j=0}^{k}\frac{x_{i}^{j}}{j!}\frac{x_{i^{'}}^{j}}{j!}+\int_{0}^{1}\frac{\left(x_{i}-t\right)_{+}^{k}}{k!}\frac{\left(x_{i^{'}}-t\right)_{+}^{k}}{k!}dt\quad\textrm{for }k>0,\label{eq:kernel_2}
\end{eqnarray}
where $\left(a\right)_{+}=a\vee0$. The kernel function $\mathcal{K}_{k+1,1}\left(w,\,z\right)=\sum_{j=0}^{k}\frac{w^{j}}{j!}\frac{z^{j}}{j!}$
generates the $k$th degree polynomial subspace and the kernel function
$\mathcal{K}_{k+1,2}\left(w,\,z\right)=\int_{0}^{1}\frac{\left(w-t\right)_{+}^{k}}{k!}\frac{\left(z-t\right)_{+}^{k}}{k!}dt$
for $k>0$ ($w\wedge z$ for $k=0$) generates the $(k+1)$th order
Sobolev subspace imposed with the restrictions that $f^{(j)}(0)=0$
for all $j=0,...,k$ and $f^{(k+1)}$ belongs to the space $\mathcal{L}^{2}\left[0,\,1\right]$.

When $\mathcal{F}=\mathcal{S}_{k+1}$ in (\ref{eq:least squares-1}),
$\hat{f}$ can be written as 
\begin{equation}
\hat{f}\left(\cdot\right)=\frac{1}{\sqrt{n}}\sum_{i^{'}=1}^{n}\hat{\pi}_{i^{'}}\mathcal{K}_{k+1}\left(\cdot,\,x_{i^{'}}\right)\label{eq:KRR estimator}
\end{equation}
where 
\begin{align}
\hat{\pi} & :=\left\{ \hat{\pi}_{i^{'}}\right\} _{i^{'}=1}^{n}=\arg\min_{\pi\in\mathbb{R}^{n}}\frac{1}{2n}\sum_{i=1}^{n}\left(y_{i}-\frac{1}{\sqrt{n}}\sum_{i^{'}=1}^{n}\pi_{i^{'}}\mathcal{K}_{k+1}\left(x_{i},\,x_{i^{'}}\right)\right)^{2}\label{eq:krr}\\
\textrm{s.t.} & \pi^{T}\mathbb{K}_{k+1}\pi\leq\overline{C}^{2}.\label{eq:constraint}
\end{align}
In particular, (\ref{eq:constraint}) comes from the representation
$\left|\check{f}\right|_{\mathcal{H},k+1}^{2}=\pi^{T}\mathbb{K}_{k+1}\pi$
when $\mathcal{F}=\mathcal{S}_{k+1}$ in (\ref{eq:least squares-1})
and $\check{f}$ takes the form $\check{f}\left(\cdot\right)=\frac{1}{\sqrt{n}}\sum_{i^{'}=1}^{n}\pi_{i^{'}}\mathcal{K}_{k+1}\left(\cdot,\,x_{i^{'}}\right)$. 

The program above is convex. Hence, by the Lagrangian duality, solving
(\ref{eq:krr})--(\ref{eq:constraint}) is equivalent to solving
\begin{equation}
\hat{\pi}=\arg\min_{\pi\in\mathbb{R}^{n}}\frac{1}{2n}\sum_{i=1}^{n}\left(y_{i}-\frac{1}{\sqrt{n}}\sum_{i^{'}=1}^{n}\pi_{i^{'}}\mathcal{K}_{k+1}\left(x_{i},\,x_{i^{'}}\right)\right)^{2}+\lambda\pi^{T}\mathbb{K}_{k+1}\pi\label{eq:constraints2}
\end{equation}
for a properly chosen regularization parameter $\lambda>0$ (detailed
in the theorems concerning KRR). Consequently, the optimal weight
vector $\hat{\pi}$ takes the form 
\begin{equation}
\hat{\pi}=\left(\mathbb{K}_{k+1}+\lambda I_{n}\right)^{-1}\frac{Y}{\sqrt{n}}\label{eq:optimal weight}
\end{equation}
where $Y=\left\{ Y_{i}\right\} _{i=1}^{n}$. The KRR estimators associated
with $\mathcal{S}_{k+1}$ are related to the smoothing splines methods
and Gaussian process regressions in machine learning. We refer interested
readers to Wahba (1990), Schölkopf and Smola (2002), as well as Rasmussen
and Williams (2006) for more details.

\begin{theorem}[upper bounds for standard Sobolev]\label{thm:MISE-upper-sobolev-std}\textit{
Suppose Assumption \ref{Assumption 1} holds. }

\textit{(i) If 
\begin{equation}
\frac{n}{\sigma^{2}}>\left(\gamma+1\right)^{2\gamma+3},\label{eq:n-lower-Theorem3.2(ii)}
\end{equation}
let }$\hat{f}\left(\cdot\right)=\frac{1}{\sqrt{n}}\sum_{i^{'}=1}^{n}\hat{\pi}_{i^{'}}\mathcal{K}_{\gamma+1}\left(\cdot,\,x_{i^{'}}\right)$\textit{
where $\left\{ \hat{\pi}_{i^{'}}\right\} _{i^{'}=1}^{n}=\left(\mathbb{K}_{\gamma+1}+\lambda I_{n}\right)^{-1}\frac{Y}{\sqrt{n}}$
and $\lambda\asymp\left(\frac{1}{n}\right)^{\frac{2\left(\gamma+1\right)}{2\left(\gamma+1\right)+1}}$.
Then we have 
\begin{eqnarray*}
\sup_{f\in\mathcal{S}_{\gamma+1}}\mathbb{E}\left(\left|\hat{f}-f\right|_{2,\mathbb{P}}^{2}\right) & \leq & \overline{c}\left[r_{1}^{2}+\exp\left\{ -cnr_{1}^{2}\right\} \right]
\end{eqnarray*}
for some positive universal constants $\overline{c}\in\left(1,\,\infty\right)$
and $c\asymp1$ (both independent of $n$ and $\gamma$ and bounded
away from zero and from above), where $r_{1}^{2}=\left(\frac{\sigma^{2}}{n}\right)^{\frac{2\left(\gamma+1\right)}{2\left(\gamma+1\right)+1}}$. }

\textit{(ii) If 
\begin{equation}
\frac{n}{\sigma^{2}}\leq\left(\gamma+1\right)^{2\gamma+3},\label{eq:n-lower-Theorem3.2(i)}
\end{equation}
let }$\hat{f}=\frac{1}{\sqrt{n}}\sum_{i^{'}=1}^{n}\hat{\pi}_{i^{'}}\mathcal{K}_{\gamma^{*}+1}\left(\cdot,\,x_{i^{'}}\right)$\textit{
where $\left\{ \hat{\pi}_{i^{'}}\right\} _{i^{'}=1}^{n}=\left(\mathbb{K}_{\gamma^{*}+1}+\lambda I_{n}\right)^{-1}\frac{Y}{\sqrt{n}}$
with $\gamma^{*}$ defined in Theorem \ref{thm:MISE-lower-standard},
and $\lambda\asymp\frac{\gamma^{*}+1}{n}$. Then we have 
\[
\sup_{f\in\mathcal{S}_{\gamma+1}}\mathbb{E}\left(\left|\hat{f}-f\right|_{2,\mathbb{P}}^{2}\right)\leq\overline{c}\left[r_{2}^{2}+\exp\left\{ -cnr_{2}^{2}\right\} \right]
\]
where $r_{2}^{2}=\frac{\sigma^{2}\left(\gamma^{*}+1\right)}{n}$. }

\end{theorem}

The proof for Theorem \ref{thm:MISE-upper-sobolev-std} is given in
Appendix \ref{subsec:Proof-for-MISE-upper-sobolev-std}.\textbf{ }

\begin{remark}

Procedures in this paper are standard formulations in modern textbooks
such as Wainwright (2019). Unresolved issues related to practical
implementations (such as the choice of the underived universal constant
in a regularization parameter) persisting in the existing literature
are not addressed in this paper.\footnote{In practical implementation, the specification of the constant in
a regularization parameter is somewhat arbitrary (without a theoretical
justification) in the existing literature. For example, Wainwright
(2019) writes ``user defined radius'' for constants similar to $\overline{C}$,
which becomes part of the universal constant in a regularization parameter.} Obviously, even with a good practical recommendation for the universal
constant in the regularization parameter $\lambda$, the orders of
$\lambda$ and the convergence rates of the upper bounds in Theorem
\ref{thm:MISE-upper-sobolev-std} stay the same. Theorem \ref{thm:MISE-upper-sobolev-std}
contributes to the existing literature in the way that, we show a
phase transition phenomenon in the upper bounds, which match the minimax
lower bounds in Theorem \ref{thm:MISE-lower-standard} in terms of
the rates (the part involving $n$ and $\gamma$ only). This phase
transition phenomenon has been overlooked in both the minimax lower
bounds and the upper bounds in the existing literature.

\end{remark}

If the smoothness degree $\gamma+1$ is known, there is a conjecture
in the literature that cross-validation will yield a regularization
parameter with the optimal order $\left(\frac{1}{n}\right)^{\frac{2\left(\gamma+1\right)}{2\left(\gamma+1\right)+1}}$
that $\lambda$ has in Theorem \ref{thm:MISE-upper-sobolev-std} under
the classical setting where the smoothness degree $\gamma+1$ is finite
and $n\rightarrow\infty$ (e.g., van de Geer, 2000).\footnote{van de Geer (2000) states ``the optimal order'' as $\left(\frac{1}{n}\right)^{\frac{\gamma+1}{2\left(\gamma+1\right)+1}}$
as her $\lambda$ corresponds to our $\sqrt{\lambda}$.} To the best of our knowledge, this conjecture has not been proved
yet. Our results prompt a much harder question of whether cross validation
can also be used to choose a smoothness degree that \textit{adapts
to the phase transition}; that is, in the regime $\frac{n}{\sigma^{2}}\leq\left(\gamma+1\right)^{2\gamma+3}$,
choosing a smoothness degree with the order $\gamma^{*}+1$ has in
Theorem \ref{thm:MISE-upper-sobolev-std}, and in the regime $\frac{n}{\sigma^{2}}>\left(\gamma+1\right)^{2\gamma+3}$,
choosing $\gamma+1$ (also see Section \ref{subsec:Open-questions}
for a related discussion about rate adaptive estimators). 

\textbf{The takeaway for applied researchers}. The conjecture and
question above are mostly of theoretical interest. Even if they are
proved eventually, using cross-validation to choose a smoothness degree
and the corresponding regularization parameter is computationally
costly. For applied researchers, the simplest takeaway from our theoretical
results is similar to the one given in the applied literature --
exploit no more than two degrees of smoothness (see Section \ref{subsec:Practical-implications}
for more discussions).

Like the Sobolev class, we can show a similar achievability result
for the Hölder class. 

\begin{theorem}[upper bounds for standard Hölder]\label{thm:MISE-upper-holder-std}\textit{
Suppose Assumption \ref{Assumption 1} holds. }

\textit{(i) If 
\begin{equation}
\frac{n}{\sigma^{2}}>\left(\gamma+1\right)^{2\gamma+3},\label{eq:sample size-lower-Theorem3.3}
\end{equation}
then we have 
\begin{eqnarray*}
\sup_{f\in\mathcal{U}_{\gamma+1}}\mathbb{E}\left(\left|\hat{f}-f\right|_{2,\mathbb{P}}^{2}\right) & \leq & \overline{c}\left[r_{1}^{2}+\exp\left\{ -cnr_{1}^{2}\right\} \right]
\end{eqnarray*}
for some positive universal constants $\overline{c}\in\left(1,\,\infty\right)$
and $c\asymp1$ (both independent of $n$ and $\gamma$ and bounded
away from zero and from above), where $r_{1}^{2}=\left(\frac{\sigma^{2}}{n}\right)^{\frac{2\left(\gamma+1\right)}{2\left(\gamma+1\right)+1}}$.}

\textit{(ii) If 
\begin{equation}
\frac{n}{\sigma^{2}}\leq\left(\gamma+1\right)^{2\gamma+3},\label{eq:sample size-upper-Theorem3.3}
\end{equation}
let $\hat{f}$ be (\ref{eq:least squares-1}) with $\mathcal{F}=\mathcal{U}_{\gamma^{*}+1}$
with $\gamma^{*}$ defined in Theorem \ref{thm:MISE-lower-standard}.
Then we have 
\[
\sup_{f\in\mathcal{U}_{\gamma+1}}\mathbb{E}\left(\left|\hat{f}-f\right|_{2,\mathbb{P}}^{2}\right)\leq\overline{c}\left[r_{2}^{2}+\exp\left\{ -cnr_{2}^{2}\right\} \right]
\]
where $r_{2}^{2}=\frac{\sigma^{2}\left(\gamma^{*}+1\right)}{n}$. }

\end{theorem}

The proof for Theorem \ref{thm:MISE-upper-holder-std} is given in
Appendix \ref{subsec:Proof-for-MISE-upper-holder-std}.\textbf{ }

Table 1 summarizes the results in Section \ref{subsec:Mean-integrated-squared}
for easy reference. 

\begin{table}[t]
{\small{}\caption{Minimax optimal MISE bounds of the commonly seen $\left(\gamma+1\right)$th
degree Sobolev and Hölder classes}
\medskip{}
}{\small\par}
\centering{}%
\begin{tabular}{c|c|c}
\hline 
 & {\scriptsize{}$\frac{n}{\sigma^{2}}\leq\left(\gamma+1\right)^{2\gamma+3}$} & {\scriptsize{}$\frac{n}{\sigma^{2}}>\left(\gamma+1\right)^{2\gamma+3}$}\tabularnewline
\hline 
{\scriptsize{}MISE} & {\scriptsize{}$\in\left[\underline{c}\frac{\sigma^{2}\left(\gamma^{*}+1\right)}{n},\:c^{'}\frac{\sigma^{2}\left(\gamma^{*}+1\right)}{n}\right]$} & {\scriptsize{}$\in\left[\underline{c}_{0}\left(\frac{\sigma^{2}}{n}\right)^{\frac{2\gamma+2}{2\gamma+3}},\,c^{'}\left(\frac{\sigma^{2}}{n}\right)^{\frac{2\gamma+2}{2\gamma+3}}\right]$}\tabularnewline
\hline 
\end{tabular}\\
{\scriptsize{}where: $\underline{c},\,\underline{c}_{0},\,c^{'}\asymp1$
are positive universal constants that are independent of $n$ and
$\gamma$, and $\underline{c}\leq\frac{\underline{c}_{0}}{2}$; $\gamma^{*}$
is the smallest integer in $\left\{ 0,...,\gamma\right\} $ such that
$\frac{n}{\sigma^{2}}\leq\left(\gamma^{*}+1\right)^{2\gamma^{*}+3}$.}{\scriptsize\par}
\end{table}

\subsection{The sample mean squared error rates}

When deriving the upper bounds in Theorems \ref{thm:MISE-upper-sobolev-std}
and \ref{thm:MISE-upper-holder-std}, we obtain the following bounds
on the sample mean squared error (SMSE) as intermediate results. 

\begin{corollary} \textit{Suppose the conditions in Theorem }\ref{thm:MISE-upper-sobolev-std}\textit{
hold. Under (\ref{eq:n-lower-Theorem3.2(ii)}), we have 
\[
\left|\hat{f}-f\right|_{n}^{2}\precsim r_{1}^{2}\quad\textrm{for any }f\in\mathcal{S}_{\gamma+1},
\]
with probability at least $1-c_{0}\exp\left\{ -cnr_{1}^{2}\right\} $.
Under (\ref{eq:n-lower-Theorem3.2(i)}), we have 
\[
\left|\hat{f}-f\right|_{n}^{2}\precsim r_{2}^{2}\quad\textrm{for any }f\in\mathcal{S}_{\gamma+1},
\]
with probability at least $1-c_{0}\exp\left\{ -cnr_{2}^{2}\right\} $. }

\end{corollary}

\begin{corollary}\textit{ Suppose the conditions in Theorem \ref{thm:MISE-upper-holder-std}
hold. Under (\ref{eq:sample size-lower-Theorem3.3}), we have 
\[
\left|\hat{f}-f\right|_{n}^{2}\precsim r_{1}^{2}\quad\textrm{for any }f\in\mathcal{U}_{\gamma+1},
\]
with probability at least $1-c_{0}\exp\left\{ -cnr_{1}^{2}\right\} $.
Under (\ref{eq:sample size-upper-Theorem3.3}), we have 
\[
\left|\hat{f}-f\right|_{n}^{2}\precsim r_{2}^{2}\quad\textrm{for any }f\in\mathcal{U}_{\gamma+1},
\]
with probability at least $1-c_{0}\exp\left\{ -cnr_{2}^{2}\right\} $. }

\end{corollary}

\subsection{Infinitely smooth functions \label{subsec:Analytic-functions-as}}

In this subsection, we consider infinitely differentiable functions
with the following form
\begin{equation}
f(x)=f(0)+\sum_{k=1}^{\infty}\frac{x^{k}}{k!}f^{(k)}(0).\label{eq:expansion-1}
\end{equation}
The standard $\mathcal{U}_{\infty}$ for infinitely smooth functions
follows the definition of the standard $\mathcal{U}_{\gamma+1}$ with
$\gamma+1=\infty$. The analog of $\mathcal{S}_{\gamma+1}$ in (\ref{eq:sobolev_std})
for infinitely smooth functions can be defined as 
\begin{align*}
\mathcal{S}_{\infty} & :=\{f:\,\left[0,\,1\right]\rightarrow\mathbb{R}\vert\,f\textrm{ is infinitely differentiable,}\\
 & \sum_{j=0}^{\infty}\left(f^{(j)}(0)\right)^{2}\leq\overline{C}^{2}\}
\end{align*}
where $\overline{C}\asymp1$ is a universal constant independent of
$n$ and $\gamma$. 

We have the following result for analytic functions. 

\begin{corollary}[analytic functions] \label{corr:analytic functions}\textit{Suppose
Assumption \ref{Assumption 1} holds. For the minimax lower bounds,
we further assume that the density $p(x)$ is bounded away from zero;
that is, $p(x)\geq c>0$ for some universal constant $c$. Given a
sample size $n$, let $\gamma^{*}\in\left\{ 0,...,\gamma\right\} $
be the smallest integer such that $\frac{n}{\sigma^{2}}\leq\left(\gamma^{*}+1\right)^{2\gamma^{*}+3}$.
Then we have 
\begin{eqnarray*}
\inf_{\tilde{f}}\sup_{f\in\mathcal{S}_{\infty}}\mathbb{E}\left(\left|\tilde{f}-f\right|_{2,\mathbb{P}}^{2}\right) & \geq & \underline{c}r_{2}^{2},\\
\sup_{f\in\mathcal{S}_{\infty}}\mathbb{E}\left(\left|\hat{f}-f\right|_{2,\mathbb{P}}^{2}\right) & \leq & \overline{c}\left[r_{2}^{2}+\exp\left\{ -cnr_{2}^{2}\right\} \right],
\end{eqnarray*}
for some positive universal constants $\underline{c},\,\overline{c},\,c\asymp1$
that are independent of $n$ and $\gamma$, where $\hat{f}$ and }$r_{2}$\textit{
are defined in Theorem \ref{thm:MISE-upper-sobolev-std}(ii).}

\textit{Let $\hat{f}$ be (\ref{eq:least squares-1}) with $\mathcal{F}=\mathcal{U}_{\gamma^{*}+1}$.
Then we also have }

\textit{
\begin{eqnarray*}
\inf_{\tilde{f}}\sup_{f\in\mathcal{U}_{\infty}}\mathbb{E}\left(\left|\tilde{f}-f\right|_{2,\mathbb{P}}^{2}\right) & \geq & \underline{c}r_{2}^{2},\\
\sup_{f\in\mathcal{U}_{\infty}}\mathbb{E}\left(\left|\hat{f}-f\right|_{2,\mathbb{P}}^{2}\right) & \leq & \overline{c}\left[r_{2}^{2}+\exp\left\{ -cnr_{2}^{2}\right\} \right],
\end{eqnarray*}
for some positive universal constants $\underline{c},\,\overline{c},\,c\asymp1$
that are independent of $n$ and $\gamma$.}

\end{corollary}

The proof for the lower bound in Corollary \ref{corr:analytic functions}
is identical to the proof for part (ii) of Theorem \ref{thm:MISE-lower-standard}
(Appendix \ref{subsec:Proof-for-MISE-lower-standard}). The proofs
for the upper bounds in Corollary \ref{corr:analytic functions} involve
only slight modifications of Appendix \ref{subsec:Proof-for-MISE-upper-sobolev-std}
and Appendix \ref{subsec:Proof-for-MISE-upper-holder-std}. See Appendix
\ref{subsec:Proof-for-analytic_MISE}.

There is no phase transition in the optimal rate associated with $\mathcal{U}_{\infty}$
and $\mathcal{S}_{\infty}$. As discussed earlier, $\frac{\gamma^{*}+1}{n}\asymp\frac{\log n}{n\log(\log n)}$.
Relatedly, see a review article about the estimation of analytic functions
by Ibragimov (2001), in particular, Theorem 4.9 where the minimax
optimal rate for the MISE associated with analytic functions is $\frac{\log n}{n\log(\log n)}$.\footnote{The author thanks an anonymous referee for pointing out this reference.} 

Corollary \ref{corr:analytic functions} highlights the importance
of our results on phase transition in this paper. Stating $\left(\frac{1}{n}\right)^{\frac{2\gamma+2}{2\gamma+3}}$
as the minimax optimal rate for MISE associated with a standard smoothness
class without a condition on $n$ raises the following ``paradox''.
Recall (\ref{eq:asymptotic interpretation}); on the other hand, $\mathcal{U}_{\infty}\subseteq\mathcal{U}_{\gamma+1}$
so (\ref{eq:asymptotic interpretation}) is not plausible if $\left(\frac{1}{n}\right)^{\frac{2\gamma+2}{2\gamma+3}}$
is stated \textit{alone}, without a condition on $n$, as the minimax
optimal rate in MISE for the standard $\mathcal{U}_{\gamma+1}$. Our
results in Section \ref{subsec:Mean-integrated-squared} are crucial
for reconciling the aforementioned ``paradox'', as the choice of
our critical smoothness parameter $\gamma^{*}+1$ guarantees that: 
\begin{itemize}
\item in the ``small $n$'' regime, the minimax optimal rate associated
with $\mathcal{U}_{\gamma+1}$ is $\frac{\gamma^{*}+1}{n}\asymp\frac{\log n}{n\log(\log n)}$; 
\item in the ``large $n$'' regime, the minimax optimal rate associated
with $\mathcal{U}_{\gamma+1}$ is $\left(\frac{1}{n}\right)^{\frac{2\gamma+2}{2\gamma+3}}$
and 
\[
\left(\frac{1}{n}\right)^{\frac{2\gamma+2}{2\gamma+3}}\succsim\frac{\gamma+1}{n}\geq\frac{\gamma^{*}+1}{n}\asymp\frac{\log n}{n\log(\log n)}.
\]
\end{itemize}
Assuming $\sigma\asymp1$, a better practice would be to state the
rate $\left(\frac{\sigma^{2}}{n}\right)^{\frac{2\gamma+2}{2\gamma+3}}$
under the condition $\frac{n}{\sigma^{2}}>\left(\gamma+1\right)^{2\gamma+3}$
and the rate $\frac{\sigma^{2}\log n}{n\log(\log n)}$ under the condition
$\frac{n}{\sigma^{2}}\leq\left(\gamma+1\right)^{2\gamma+3}$, as shown
in this paper. 

\section{Covering and packing numbers \label{sec:Covering-and-packing}}

In this section, we present bounds on covering and packing numbers
associated with $\mathcal{U}_{\gamma+1,1}$, $\mathcal{U}_{\gamma+1,2}$,
and $\mathcal{H}_{\gamma+1}$. Various $c$ and $C$ letters in this
section denote positive universal constants that are: finite and bounded
away from zero (denoted by $\asymp1$) and independent of $\gamma$
and $\left\{ R_{k}\right\} _{k=0}^{\gamma+1}$ (parameters of the
function classes); these constants may vary from place to place. 

Table 2 summarizes the results in this section for easy reference.
\begin{table}
{\footnotesize{}\caption{Upper and lower bounds on the $\log(\delta-\textrm{covering number})$
and $\log(\delta-\textrm{packing number})$ of the generalized $\mathcal{U}_{\gamma+1,1}$,
$\mathcal{U}_{\gamma+1,2}$ and $\mathcal{H}_{\gamma+1}$ in $L^{q}-$norm}
}{\small{}\medskip{}
}{\small\par}
\begin{centering}
{\scriptsize{}}%
\begin{tabular}{c|c||c||c||c||c}
\hline 
\multicolumn{1}{c}{} & \multicolumn{2}{c||}{{\scriptsize{}$\mathcal{U}_{\gamma+1,1}$ ($q\in\left\{ 2,\,\infty\right\} $)}} & \multicolumn{2}{c||}{{\scriptsize{}$\mathcal{U}_{\gamma+1,2}$ ($q\in\left\{ 2,\,\infty\right\} $)}} & {\scriptsize{}$\mathcal{H}_{\gamma+1}$ ($q=2$)}\tabularnewline
\hline 
{\scriptsize{}$\precsim$} & \multicolumn{2}{c||}{{\scriptsize{}$\begin{cases}
\overline{B}_{1}\left(\delta\right) & \textrm{if }\min_{k\in\left\{ 0,...,\gamma\right\} }\log\frac{4\left(\gamma+1\right)R_{k}}{k!\delta}\geq0\\
\overline{B}_{2}\left(\delta\right) & \textrm{otherwise}
\end{cases}$}} & \multicolumn{2}{c||}{{\scriptsize{}$R^{*\frac{1}{\gamma+1}}\delta^{\frac{-1}{\gamma+1}}$}} & {\scriptsize{}$\begin{cases}
R_{\gamma+1}^{\frac{1}{\gamma+1}}\delta^{\frac{-1}{\gamma+1}} & \textrm{if }R_{\gamma+1}\succsim\gamma+1\\
\delta^{\frac{-1}{\gamma+1}} & \textrm{if }R_{\gamma+1}\precsim\gamma+1
\end{cases}$}\tabularnewline
\hline 
{\scriptsize{}$\succsim$} & \multicolumn{2}{c||}{{\scriptsize{}$\max\left\{ \underline{B}_{1}\left(\delta\right),\,\underline{B}_{2}\right\} $}} & \multicolumn{2}{c||}{{\scriptsize{}$\begin{cases}
R^{*\frac{1}{\gamma+1}}\delta^{\frac{-1}{\gamma+1}} & \textrm{if }R_{0}\succsim1\\
\left(R^{*}R_{0}\right)^{\frac{1}{\gamma+1}}\delta^{\frac{-1}{\gamma+1}} & \textrm{if }R_{0}\precsim1
\end{cases}$}} & {\scriptsize{}$R_{\gamma+1}^{\frac{1}{\gamma+1}}\delta^{\frac{-1}{\gamma+1}}$}\tabularnewline
\hline 
\end{tabular}{\scriptsize\par}
\par\end{centering}
\centering{}{\scriptsize{}where: $\overline{B}_{1}\left(\delta\right)=\sum_{k=0}^{\gamma}\log\frac{4\left(\gamma+1\right)R_{k}}{k!\delta}$;
$\overline{B}_{2}\left(\delta\right)=\left(\frac{\gamma}{2}+1\right)\log\frac{1}{\delta}+\sum_{k=0}^{\gamma}\log R_{k}$;
$\underline{B}_{1}\left(\delta\right)=\sum_{k=0}^{\gamma}\log\left(9^{-\gamma}\gamma^{-\gamma}\right)+\sum_{k=0}^{\gamma}\log\frac{C\sum_{m=0}^{\left\lfloor \gamma/2\right\rfloor }R_{k+2m}}{\delta}$
(with $R_{k+2m}=0$ for $k+2m>\gamma$); $\underline{B}_{2}=C^{'}\left(\gamma+1\right)$
(valid for all $\delta$ below a threshold detailed in Lemma 3.1);
$R^{*}=\left(\max_{k\in\left\{ 1,...,\gamma+1\right\} }\frac{R_{k}}{\left(k-1\right)!}\right)\vee1$;
$C$ and $C^{'}$ are positive universal constants that are: $\asymp1$,
independent of $\gamma$ and $\left\{ R_{k}\right\} _{k=0}^{\gamma+1}$.}{\scriptsize\par}
\end{table}

\subsection{The generalized polynomial subclass, $\mathcal{U}_{\gamma+1,1}$}

\begin{lemma}\label{lm:entropy-generalized-polynomial} \textit{(i)
If $\delta$ is small enough such that $\min_{k\in\left\{ 0,...,\gamma\right\} }\log\frac{4\left(\gamma+1\right)R_{k}}{k!\delta}\geq0$,
we have
\begin{align}
\log N_{2,\mathbb{P}}\left(\delta,\,\mathcal{U}_{\gamma+1,1}\right) & \leq\log N_{\infty}\left(\delta,\,\mathcal{U}_{\gamma+1,1}\right)\leq\underset{\overline{B}_{1}\left(\delta\right)}{\underbrace{\sum_{k=0}^{\gamma}\log\frac{4\left(\gamma+1\right)R_{k}}{k!\delta}}};\label{eq:b1_upper_Lemma3.1}
\end{align}
if $\delta$ is large enough such that $\min_{k\in\left\{ 0,...,\gamma\right\} }\log\frac{4\left(\gamma+1\right)R_{k}}{k!\delta}<0$,
we have
\begin{equation}
\log N_{2,\mathbb{P}}\left(\delta,\,\mathcal{U}_{\gamma+1,1}\right)\leq\log N_{\infty}\left(\delta,\,\mathcal{U}_{\gamma+1,1}\right)\leq\underset{\overline{B}_{2}\left(\delta\right)}{\underbrace{\left(\frac{\gamma}{2}+1\right)\log\frac{1}{\delta}+\sum_{k=0}^{\gamma}\log R_{k}}}.\label{eq:Kolmogorov_upper-1-1}
\end{equation}
}

\textit{(ii) In terms of the lower bounds, we have 
\begin{eqnarray*}
\log M_{2}\left(\delta,\,\mathcal{U}_{\gamma+1,1}\right) & \geq & \underline{B}_{1}\left(\delta\right),\\
\log M_{\infty}\left(\delta,\,\mathcal{U}_{\gamma+1,1}\right) & \succsim & \underline{B}_{1}\left(\delta\right),
\end{eqnarray*}
where $\underline{B}_{1}\left(\delta\right)=\sum_{k=0}^{\gamma}\log\left(9^{-\gamma}\gamma^{-\gamma}\right)+\sum_{k=0}^{\gamma}\log\frac{C\sum_{m=0}^{\left\lfloor \gamma/2\right\rfloor }R_{k+2m}}{\delta}$
(with $R_{k+2m}=0$ for $k+2m>\gamma$) for some positive universal
constant $C\asymp1$ independent of $\gamma$ and $\left\{ R_{k}\right\} _{k=0}^{\gamma}$.
Let $\tilde{k}\in\arg\max_{k\in\left\{ 0,...,\gamma\right\} }\frac{R_{k}}{k!}$.
If 
\begin{equation}
\frac{R_{\tilde{k}}}{\tilde{k}!\delta\left[\left(\tilde{k}+1\right)\vee\sum_{k=0}^{\gamma}\frac{R_{k}}{k!}\right]}\succsim2^{\gamma+1},\label{eq:6}
\end{equation}
we also have 
\begin{eqnarray}
\log M_{2}\left(\delta,\,\mathcal{U}_{\gamma+1,1}\right) & \geq & \underline{B}_{2}=C^{'}\left(\gamma+1\right),\label{eq:10}\\
\log M_{\infty}\left(\delta,\,\mathcal{U}_{\gamma+1,1}\right) & \succsim & \underline{B}_{2},\nonumber 
\end{eqnarray}
for some positive universal constant $C^{'}\asymp1$ independent of
$\gamma$ and $\left\{ R_{k}\right\} _{k=0}^{\gamma}$.}

\textit{(iii) If the density function $p(x)$ on $\left[-1,\,1\right]$
is bounded away from zero, i.e., $p(x)\geq c>0$, then
\begin{equation}
\log M_{2,\mathbb{P}}\left(\delta,\,\mathcal{U}_{\gamma+1,1}\right)\succsim\underline{B}_{1}\left(\delta\right);\label{eq:b1_lower_Lemma3.1}
\end{equation}
under (\ref{eq:6}), we also have 
\begin{equation}
\log M_{2,\mathbb{P}}\left(\delta,\,\mathcal{U}_{\gamma+1,1}\right)\succsim\underline{B}_{2}.\label{eq:b2_lower_Lemma3.1}
\end{equation}
}

\end{lemma}

\begin{remark} When $R_{k}=\overline{C}$ for $k=0,...,\gamma$,
$\left(\tilde{k}+1\right)\vee\sum_{k=0}^{\gamma}\frac{R_{k}}{k!}\asymp1$;
when $R_{0}=\overline{C}$ and $R_{k}\leq\overline{C}\left(k-1\right)!$
for $k=1,...,\gamma$, $\left(\tilde{k}+1\right)\vee\sum_{k=0}^{\gamma}\frac{R_{k}}{k!}\precsim\log\left(\gamma\vee2\right)$;
when $R_{k}=\overline{C}k!$ for all $k=0,...,\gamma$, $\left(\tilde{k}+1\right)\vee\sum_{k=0}^{\gamma}\frac{R_{k}}{k!}\asymp\left(\gamma\vee1\right)$.

\end{remark}

The proof for Lemma \ref{lm:entropy-generalized-polynomial} is given
in Appendix \ref{subsec:Proof-for-entropy-poly}. 

The lower bounds $\underline{B}_{1}\left(\delta\right)$ and $\underline{B}_{2}$,
as well as the upper bound $\overline{B}_{1}\left(\delta\right)$
are original. The (less original) bound $\overline{B}_{2}\left(\delta\right)$
generalizes the upper bound associated with the polynomial subclass
in Kolmogorov and Tikhomirov (1959), which takes the form $\left(\gamma+1\right)\log\frac{1}{\delta}$.
It is worth pointing out that $\overline{B}_{2}\left(\delta\right)$
holds for all $\delta\in(0,\,1)$ (not just $\delta$ such that \textit{$\min_{k\in\left\{ 0,...,\gamma\right\} }\log\frac{4\left(\gamma+1\right)R_{k}}{k!\delta}<0$})
but is far from being tight when $\min_{k\in\left\{ 0,...,\gamma\right\} }\log\frac{4\left(\gamma+1\right)R_{k}}{k!\delta}\geq0$.
Obviously $\overline{B}_{1}\left(\delta\right)\precsim\overline{B}_{2}\left(\delta\right)$.
When it comes to deriving the minimax optimal rates for the MISE under
large enough $R_{k}$ (e.g., Theorems \ref{thm:MISE_(k-1)!} and \ref{thm:MISE_k!}
in Appendix \ref{sec:MISE-nonstandard}), $\overline{B}_{1}\left(\delta\right)$
will be very useful. 

When deriving a lower bound for the packing number of the \textit{standard}
$\mathcal{U}_{\gamma+1}$ under the assumption that $R_{k}=\overline{C}\asymp1$,
Kolmogorov and Tikhomirov (1959) constructs a set of functions where
$\mathcal{U}_{\gamma+1,1}$ is a singleton; therefore, the cardinality
of this set only gives a lower bound for the packing number of $\mathcal{U}_{\gamma+1,2}$
and overlooks $\mathcal{U}_{\gamma+1,1}$. This issue is further discussed
in Section \ref{subsec:entropy-generalized-H=0000F6lder}. Our Lemma
\ref{lm:entropy-generalized-polynomial} establishes two different
lower bounds for the packing number of $\mathcal{U}_{\gamma+1,1}$.
In particular, the construction of $\underline{B}_{2}$ in (\ref{eq:b2_lower_Lemma3.1})
will be useful for deriving the minimax lower bounds for the MISE
when $R_{k}$ is relatively small (e.g., Theorem \ref{thm:MISE-lower-standard}
in Section \ref{sec:Minimax-standard} and Theorem \ref{thm:MISE_(k-1)!}
in Appendix \ref{sec:MISE-nonstandard}), while $\underline{B}_{1}\left(\delta\right)$
in (\ref{eq:b1_upper_Lemma3.1}) will be useful when $R_{k}$ is relatively
large (Theorem \ref{thm:MISE_k!} in Appendix \ref{sec:MISE-nonstandard}). 

To establish $\overline{B}_{1}\left(\delta\right)$, $\underline{B}_{1}\left(\delta\right)$
and $\underline{B}_{2}$, we discard the argument in Kolmogorov and
Tikhomirov (1959) and develop our own. The derivation of $\underline{B}_{2}$
is based on a constructive proof. To derive $\overline{B}_{1}\left(\delta\right)$
and $\underline{B}_{1}\left(\delta\right)$, we consider two classes
(equivalent to $\mathcal{U}_{\gamma+1,1}$), each involving a $\left(\gamma+1\right)-$dimensional
polyhedron. The lower bound $\underline{B}_{1}\left(\delta\right)$
is the more delicate part. In particular, for any $f\in\mathcal{U}_{\gamma+1,1}$,
we write $f\left(x\right)=\sum_{k=0}^{\gamma}\tilde{\theta}_{k}\phi_{k}\left(x\right)$,
where $\left(\phi_{k}\right)_{k=0}^{\gamma}$ are the Legendre polynomials. 

\subsection{The generalized Hölder subclass, $\mathcal{U}_{\gamma+1,2}$\label{subsec:entropy-generalized-H=0000F6lder}}

\begin{lemma} \textit{\label{lm:entropy-generalized-H=0000F6lder}Let
$R^{*}=\left(\max_{k\in\left\{ 1,...,\gamma+1\right\} }\frac{R_{k}}{\left(k-1\right)!}\right)\vee1$.
We have 
\begin{align*}
\log N_{2,\mathbb{P}}\left(\delta,\,\mathcal{U}_{\gamma+1,2}\right) & \leq\log N_{\infty}\left(\delta,\,\mathcal{U}_{\gamma+1,2}\right)\precsim R^{*\frac{1}{\gamma+1}}\delta^{\frac{-1}{\gamma+1}}.
\end{align*}
We also have
\begin{eqnarray*}
\log M_{\infty}\left(\delta,\,\mathcal{U}_{\gamma+1,2}\right) & \succsim & \log M_{2}\left(\delta,\,\mathcal{U}_{\gamma+1,2}\right)\succsim R^{*\frac{1}{\gamma+1}}\delta^{\frac{-1}{\gamma+1}},\quad\textrm{if }R_{0}\succsim1,\,\delta\in(0,\,1);\\
\log M_{\infty}\left(\delta,\,\mathcal{U}_{\gamma+1,2}\right) & \succsim & \log M_{2}\left(\delta,\,\mathcal{U}_{\gamma+1,2}\right)\succsim\left(R^{*}R_{0}\right)^{\frac{1}{\gamma+1}}\delta^{\frac{-1}{\gamma+1}},\quad\textrm{if }R_{0}\precsim1,\,\delta\in(0,\,1).
\end{eqnarray*}
If the density function $p(x)$ on $\left[-1,\,1\right]$ is bounded
away from zero, i.e., $p(x)\geq c>0$, then 
\begin{eqnarray*}
\log M_{2,\mathbb{P}}\left(\delta,\,\mathcal{U}_{\gamma+1,2}\right) & \succsim & R^{*\frac{1}{\gamma+1}}\delta^{\frac{-1}{\gamma+1}},\quad\textrm{if }R_{0}\succsim1,\,\delta\in(0,\,1);\\
\log M_{2,\mathbb{P}}\left(\delta,\,\mathcal{U}_{\gamma+1,2}\right) & \succsim & \left(R^{*}R_{0}\right)^{\frac{1}{\gamma+1}}\delta^{\frac{-1}{\gamma+1}},\quad\textrm{if }R_{0}\precsim1,\,\delta\in(0,\,1).
\end{eqnarray*}
}

\end{lemma}

The proof for Lemma \ref{lm:entropy-generalized-H=0000F6lder} is
given in Appendix \ref{subsec:Proof-for-entropy-holder-sub}. Lemma
\ref{lm:entropy-generalized-polynomial} and Lemma \ref{lm:entropy-generalized-H=0000F6lder}
together are used to prove Theorems \ref{thm:MISE-lower-standard}
and \ref{thm:MISE-upper-holder-std} in Section \ref{sec:Minimax-standard},
as well as Theorems \ref{thm:MISE_(k-1)!} and \ref{thm:MISE_k!}
in Appendix \ref{sec:MISE-nonstandard}; in addition, Lemma \ref{lm:entropy-generalized-H=0000F6lder}
is also used to prove Theorem \ref{thm:MISE-gen-sub-ellipsoid}(ii)
in Appendix \ref{sec:MISE-nonstandard}.

Lemma \ref{lm:entropy-generalized-H=0000F6lder} extends Kolmogorov
and Tikhomirov (1959) to allow for general $R_{k}$s. When $R_{k}\leq\overline{C}k!$
for all $k=1,...,\gamma+1$, $R^{*\frac{1}{\gamma+1}}\asymp1$. If
$R_{k}\geq\overline{C}k!$ for all $k=1,...,\gamma+1$, $R^{*\frac{1}{\gamma+1}}\succsim1$;
for example, taking $R_{k}\geq\overline{C}\left(k!\right)^{2}$ for
all $k=1,...,\gamma+1$ yields $R^{*\frac{1}{\gamma+1}}\succsim\gamma$. 

Given Lemmas \ref{lm:entropy-generalized-polynomial} and \ref{lm:entropy-generalized-H=0000F6lder},
(\ref{eq:inner}) and (\ref{eq:outter}), we have 
\begin{align}
\log N_{2,\mathbb{P}}\left(2\delta,\,\mathcal{U}_{\gamma+1}\right) & \leq\log N_{\infty}\left(2\delta,\,\mathcal{U}_{\gamma+1}\right)\nonumber \\
 & \leq\begin{cases}
\overline{B}_{1}\left(\delta\right)+R^{*\frac{1}{\gamma+1}}\delta^{\frac{-1}{\gamma+1}} & \textrm{if }\min_{k\in\left\{ 0,...,\gamma\right\} }\log\frac{4\left(\gamma+1\right)R_{k}}{k!\delta}\geq0\\
\overline{B}_{2}\left(\delta\right)+R^{*\frac{1}{\gamma+1}}\delta^{\frac{-1}{\gamma+1}} & \textrm{if }\min_{k\in\left\{ 0,...,\gamma\right\} }\log\frac{4\left(\gamma+1\right)R_{k}}{k!\delta}<0
\end{cases}\label{eq:upper_together}
\end{align}
and
\[
\log M_{\infty}\left(\delta,\,\mathcal{U}_{\gamma+1}\right)\succsim\log M_{2}\left(\delta,\,\mathcal{U}_{\gamma+1}\right)\succsim\begin{cases}
\max\left\{ \underline{B}_{1}\left(\delta\right),\,\underline{B}_{2},\,R^{*\frac{1}{\gamma+1}}\delta^{\frac{-1}{\gamma+1}}\right\}  & \textrm{if }R_{0}\succsim1\\
\max\left\{ \underline{B}_{1}\left(\delta\right),\,\underline{B}_{2},\,\left(R^{*}R_{0}\right)^{\frac{1}{\gamma+1}}\delta^{\frac{-1}{\gamma+1}}\right\}  & \textrm{if }R_{0}\precsim1
\end{cases}.
\]

Our lower bounds above sharpen the classical result in Kolmogorov
and Tikhomirov (1959). In particular, the lower bound for $\mathcal{U}_{\gamma+1}$
in Kolmogorov and Tikhomirov (1959) (derived under the assumption
that $R_{k}=\overline{C}\asymp1$) takes the form $\delta^{\frac{-1}{\gamma+1}}$.
This result and its proof are inherited later in papers and textbooks
including the more recent textbook on nonasymptotic statistics by
Wainwright (2019, Example 5.11), where in the derivation of the lower
bound, a set of functions are constructed in the way such that their
$k$th order derivatives evaluated at zero are zero for all $k=0,...,\gamma$.
In other words, $\mathcal{U}_{\gamma+1,1}$ is a singleton in this
construction and the cardinality of this set only gives a lower bound
for $\mathcal{U}_{\gamma+1,2}$. In particular, the lower bound $\delta^{\frac{-1}{\gamma+1}}$
is not sharp when $\gamma+1$ and $\delta$ are large enough. 

\subsection{The ellipsoid subclass, $\mathcal{H}_{\gamma+1}$}

\begin{lemma}\label{lm:entropy-ellipsoid-subclass} \textit{Assume
$\mu_{m}=\left(cm\right)^{-2\left(\gamma+1\right)}$ in (\ref{eq:Ellipsoid-1})
for a positive constant $c\asymp1$ independent of $\gamma$ and $R_{\gamma+1}$.
If $R_{\gamma+1}\succsim\gamma+1$, we have
\[
\log N_{2}\left(\delta,\,\mathcal{H}_{\gamma+1}\right)\asymp\left(R_{\gamma+1}\delta^{-1}\right)^{\frac{1}{\gamma+1}}.
\]
If $R_{\gamma+1}\precsim\gamma+1$, we have
\begin{align}
\log N_{2}\left(\delta,\,\mathcal{H}_{\gamma+1}\right) & \precsim\delta^{\frac{-1}{\gamma+1}},\label{eq:15-1}\\
\log N_{2}\left(\delta,\,\mathcal{H}_{\gamma+1}\right) & \succsim\left(R_{\gamma+1}\delta^{-1}\right)^{\frac{1}{\gamma+1}}.\label{eq:15}
\end{align}
If the density function $p(x)$ on $\left[0,\,1\right]$ is bounded
away from zero, i.e., $p(x)\geq c>0$, then the bounds above also
hold for $\log N_{2,\mathbb{P}}\left(\delta,\,\mathcal{H}_{\gamma+1}\right)$.}

\end{lemma}

The proof for Lemma \ref{lm:entropy-ellipsoid-subclass} is given
in Appendix \ref{subsec:Proof-for-entropy_ellipsoid_sub}. Lemma \ref{lm:entropy-generalized-polynomial}
and Lemma \ref{lm:entropy-ellipsoid-subclass} together are used to
prove Theorem \ref{thm:MISE-lower-standard}. Lemma \ref{lm:entropy-ellipsoid-subclass}
is also used to prove Theorem \ref{thm:MISE-gen-sub-ellipsoid}(i)
in Appendix \ref{sec:MISE-nonstandard}. 

When $R_{\gamma+1}=1$, Lemma \ref{lm:entropy-ellipsoid-subclass}
sharpens the upper bound for $\log N_{2}\left(\delta,\,\mathcal{H}_{\gamma+1}\right)$
in Wainwright (2019) from $\left(\gamma\vee1\right)\delta^{-\frac{1}{\gamma+1}}$
to $\delta^{-\frac{1}{\gamma+1}}$; in particular, the upper and lower
bounds in Wainwright (2019) (the last two inequalities on p.131) scale
as $\left(\gamma\vee1\right)\delta^{\frac{-1}{\gamma+1}}$ and $\delta^{\frac{-1}{\gamma+1}}$,
respectively, while our upper and lower bounds in Lemma \ref{lm:entropy-ellipsoid-subclass}
have the same scaling $\delta^{\frac{-1}{\gamma+1}}$. We discover
the cause of the gap lies in that the ``pivotal'' eigenvalue (that
balances the ``estimation error'' and the ``approximation error''
from truncating for a given resolution $\delta$) in Wainwright (2019)
is not optimal. The truncation in Wainwright (2019) is commonly used
in the existing literature and seems to originate from Theorem 3 in
Mityagin (1961). We close the gap by finding the optimal ``pivotal''
eigenvalue.

More generally, for the case of $R_{\gamma+1}\precsim\gamma+1$, we
consider two different truncations, one giving the upper bound $\delta^{\frac{-1}{\gamma+1}}$
and the other giving the lower bound $\left(R_{\gamma+1}\delta^{-1}\right)^{\frac{1}{\gamma+1}}$.
Note that $\left(R_{\gamma+1}\delta^{-1}\right)^{\frac{1}{\gamma+1}}\asymp\delta^{\frac{-1}{\gamma+1}}$
when $R_{\gamma+1}\asymp1$. For the case of $R_{\gamma+1}\succsim\gamma+1$,
we use only one truncation to show that both the upper bound and the
lower bound scale as $\left(R_{\gamma+1}\delta^{-1}\right)^{\frac{1}{\gamma+1}}$.

\section{Discussions\label{sec:Discussions}}

\subsection{Recommendation from the applied literature \label{subsec:Practical-implications}}

In the discussion following Theorem \ref{thm:MISE-lower-standard},
we have brought up the trade-off between the polynomial subclass and
the nonparametric subclass. It is worth mentioning the connection
between our theoretical results and the empirical findings from Gelman
and Imbens (2019), which implements high order polynomials in regression
discontinuity designs (RDD) analyses. In particular, when applying
RDD to perform causal inference, two conditional mean functions of
a pretreatment variable are estimated from (\ref{eq:model}). There
are several empirical issues of using high order polynomials raised
in Gelman and Imbens (2019). The implication of our results is most
related to their paper on the issue of mean squared errors (MSE).

Regarding the data sets studied in Gelman and Imbens (2019), the Jacob-Lefgren
data (Jacob and Lefgren, 2004), Lee data (Lee, 2008), Matsudaira data
(Matsudaira, 2008), the LaLonde data (LaLonde, 1986), and the census
data in 1974, 1975 and 1978, the sample sizes used in the implementation
of Gelman and Imbens (2019) range from thousands to at most thirties
of thousands. Based on their empirical evidence from studying these
data sets, Gelman and Imbens (2019) recommend researchers to avoid
using high order polynomials but use local linear or local quadratic
polynomials. In view of the asymptotic condition $\left(\gamma+1\right)^{2\gamma+3}=o(n)$
for the rate $\left(\frac{1}{n}\right)^{\frac{2\gamma+2}{2\gamma+3}}$
to kick in (see (\ref{eq:asymptotic condition_1})), let us solve
for $\overline{\gamma}+1$ roughly from $n\geq\left(\overline{\gamma}+1\right)^{2\overline{\gamma}+3}$.
This \textit{heuristic} gives a rough degree of smoothness to exploit,
$\max\left\{ \left\lfloor \frac{\log n}{2\log\left(\log n\right)}\right\rfloor ,\,1\right\} $,
which is quite close to the recommended smoothness degrees in Gelman
and Imbens (2019).

\subsection{Open questions\label{subsec:Open-questions}}

We conclude the paper by discussing a few open questions motivated
by this work. First, our focus in this paper is on the global criterion
MISE while the concern of Gelman and Imbens (2019) is about the MSE
of the high order polynomial implementation at a point. The minimax
optimality of pointwise MSE and the global MISE (concerning an entire
function) would involve different proofs. It is well known that the
minimax optimal MISE rate coincides with the minimax optimal pointwise
MSE rate in the regime where $\gamma$ is finite and $n\rightarrow\infty$
(see Tsybakov, 2009). We conjecture that the minimax optimal rates
would be the same for the MISE and pointwise MSE in our ``small $n$''
regime, simply because the trade-off between the polynomial subclass
and the nonparametric subclass exists whether the interest is the
MISE or pointwise MSE.

Second, discussions of related literature in Section \ref{sec:Introduction}
indicate that the classical rate $\left(\frac{1}{n}\right)^{\frac{2\gamma+2}{2\gamma+3}}$
is an underestimate of the MISE for local smoothing methods such as
kernel density estimators and local polynomials when $n$ is not large
enough. For this problem, we could consider (\ref{eq:model}) where
$X_{i}=\frac{i}{n}$ for $i=1,...,n$ and $\left\{ \varepsilon_{i}\right\} _{i=1}^{n}$
satisfies the assumptions in Corollary 2.3 of Tsybakov (2009). This
setup is simpler than the one considered in this paper, but serves
a good starting point. Like how we establish the results in this paper,
we would first show the minimax lower bound under the ``small $n$''
regime, and then show that the MISE of a local smoothing method has
an upper bound that matches the lower bound up to some universal constant
independent of $n$ and $\gamma$. The proofs would be different from
the ones in this paper. There is some theoretical evidence (although
not a proof) suggesting that it would require a large $n$ for higher
order local polynomials to become beneficial; for example, Tsybakov
(2009) requires the smallest eigenvalue associated with the local
polynomials to be bounded away from zero (Assumption LP1) to establish
the upper bound $\left(\frac{1}{n}\right)^{\frac{2\gamma+2}{2\gamma+3}}$.
This eigenvalue condition in Tsybakov (2009) requires a large enough
$n$ and a sufficient condition given in Tsybakov (2009) is that $n\rightarrow\infty$.

Third, our results prompt a hard question of whether a practical estimator
can be developed to \textit{adapt to the phase transition} shown in
this paper. Cross validation is one possibility as discussed in Section
\ref{subsec:Mean-integrated-squared}. In the literature, an alternative
construction of adaptive estimators uses the Lepski's method; see,
for example, Chen et. al (2021) where the unknown degree of smoothness
is assumed to be fixed while $n\rightarrow\infty$. Hence, the construction
of smoothness degrees in Chen et. al (2021) does not depend on $n$.
Chen raises an interesting point: a better construction of smoothness
degrees should depend on $n$, based on our results in this paper
(personal communication, October and November 2022). However, it is
not clear if this approach will yield an estimator that automatically
\textit{switches} the degree of smoothness to adapt to the phase transition
in the rates, which, after all, is a nonlinear phenomenon. 

\newpage{}

\appendix

\part*{Appendices for ``Phase transitions in nonparametric regressions''}

\section{Minimax optimal rates in non-standard cases\label{sec:MISE-nonstandard}}

In this section, we explore the minimax optimal MISE rates associated
with several non-standard smoothness classes motivated in Zhu and
Mirzaei (2021). The results in this section rely on the bounds $\overline{B}_{1}\left(\delta\right)$,
$\underline{B}_{1}\left(\delta\right)$ and $\underline{B}_{2}$ in
Lemma \ref{lm:entropy-generalized-polynomial}, as well as the bounds
in Lemmas \ref{lm:entropy-generalized-H=0000F6lder} and \ref{lm:entropy-ellipsoid-subclass}.

\begin{theorem}\label{thm:MISE-gen-sub-ellipsoid}\textit{ Suppose
Assumption \ref{Assumption 1} holds. For the minimax lower bounds,
we further assume that the density $p(x)$ is bounded away from zero;
that is, $p(x)\geq c>0$ for some universal constant $c$.}

\textit{(i) Let $\hat{f}$ be (\ref{eq:KRR estimator}) based on the
kernel function $\mathcal{K}\left(\cdot,\cdot\right)$ associated
with $\mathcal{H}_{\gamma+1}$ in (\ref{eq:Ellipsoid-1}) such that
$\mu_{m}=\left(cm\right)^{-2\left(\gamma+1\right)}$ for a positive
constant $c\asymp1$ independent of $\gamma$ and $R_{\gamma+1}$,
where $\hat{\pi}$ is given by (\ref{eq:optimal weight}) with $\lambda\asymp R_{\gamma+1}^{\frac{-4\left(\gamma+1\right)}{2\left(\gamma+1\right)+1}}\left(\frac{1}{n}\right)^{\frac{2\left(\gamma+1\right)}{2\gamma+3}}$.
Suppose $\mathcal{K}$ is continuous, positive semidefinite, and satisfies
$\mathcal{K}\left(x,x^{'}\right)\precsim1$ for all $x,\,x^{'}\in\left[0,\,1\right]$.
If $R_{\gamma+1}\succsim1$, we have 
\begin{eqnarray*}
\inf_{\tilde{f}}\sup_{f\in\mathcal{H}_{\gamma+1}}\mathbb{E}\left(\left|\tilde{f}-f\right|_{2,\mathbb{P}}^{2}\right) & \succsim & r^{2},\\
\sup_{f\in\mathcal{H}_{\gamma+1}}\mathbb{E}\left(\left|\hat{f}-f\right|_{2,\mathbb{P}}^{2}\right) & \precsim & r^{2}+\exp\left\{ -cnr^{2}\right\} ,
\end{eqnarray*}
where $r^{2}=R_{\gamma+1}^{\frac{2}{2\left(\gamma+1\right)+1}}\left(\frac{\sigma^{2}}{n}\right)^{\frac{2\left(\gamma+1\right)}{2\left(\gamma+1\right)+1}}$. }

\textit{(ii) Let $\hat{f}$ be (\ref{eq:least squares-1}) with $\mathcal{F}=\mathcal{U}_{\gamma+1,2}$.
If $R_{0}\succsim1$, we have 
\begin{eqnarray*}
\inf_{\tilde{f}}\sup_{f\in\mathcal{U}_{\gamma+1,2}}\mathbb{E}\left(\left|\tilde{f}-f\right|_{2,\mathbb{P}}^{2}\right) & \succsim & r^{2},\\
\sup_{f\in\mathcal{U}_{\gamma+1,2}}\mathbb{E}\left(\left|\hat{f}-f\right|_{2,\mathbb{P}}^{2}\right) & \precsim & r^{2}+\exp\left\{ -cnr^{2}\right\} ,
\end{eqnarray*}
where $r^{2}=\left(R^{*}\right)^{\frac{2}{2\left(\gamma+1\right)+1}}\left(\frac{\sigma^{2}}{n}\right)^{\frac{2\left(\gamma+1\right)}{2\left(\gamma+1\right)+1}}$
and }$R^{*}=\left(\max_{k\in\left\{ 1,...,\gamma+1\right\} }\frac{R_{k}}{\left(k-1\right)!}\right)\vee1$.

\end{theorem}

The proof for Theorem \ref{thm:MISE-gen-sub-ellipsoid} is given in
Appendix \ref{subsec:Proof-for-MISE-gen-sub-ellipsoid}.

Theorems \ref{thm:MISE-lower-standard}, \ref{thm:MISE-upper-sobolev-std}
and \ref{thm:MISE-upper-holder-std} suggest that the blessing of
exploiting higher degree of smoothness arises when $\frac{n}{\sigma^{2}}>\left(\gamma+1\right)^{2\gamma+3}$
(here, $\sigma\asymp1$), which clearly includes the case of $\gamma$
being finite and $n$ tending to $\infty$. In these cases, the minimax
optimal rate for the MISE is $\left(\frac{\sigma^{2}}{n}\right)^{\frac{2\gamma+2}{2\gamma+3}}$,
which decreases in $\gamma$. Theorem \ref{thm:MISE-gen-sub-ellipsoid}
suggests that the blessing of exploiting higher degree of smoothness
may also arise when a smoothness class is imposed with the restrictions
that $f^{(k)}(0)=0$ for all $k=0,...,\gamma$. Note that if $R_{\gamma+1}^{\frac{2}{2\gamma+3}}\asymp1$
in the case of $\mathcal{H}_{\gamma+1}$ and $\left(R^{*}\right)^{\frac{2}{2\gamma+3}}\asymp1$
in the case of $\mathcal{U}_{\gamma+1,2}$, the minimax optimal rate
for the MISE is $\left(\frac{\sigma^{2}}{n}\right)^{\frac{2\gamma+2}{2\gamma+3}}$. 

A somewhat counter-intuitive finding from Theorem \ref{thm:MISE-gen-sub-ellipsoid}
is that the parameters $R_{\gamma+1}$ and $R^{*}$ only scale the
standard minimax optimal MISE rate $\left(\frac{\sigma^{2}}{n}\right)^{\frac{2\gamma+2}{2\gamma+3}}$
by $R^{\frac{2}{2\gamma+3}}$ instead of $R$ (where $R=R_{\gamma+1}$
in the case of \textit{$\mathcal{H}_{\gamma+1}$}, and $R=R^{*}$
in the case of \textit{$\mathcal{U}_{\gamma+1,2}$}). Because of the
different forms $R_{\gamma+1}$ and $R^{*}$ take, the optimal rates
can differ between $\mathcal{H}_{\gamma+1}$ and $\mathcal{U}_{\gamma+1,2}$.
For example, when $R_{k}=\left(k-1\right)!$ for all $k=0,...,\gamma+1$,
$R^{*}=1$ and $r^{2}=\left(\frac{\sigma^{2}}{n}\right)^{\frac{2\left(\gamma+1\right)}{2\left(\gamma+1\right)+1}}$
in Theorem \ref{thm:MISE-gen-sub-ellipsoid}(ii). Meanwhile, when
$R_{\gamma+1}=\gamma!$, $r^{2}\asymp\gamma\left(\frac{\sigma^{2}}{n}\right)^{\frac{2\left(\gamma+1\right)}{2\left(\gamma+1\right)+1}}$
in Theorem \ref{thm:MISE-gen-sub-ellipsoid}(i). Note that this difference
cannot be revealed by minimax optimal rates derived based on the metric
entropy $\delta^{\frac{-1}{\gamma+1}}$ (which would simply yield
$\left(\frac{\sigma^{2}}{n}\right)^{\frac{2\left(\gamma+1\right)}{2\left(\gamma+1\right)+1}}$). 

The next two theorems explore the minimax optimal MISE rates for cases
motivated in Zhu and Mirzaei (2021). There are many interesting results
that can be explored using the bounds in Section \ref{sec:Covering-and-packing}.
We focus on the Hölder classes which reveal an interesting contrast
coming from the polynomial subclass when $R_{k}$ is increased from
$\overline{C}\left(k-1\right)!$ to $\overline{C}k!$. 

\begin{theorem}\label{thm:MISE_(k-1)!}

\textit{Suppose Assumption \ref{Assumption 1} holds, $R_{0}=\overline{C}\asymp1$
and $R_{k}$ can be any value in $\left[\overline{C},\,\overline{C}\left(k-1\right)!\right]$
for $k=1,...,\gamma+1$. Let $R^{\dagger}:=1\vee\sum_{k=0}^{\gamma}\frac{R_{k}}{k!}$.
For the minimax lower bounds, we further assume that the density $p(x)$
is bounded away from zero; that is, $p(x)\geq c>0$ for some universal
constant $c$.}

\textit{(i) If 
\begin{equation}
\frac{n}{\sigma^{2}}>\left(\gamma+1\right)^{2\gamma+3},\label{eq:9-2}
\end{equation}
then we have 
\begin{eqnarray*}
\inf_{\tilde{f}}\sup_{f\in\mathcal{U}_{\gamma+1}}\mathbb{E}\left(\left|\tilde{f}-f\right|_{2,\mathbb{P}}^{2}\right) & \geq & \underline{c}_{0}r_{1}^{2},\\
\sup_{f\in\mathcal{U}_{\gamma+1}}\mathbb{E}\left(\left|\hat{f}-f\right|_{2,\mathbb{P}}^{2}\right) & \leq & \overline{c}\left[r_{1}^{2}+\exp\left\{ -cnr_{1}^{2}\right\} \right],
\end{eqnarray*}
for some positive universal constants $\underline{c}_{0}\in(0,\,1]$,
$\overline{c}\in\left(1,\,\infty\right)$ and $c\asymp1$ (all independent
of $n$ and $\gamma$ and bounded away from zero and from above),
where $r_{1}^{2}=\left(\frac{\sigma^{2}}{n}\right)^{\frac{2\left(\gamma+1\right)}{2\left(\gamma+1\right)+1}}$.}

\textit{(ii) If 
\[
\frac{n}{\sigma^{2}}\leq\left(\gamma+1\right)^{2\gamma+3},
\]
we let $\gamma^{*}\in\left\{ 0,...,\gamma\right\} $ be the smallest
integer such that $\frac{n}{\sigma^{2}}\leq\left(\gamma^{*}+1\right)^{2\gamma^{*}+3}$.
Then we have 
\[
\inf_{\tilde{f}}\sup_{f\in\mathcal{U}_{\gamma+1}}\mathbb{E}\left(\left|\tilde{f}-f\right|_{2,\mathbb{P}}^{2}\right)\geq\underline{c}\frac{\sigma^{2}\left(\gamma^{*}+1\right)}{n}
\]
for some universal constant $\underline{c}\in(0,\,1]$ independent
of $n$ and $\gamma$ such that $\underline{c}\leq\frac{\underline{c}_{0}}{2}$,
and bounded away from zero. }

\textit{(iii) If 
\begin{equation}
\frac{n}{\sigma^{2}}\leq\left(\gamma+1\right)^{2\gamma+3},\label{eq:sample size-upper-Theorem5.2}
\end{equation}
let $\hat{f}$ be (\ref{eq:least squares-1}) with $\mathcal{F}=\mathcal{U}_{\gamma^{*}+1}$.
Then we have 
\begin{equation}
\sup_{f\in\mathcal{U}_{\gamma+1}}\mathbb{E}\left(\left|\hat{f}-f\right|_{2,\mathbb{P}}^{2}\right)\leq\overline{c}\left[r_{2}^{2}+\exp\left\{ -cnr_{2}^{2}\right\} \right]\label{eq:34}
\end{equation}
where $r_{2}^{2}=\frac{\sigma^{2}\left(\gamma^{*}+1\right)}{n}$. }

\end{theorem}

The proof for Theorem \ref{thm:MISE_(k-1)!} is given in Appendix
\ref{subsec:Proof-for-MISE_(k-1)!}.\textbf{ }

Theorem \ref{thm:MISE_(k-1)!} is another example (besides Theorems
\ref{thm:MISE-lower-standard}, \ref{thm:MISE-upper-sobolev-std}
and \ref{thm:MISE-upper-holder-std}) that illustrates the importance
of Lemma \ref{lm:entropy-generalized-polynomial}. In particular,
Zhu and Mirzaei (2021) applies the counting argument in Kolmogorov
and Tikhomirov (1959) to derive an upper bound for the covering number
of the polynomial subclass under $R_{k}=\left(k-1\right)!$. If this
result is used to derive an upper bound for the MISE, we would have
obtained (\ref{eq:34}) with $r_{1}^{2}=\frac{\sigma^{2}\left(\gamma^{*}+1\right)^{2}\log\left(\gamma^{*}\vee2\right)}{n}$.
With the new bound $\overline{B}_{1}\left(\delta\right)$ developed
in our Lemma \ref{lm:entropy-generalized-polynomial}, the upper bound
for the MISE shown in Theorem \ref{thm:MISE_(k-1)!} improves $\frac{\sigma^{2}\left(\gamma^{*}+1\right)^{2}\log\left(\gamma^{*}\vee2\right)}{n}$
by a factor of $\left(\gamma^{*}+1\right)\log\left(\gamma^{*}\vee2\right)$.
For the lower bound on the covering number of the Hölder class under
$R_{k}=\left(k-1\right)!$, Zhu and Mirzaei (2021) simply takes the
lower bound $\delta^{\frac{-1}{\gamma+1}}$ from Kolmogorov and Tikhomirov
(1959). As discussed in Section \ref{sec:Covering-and-packing}, this
result ignores the polynomial subclass and is not sharp unless $\delta$
is small enough.

\begin{theorem}\label{thm:MISE_k!} \textit{Suppose Assumption \ref{Assumption 1}
holds and $R_{k}=\overline{C}k!$ for $k=0,...,\gamma+1$, where $\overline{C}\asymp1$.
For the minimax lower bounds, we further assume that the density $p(x)$
is bounded away from zero; that is, $p(x)\geq c>0$ for some universal
constant $c$.}

\textit{(i) If 
\begin{equation}
\frac{n}{\sigma^{2}}>\left(\left(\gamma+1\right)\log\left(\gamma\vee2\right)\right)^{2\gamma+3},\label{eq:sample size-lower-Theorem5.3-1}
\end{equation}
then we have 
\begin{eqnarray*}
\inf_{\tilde{f}}\sup_{f\in\mathcal{U}_{\gamma+1}}\mathbb{E}\left(\left|\tilde{f}-f\right|_{2,\mathbb{P}}^{2}\right) & \geq & \underline{c}_{0}r_{1}^{2},\\
\sup_{f\in\mathcal{U}_{\gamma+1}}\mathbb{E}\left(\left|\hat{f}-f\right|_{2,\mathbb{P}}^{2}\right) & \leq & \overline{c}\left[r_{1}^{2}+\exp\left\{ -cnr_{1}^{2}\right\} \right],
\end{eqnarray*}
for some positive universal constants $\underline{c}_{0}\in(0,\,1]$,
$\overline{c}\in\left(1,\,\infty\right)$ and $c\asymp1$ (all independent
of $n$ and $\gamma$ and bounded away from zero and from above),
where $r_{1}^{2}=\left(\frac{\sigma^{2}}{n}\right)^{\frac{2\left(\gamma+1\right)}{2\left(\gamma+1\right)+1}}$.}

\textit{(ii) If 
\[
\frac{n}{\sigma^{2}}\leq\left(\left(\gamma+1\right)\log\left(\gamma\vee2\right)\right)^{2\gamma+3},
\]
we let $\gamma^{*}\in\left\{ 0,...,\gamma\right\} $ be the smallest
integer such that $\frac{n}{\sigma^{2}}\leq\left(\left(\gamma^{*}+1\right)\log\left(\gamma^{*}\vee2\right)\right)^{2\gamma^{*}+3}$.
Then we have 
\[
\inf_{\tilde{f}}\sup_{f\in\mathcal{U}_{\gamma+1}}\mathbb{E}\left(\left|\tilde{f}-f\right|_{2,\mathbb{P}}^{2}\right)\geq\underline{c}\frac{\sigma^{2}\left(\gamma^{*}+1\right)\log\left(\gamma^{*}\vee2\right)}{n}
\]
for some universal constant $\underline{c}\in(0,\,1]$ independent
of $n$ and $\gamma$ such that $\underline{c}\leq\frac{\underline{c}_{0}}{3}$,
and bounded away from zero. }

\textit{(iii) If 
\begin{equation}
\frac{n}{\sigma^{2}}\leq\left(\left(\gamma+1\right)\log\left(\gamma\vee2\right)\right)^{2\gamma+3},\label{eq:sample size-upper-Theorem5.3}
\end{equation}
let $\hat{f}$ be (\ref{eq:least squares-1}) with $\mathcal{F}=\mathcal{U}_{\gamma^{*}+1}$.
Then we have 
\[
\sup_{f\in\mathcal{U}_{\gamma+1}}\mathbb{E}\left(\left|\hat{f}-f\right|_{2,\mathbb{P}}^{2}\right)\leq\overline{c}\left[r_{2}^{2}+\exp\left\{ -cnr_{2}^{2}\right\} \right]
\]
where $r_{2}^{2}=\frac{\sigma^{2}\left(\gamma^{*}+1\right)\log\left(\gamma^{*}\vee2\right)}{n}$. }

\end{theorem} 

The proof for Theorem \ref{thm:MISE_k!} is given in Appendix \ref{subsec:Proof-for-MISE_k!}.\textbf{ }

A couple of interesting facts are revealed by Theorems \ref{thm:MISE_(k-1)!}
and \ref{thm:MISE_k!}. First, the minimax optimal MISE rates are
the same when $R_{k}$ takes any value in $\left[\overline{C},\,\overline{C}\left(k-1\right)!\right]$.
Second, as $R_{k}$ is increased from $\overline{C}\left(k-1\right)!$
to $\overline{C}k!$, the minimax optimal rate is increased from $\frac{\sigma^{2}\left(\gamma^{*}+1\right)}{n}$
to $\frac{\sigma^{2}\left(\gamma^{*}+1\right)\log\left(\gamma^{*}\vee2\right)}{n}$
when the sample size is not large enough. Once $\frac{n}{\sigma^{2}}>\left(\gamma+1\right)^{2\gamma+3}$
in the case of $R_{k}=\overline{C}\left(k-1\right)!$, and $\frac{n}{\sigma^{2}}>\left(\left(\gamma+1\right)\log\left(\gamma\vee2\right)\right)^{2\gamma+3}$
in the case of $R_{k}=\overline{C}k!$, the optimal rate becomes $\left(\frac{\sigma^{2}}{n}\right)^{\frac{2\gamma+2}{2\gamma+3}}$. 

The terms $\left(k!\right)_{k=0}^{\gamma}$ in (\ref{eq:expansion})
play a more important role on the size of the polynomial subclass
when $R_{k}$ becomes large enough, which is why $\overline{B}_{1}\left(\delta\right)$
and $\underline{B}_{1}\left(\delta\right)$ in Lemma \ref{lm:entropy-generalized-polynomial}
are very useful for deriving the minimax optimal rate for the MISE
under large $R_{k}$. In particular, in dealing with the polynomial
subclass, Theorems \ref{thm:MISE-lower-standard}, \ref{thm:MISE-upper-sobolev-std},
and \ref{thm:MISE-upper-holder-std} rely on the bounds $\overline{B}_{2}\left(\delta\right)$
and $\underline{B}_{2}$ in Lemma \ref{lm:entropy-generalized-polynomial},
Theorem \ref{thm:MISE_(k-1)!} relies on $\overline{B}_{1}\left(\delta\right)$
and $\underline{B}_{2}$, and Theorem \ref{thm:MISE_k!} relies on
$\overline{B}_{1}\left(\delta\right)$ and $\underline{B}_{1}\left(\delta\right)$.
The use of various bounds in Lemma \ref{lm:entropy-generalized-polynomial}
to control for the polynomial component suggests the intricacy of
deriving the minimax optimal MISE rates in the ``small $n$'' regimes.
For example, the intuition for ``$\left(\gamma^{*}+1\right)\log\left(\gamma^{*}\vee2\right)$''
in Theorem \ref{thm:MISE_k!} may be explained by that the metric
entropy of the polynomial subclass with respect to the $L^{q}-$norm
($q\in\left\{ 2,\,\infty\right\} $) under $R_{k}=\overline{C}k!$
behaves the same way as the metric entropy of an $l_{1}-$ball with
respect to the $l_{1}-$norm. However, if this intuition is applied
to the case where $R_{0}=\overline{C}$ and $R_{k}=\overline{C}\left(k-1\right)!$,
then we would expect to see ``$\left(\gamma^{*}+1\right)\log\left(2\vee\log\left(\gamma^{*}\vee2\right)\right)$''
in the minimax optimal rates. It turns out that the correct order
is ``$\gamma^{*}+1$'' (as in Theorem \ref{thm:MISE_(k-1)!}) rather
than ``$\left(\gamma^{*}+1\right)\log\left(2\vee\log\left(\gamma^{*}\vee2\right)\right)$''.

Tables 3 and 4 summarize the results in this section for easy reference.

\begin{table}
{\footnotesize{}\caption{Minimax optimal MISE bounds of the generalized $\mathcal{H}_{\gamma+1}$
and $\mathcal{U}_{\gamma+1,2}$}
}{\scriptsize{}\medskip{}
}{\scriptsize\par}
\centering{}{\scriptsize{}}%
\begin{tabular}{c|c||c||c||c}
\hline 
\multicolumn{1}{c}{} & \multicolumn{2}{c||}{{\scriptsize{}$\mathcal{H}_{\gamma+1}$}} & \multicolumn{2}{c}{{\scriptsize{}$\mathcal{U}_{\gamma+1,2}$}}\tabularnewline
\hline 
{\scriptsize{}MISE} & \multicolumn{2}{c||}{{\scriptsize{}$\in\left[\underline{c}R_{\gamma+1}^{\frac{2}{2\gamma+3}}\left(\frac{\sigma^{2}}{n}\right)^{\frac{2\left(\gamma+1\right)}{2\gamma+3}},\,\overline{c}R_{\gamma+1}^{\frac{2}{2\gamma+3}}\left(\frac{\sigma^{2}}{n}\right)^{\frac{2\left(\gamma+1\right)}{2\gamma+3}}\right]$}} & \multicolumn{2}{c}{{\scriptsize{}$\in\left[\underline{c}^{'}\left(R^{*}\right)^{\frac{2}{2\gamma+3}}\left(\frac{\sigma^{2}}{n}\right)^{\frac{2\left(\gamma+1\right)}{2\gamma+3}},\,\overline{c}^{'}\left(R^{*}\right)^{\frac{2}{2\gamma+3}}\left(\frac{\sigma^{2}}{n}\right)^{\frac{2\left(\gamma+1\right)}{2\gamma+3}}\right]$}}\tabularnewline
\hline 
\end{tabular}{\scriptsize{}}\\
{\scriptsize{}where: $R^{*}=\left(\max_{k\in\left\{ 1,...,\gamma+1\right\} }\frac{R_{k}}{\left(k-1\right)!}\right)\vee1$;
$0<\underline{c}<\overline{c}<\infty$ and $0<\underline{c}^{'}<\overline{c}^{'}<\infty$
are universal constants independent of $n$, $\gamma$ and $R_{k}$s,
and bounded away from zero and from above.}{\scriptsize\par}
\end{table}

\begin{table}[t]
{\small{}\caption{Minimax optimal MISE bounds of non-standard $\mathcal{U}_{\gamma+1}$}
\medskip{}
}{\small\par}
\begin{centering}
\begin{tabular}{c|c|c}
\hline 
\multicolumn{1}{c}{} & \multicolumn{2}{c}{{\scriptsize{}$R_{0}=\overline{C},\,R_{k}\in\left[\overline{C},\,\overline{C}(k-1)!\right]\:\forall k=1,...,\gamma+1$}}\tabularnewline
\hline 
 & {\scriptsize{}$\frac{n}{\sigma^{2}}\leq\left(\gamma+1\right)^{2\gamma+3}$} & {\scriptsize{}$\frac{n}{\sigma^{2}}>\left(\gamma+1\right)^{2\gamma+3}$}\tabularnewline
{\scriptsize{}MISE} & {\scriptsize{}$\in\left[\underline{c}\frac{\sigma^{2}\left(\gamma^{*}+1\right)}{n},\:c^{'}\frac{\sigma^{2}\left(\gamma^{*}+1\right)}{n}\right]$} & {\scriptsize{}$\in\left[\underline{c}_{0}\left(\frac{\sigma^{2}}{n}\right)^{\frac{2\gamma+2}{2\gamma+3}},\,c^{'}\left(\frac{\sigma^{2}}{n}\right)^{\frac{2\gamma+2}{2\gamma+3}}\right]$}\tabularnewline
\hline 
\multicolumn{1}{c}{} & \multicolumn{2}{c}{{\scriptsize{}$R_{k}=\overline{C}k!\:\forall k=0,...,\gamma+1$}}\tabularnewline
\hline 
 & {\scriptsize{}$\frac{n}{\sigma^{2}}\leq\left(\left(\gamma+1\right)\log\left(\gamma\vee2\right)\right)^{2\gamma+3}$} & {\scriptsize{}$\frac{n}{\sigma^{2}}>\left(\left(\gamma+1\right)\log\left(\gamma\vee2\right)\right)^{2\gamma+3}$}\tabularnewline
{\scriptsize{}MISE} & {\scriptsize{}$\in\left[\underline{c}\frac{\sigma^{2}\left(\gamma^{*}+1\right)\log\left(\gamma^{*}\vee2\right)}{n},\,c^{'}\frac{\sigma^{2}\left(\gamma^{*}+1\right)\log\left(\gamma^{*}\vee2\right)}{n}\right]$} & {\scriptsize{}$\in\left[\underline{c}_{0}\left(\frac{\sigma^{2}}{n}\right)^{\frac{2\gamma+2}{2\gamma+3}},\,c^{'}\left(\frac{\sigma^{2}}{n}\right)^{\frac{2\gamma+2}{2\gamma+3}}\right]$}\tabularnewline
\hline 
\end{tabular}
\par\end{centering}
\centering{}{\scriptsize{}where $\underline{c},\,\underline{c}_{0},\,c^{'}\asymp1$
are positive universal constants independent of $n$ and $\gamma$;
these constants can vary from the first class to the second class;
$\underline{c}\leq\frac{\underline{c}_{0}}{2}$ in the first class
and $\underline{c}\leq\frac{\underline{c}_{0}}{3}$ in the second
class; $\gamma^{*}$ is the smallest integer in $\left\{ 0,...,\gamma\right\} $
such that $\frac{n}{\sigma^{2}}\leq\left(\gamma^{*}+1\right)^{2\gamma^{*}+3}$
in the first class (respectively, $\frac{n}{\sigma^{2}}\leq\left(\left(\gamma^{*}+1\right)\log\left(\gamma^{*}\vee2\right)\right)^{2\gamma^{*}+3}$
in the second class).}{\scriptsize\par}
\end{table}

\section{Some insights about multivariate smooth functions \label{sec:multi-dim-extensions}}

The extension of our analysis to $d-$variate smooth functions is
a lot more complex, because of an additional interplay between the
smoothness parameter $\gamma$ and the dimension $d$. We provide
below some partial results about the higher dimensional generalized
Hölder class. 

Let $p=\left(p_{j}\right)_{j=1}^{d}$ and $P=\sum_{j=1}^{d}p_{j}$
where $p_{j}$s are non-negative integers; $x=\left(x_{j}\right)_{j=1}^{d}$
and $x^{p}=\prod_{j=1}^{d}x_{j}^{p_{j}}$. Write $D^{p}f\left(x\right)=\partial^{P}f/\partial x_{1}^{p_{1}}\cdots\partial x_{d}^{p_{d}}$. 

For a non-negative integer $\gamma$, let the \textit{generalized}
Hölder class $\mathcal{U}_{\gamma+1}\left(\left(R_{k}\right)_{k=0}^{\gamma+1},\,\left[-1,\,1\right]^{d}\right)$
be the class of functions such that any function $f\in\mathcal{U}_{\gamma+1}\left(\left(R_{k}\right)_{k=0}^{\gamma+1},\,\left[-1,\,1\right]^{d}\right)$
satisfies: (1) $f$ is continuous on $\left[-1,\,1\right]^{d}$, and
all partial derivatives of $f$ exist for all $p$ with $P:=\sum_{k=1}^{d}p_{k}\leq\gamma$;
(2) $\left|D^{p}f\left(x\right)\right|\leq R_{k}$ for all $p$ with
$P=k$ ($k=0,...,\gamma$) and $x\in\left[-1,\,1\right]^{d}$, where
$D^{0}f\left(x\right)=f\left(x\right)$; (3) $\left|D^{p}f(x)-D^{p}f(x^{'})\right|\leq R_{\gamma+1}\left|x-x^{'}\right|_{\infty}$
for all $p$ with $P=\gamma$ and $x,\,x^{'}\in\left[-1,\,1\right]^{d}$. 

Given any function $f$ in the $d-$variate Hölder class, we have
\[
f(x)=\sum_{k=0}^{\gamma}\sum_{p:P=k}\frac{x^{p}D^{p}f\left(0\right)}{k!}+\sum_{p:P=\gamma}\frac{x^{p}D^{p}f\left(z\right)}{\gamma!}-\sum_{p:P=\gamma}\frac{x^{p}D^{p}f\left(0\right)}{\gamma!}
\]
for some intermediate value $z$. Similar to Section \ref{subsec:Classes-of-smooth},
we have the following relationships: 
\begin{eqnarray*}
 & \mathcal{U}_{\gamma+1,1}^{d}\subseteq\mathcal{U}_{\gamma+1}^{d}, & \mathcal{U}_{\gamma+1,2}^{d}\subseteq\mathcal{U}_{\gamma+1}^{d},\\
\mathcal{U}_{\gamma+1}^{d}\subseteq & \mathcal{U}_{\gamma+1,1}^{d}+\mathcal{U}_{\gamma+1,2}^{d} & :=\left\{ f_{1}+f_{2}:f_{1}\in\mathcal{U}_{\gamma+1,1}^{d},\,f_{2}\in\mathcal{U}_{\gamma+1,2}^{d}\right\} 
\end{eqnarray*}
where 
\[
\mathcal{U}_{\gamma+1,1}^{d}=\left\{ f=\sum_{k=0}^{\gamma}\sum_{p:P=k}x^{p}\theta_{(p,k)}:\{\theta_{(p,k)}\}_{(p,k)}\ensuremath{\in}\ensuremath{\mathcal{P}_{\Gamma}},\,x\in\left[-1,\,1\right]^{d}\right\} 
\]
with the $\Gamma:=\sum_{k=0}^{\gamma}\left(\begin{array}{c}
d+k-1\\
d-1
\end{array}\right)-$dimensional polyhedron
\[
\ensuremath{\mathcal{P}_{\Gamma}=\left\{ \{\theta_{(p,k)}\}_{(p,k)}\ensuremath{\ensuremath{\in}}\mathbb{R}^{\Gamma}:\textrm{for any given }k\in\left\{ 0,...,\gamma\right\} ,\,\theta_{(p,k)}\in\left[\frac{-R_{k}}{k!},\,\frac{R_{k}}{k!}\right]\textrm{ for all }p\textrm{ with }P\leq k\right\} }
\]
where $\theta=\{\theta_{(p,k)}\}_{(p,k)}$ denotes the collection
of $\theta_{(p,k)}$ over all $(p,k)$ configurations. The Hölder
subclass $\mathcal{U}_{\gamma+1,2}^{d}$ is the class of functions
such that any function $f\in\mathcal{U}_{\gamma+1,2}^{d}$ satisfies:
$f\in\mathcal{U}_{\gamma+1}^{d}$ such that $D^{p}f\left(0\right)=0$
for all $p$ with $P\leq k$, $k=0,...,\gamma$.

\begin{lemma}\label{lm:entropy_multi_dim_Holder_sub}\textit{ Let
$R^{*}=\left(\max_{k\in\left\{ 1,...,\gamma+1\right\} }\frac{D_{k-1}^{*}R_{k}}{\left(k-1\right)!}\vee1\right)$.
We have 
\begin{eqnarray}
\log N_{2,\mathbb{P}}\left(\delta,\,\mathcal{U}_{\gamma+1,2}^{d}\right) & \leq & \log N_{\infty}\left(\delta,\,\mathcal{U}_{\gamma+1,2}^{d}\right)\precsim d^{d}R^{*\frac{d}{\gamma+1}}\delta^{\frac{-d}{\gamma+1}},\label{eq:upper2-1-1}\\
\log M_{\infty}\left(\delta,\,\mathcal{U}_{\gamma+1,2}^{d}\right) & \succsim & \log M_{2}\left(\delta,\,\mathcal{U}_{\gamma+1,2}^{d}\right)\succsim d^{d}R^{*\frac{d}{\gamma+1}}\delta^{\frac{-d}{\gamma+1}}.\label{eq:upper2-1b-1}
\end{eqnarray}
}\textbf{Remark}. With Lemma \ref{lm:entropy_multi_dim_Holder_sub},
we can easily establish the minimax optimal MSE rate for $\mathcal{U}_{\gamma+1,2}^{d}$,
using arguments almost identical to those for Theorem \ref{thm:MISE-gen-sub-ellipsoid}.

\end{lemma}

The proof for Lemma \ref{lm:entropy_multi_dim_Holder_sub} is given
in Appendix \ref{subsec:Proof-for-entropy-multi-dim-holder-sub}.\textbf{ }

\begin{lemma}\label{lm:entropy_multi_dim_poly}

\textit{If $\delta$ is small enough such that $\min_{k\in\left\{ 0,...,\gamma\right\} }\log\frac{4\left(\gamma+1\right)D_{k}^{*}R_{k}}{\delta k!}\geq0$,
we have 
\[
\log N_{2,\mathbb{P}}\left(\delta,\,\mathcal{U}_{\gamma+1,1}^{d}\right)\leq\log N_{\infty}\left(\delta,\,\mathcal{U}_{\gamma+1,1}^{d}\right)\leq\sum_{k=0}^{\gamma}D_{k}^{*}\log\frac{4\left(\gamma+1\right)D_{k}^{*}R_{k}}{\delta k!}
\]
where $D_{k}^{*}=\left(\begin{array}{c}
d+k-1\\
d-1
\end{array}\right)$; if $\delta$ is large enough such that $\min_{k\in\left\{ 0,...,\gamma\right\} }\log\frac{4\left(\gamma+1\right)D_{k}^{*}R_{k}}{\delta k!}<0$,
we have
\begin{equation}
\log N_{2,\mathbb{P}}\left(\delta,\,\mathcal{U}_{\gamma+1,1}^{d}\right)\leq\log N_{\infty}\left(\delta,\,\mathcal{U}_{\gamma+1,1}^{d}\right)\precsim\left(\sum_{k=0}^{\gamma}D_{k}^{*}\right)\log\frac{1}{\delta}+\sum_{k=0}^{\gamma}D_{k}^{*}\log R_{k}.\label{eq:kolmogorov_d}
\end{equation}
}

\end{lemma}

\begin{remark} The bound (\ref{eq:kolmogorov_d}) holds for all $\delta\in(0,\,1)$
(not just $\delta$ such that $\min_{k}\log\frac{4\left(\gamma+1\right)D_{k}^{*}R_{k}}{\delta k!}<0$)
but is far from being tight when $\min_{k\in\left\{ 0,...,\gamma\right\} }\log\frac{4\left(\gamma+1\right)D_{k}^{*}R_{k}}{\delta k!}\geq0$.

\end{remark}

\begin{remark} A simple upper bound on $\sum_{k=0}^{\gamma}D_{k}^{*}$
is $\sum_{k=1}^{\gamma}d^{k}\asymp d^{\gamma}$. Let us show a lower
bound on $\sum_{k=0}^{\gamma}D_{k}^{*}$ for the case of $\gamma\geq2d^{2}$
to illustrate how large $\frac{\sigma^{2}}{n}\sum_{k=0}^{\gamma}D_{k}^{*}$
can be. We can write $D_{k}^{*}=\frac{\left(k+d-1\right)!}{\left(d-1\right)!k!}=\prod_{j=1}^{d-1}\frac{k+j}{j}$.
Because $\gamma\geq2d^{2}$, we have 
\begin{eqnarray*}
\sum_{k=0}^{\gamma}D_{k}^{*} & = & \left(\sum_{k=0}^{\gamma}\prod_{j=1}^{d-1}\frac{k+j}{j}\right)\geq\left(\sum_{k=d^{2}}^{\gamma}\prod_{j=1}^{d-1}\frac{k+j}{j}\right)\\
 & \geq & \left(d^{2}\prod_{j=1}^{d-1}\left(\frac{d^{2}}{j}+1\right)\right)\geq\left(d^{2}\left(\frac{d^{2}}{d}+1\right)^{d-1}\right)\geq d^{d+1}.
\end{eqnarray*}

\end{remark}

The proof for Lemma \ref{lm:entropy_multi_dim_poly} is given in Appendix
\ref{subsec:Proof-for-entropy_multi_poly}.\textbf{ }In theory, our
arguments for $\underline{B}_{1}\left(\delta\right)$ in Lemma \ref{lm:entropy-generalized-polynomial}
can be extended for analyzing the lower bound for $\log M_{2}\left(\delta,\,\mathcal{U}_{\gamma+1,1}^{d}\right)$.
However, this extension is very intensive. Arguments similar to those
for $\underline{B}_{2}$ in Lemma \ref{lm:entropy-generalized-polynomial}
will not lead to a useful bound for $\mathcal{U}_{\gamma+1,1}^{d}$.
Despite the lack of lower bounds, we can still gain some insights
from Lemma \ref{lm:entropy_multi_dim_poly}, as it implies
\[
\mathbb{E}\left(\left|\hat{f}-f\right|_{2,\mathbb{P}}^{2}\right)\precsim r^{2}+\exp\left\{ -cnr^{2}\right\} 
\]
where $r^{2}=\frac{\sigma^{2}}{n}\sum_{k=0}^{\gamma}D_{k}^{*}$ and
$\hat{f}$ is the estimator in (\ref{eq:least squares-1}) with $\mathcal{F}=\mathcal{U}_{\gamma+1,1}^{d}$.
The quantity $\sum_{k=0}^{\gamma}D_{k}^{*}$ is the higher dimensional
analogue of $\gamma+1$ and arises from the fact that a function in
$\mathcal{U}_{\gamma+1,1}^{d}$ has $D_{k}^{*}$ distinct $k$th partial
derivatives. 

Suppose $R_{k}=1$ for all $k=0,...,\gamma+1$. If $d$ is small relative
to $\gamma$ and $n$, Lemma \ref{lm:entropy_multi_dim_Holder_sub}
implies that the minimax optimal rate concerning $\mathcal{U}_{\gamma+1,2}^{d}$
is roughly $\left(\frac{\sigma^{2}}{n}\right)^{\frac{2\gamma+2}{2\gamma+2+d}}$,
the classical rate for $\mathcal{U}_{\gamma+1}^{d}$ derived under
the regime where $\gamma$ and $d$ are finite but $n\rightarrow\infty$.
Observe that $\frac{\sigma^{2}}{n}\sum_{k=0}^{\gamma}D_{k}^{*}\geq\left(\frac{\sigma^{2}}{n}\right)^{\frac{2\gamma+2}{2\gamma+2+d}}$
whenever $\frac{n}{\sigma^{2}}\leq\left(\sum_{k=0}^{\gamma}D_{k}^{*}\right)^{\frac{2\gamma+2+d}{d}}$,
and $\frac{\sigma^{2}}{n}\sum_{k=0}^{\gamma}D_{k}^{*}<\left(\frac{\sigma^{2}}{n}\right)^{\frac{2\gamma+2}{2\gamma+2+d}}$
whenever $\frac{n}{\sigma^{2}}>\left(\sum_{k=0}^{\gamma}D_{k}^{*}\right)^{\frac{2\gamma+2+d}{d}}$.
Therefore, the classical asymptotic minimax rate $\left(\frac{\sigma^{2}}{n}\right)^{\frac{2\gamma+2}{2\gamma+2+d}}$
could be an underestimate of the MISE when $n$ is not large enough
and the optimal degree of smoothness to exploit may not be the maximum
smoothness degree $\gamma+1$, but rather depend on both $n$ and
$d$. 

\section{Proofs for Section \ref{sec:Covering-and-packing}\label{sec:proofs_entropy} }

\subsection{Proof for Lemma \ref{lm:entropy-generalized-polynomial}\label{subsec:Proof-for-entropy-poly}}

\textbf{The upper bound}. Recall the definition of $\mathcal{U}_{\gamma+1,1}$:
\[
\mathcal{U}_{\gamma+1,1}=\left\{ f=\sum_{k=0}^{\gamma}\theta_{k}x^{k}:\textrm{\ensuremath{\left(\theta_{k}\right)_{k=0}^{\gamma}}\ensuremath{\ensuremath{\in}}\ensuremath{\mathcal{P}_{\gamma}}},\,x\in\left[-1,\,1\right]\right\} 
\]
with the $\left(\gamma+1\right)-$dimensional polyhedron 
\[
\ensuremath{\mathcal{P}_{\gamma}=\left\{ \ensuremath{\left(\theta_{k}\right)_{k=0}^{\gamma}}\ensuremath{\ensuremath{\in}}\mathbb{R}^{\gamma+1}:\theta_{k}\in\left[\frac{-R_{k}}{k!},\,\frac{R_{k}}{k!}\right]\right\} }
\]
where $R_{k}$ is allowed to depend on $k\in\left\{ 0,...,\gamma\right\} $
only. We first derive an upper bound for $N_{\infty}\left(\delta,\,\mathcal{U}_{\gamma+1,1}\right)$.
Because the weighted $L^{2}\left(\mathbb{P}\right)-$norm is no greater
than the sup norm and a smallest $\delta-$cover of $\mathcal{U}_{\gamma+1,1}$
with respect to the $\left|\cdot\right|_{\infty}$ norm also covers
$\mathcal{U}_{\gamma+1,1}$ with respect to the $\left|\cdot\right|_{2,\mathbb{P}}$
norm, we have 
\[
N_{2,\mathbb{P}}\left(\delta,\,\mathcal{U}_{\gamma+1,1}\right)\leq N_{\infty}\left(\delta,\,\mathcal{U}_{\gamma+1,1}\right).
\]

To bound $\log N_{\infty}\left(\delta,\,\mathcal{U}_{\gamma+1,1}\right)$
from above, note that for $f,\,f^{'}\in\mathcal{U}_{\gamma+1,1}$,
we have 
\[
\left|f-f^{'}\right|_{\infty}\leq\sum_{k=0}^{\gamma}\left|\theta_{k}-\theta_{k}^{'}\right|
\]
where $f^{'}=\sum_{k=0}^{\gamma}\theta_{k}^{'}x^{k}$ such that $\theta^{'}\in\mathcal{P}_{\gamma}$.
Therefore, the problem is reduced to bounding $N_{1}\left(\delta,\,\mathcal{P}_{\gamma}\right)$.

Consider $\left(a_{k}\right)_{k=0}^{\gamma}$ such that $a_{k}>0$
for every $k=0,...,\gamma$ and $\sum_{k=0}^{\gamma}a_{k}=1$. To
cover $\mathcal{P}_{\gamma}$ within $\delta-$precision, we find
a smallest $a_{k}\delta-$cover of $\left[\frac{-R_{k}}{k!},\,\frac{R_{k}}{k!}\right]$
for each $k=0,...,\gamma$, $\left\{ \theta_{k}^{1},...,\theta_{k}^{N_{k}}\right\} $,
such that for any $\theta\in\mathcal{P}_{\gamma}$, there exists some
$i_{k}\in\left\{ 1,...,N_{k}\right\} $ with 
\[
\sum_{k=0}^{\gamma}\left|\theta_{k}-\theta_{k}^{i_{k}}\right|\leq\delta.
\]
As a consequence, we have 
\begin{equation}
\log N_{1}\left(\delta,\,\mathcal{P}_{\gamma}\right)\leq\sum_{k=0}^{\gamma}\log\frac{4R_{k}}{a_{k}k!\delta}=-\sum_{k=0}^{\gamma}\log a_{k}+\sum_{k=0}^{\gamma}\log\frac{4R_{k}}{k!\delta}.\label{eq:upper_poly}
\end{equation}
For $\left(a_{k}\right)_{k=0}^{\gamma}$ such that $\sum_{k=0}^{\gamma}a_{k}=1$,
the function 
\[
h(a_{0},...,a_{\gamma}):=-\sum_{k=0}^{\gamma}\log a_{k}=-\log\left(\prod_{k=0}^{\gamma}a_{k}\right)
\]
is minimized at $a_{k}=\frac{1}{\gamma+1}$. Consequently, the minimum
of $\sum_{k=0}^{\gamma}\log\frac{4R_{k}}{a_{k}k!\delta}$ equals $\sum_{k=0}^{\gamma}\log\frac{4\left(\gamma+1\right)R_{k}}{k!\delta}$
and we have 
\[
\log N_{1}\left(\delta,\,\mathcal{P}_{\gamma}\right)\leq\sum_{k=0}^{\gamma}\log\frac{4\left(\gamma+1\right)R_{k}}{k!\delta}.
\]
Therefore, 
\begin{align}
\log N_{2,\mathbb{P}}\left(\delta,\,\mathcal{U}_{\gamma+1,1}\right) & \leq\log N_{\infty}\left(\delta,\,\mathcal{U}_{\gamma+1,1}\right)\leq\sum_{k=0}^{\gamma}\log\frac{4\left(\gamma+1\right)R_{k}}{k!\delta}.\label{eq:upper}
\end{align}
If $\delta$ is large enough such that $\min_{k\in\left\{ 0,...,\gamma\right\} }\log\frac{4\left(\gamma+1\right)R_{k}}{k!\delta}<0$,
we can evoke the counting argument in Kolmogorov and Tikhomirov (1959)
to obtain 
\begin{equation}
\log N_{2,\mathbb{P}}\left(\delta,\,\mathcal{U}_{\gamma+1,1}\right)\leq\log N_{\infty}\left(\delta,\,\mathcal{U}_{\gamma+1,1}\right)\leq\left(\frac{\gamma}{2}+1\right)\log\frac{1}{\delta}+\sum_{k=0}^{\gamma}\log R_{k}.\label{eq:Kolmogorov_upper-1}
\end{equation}
\textbf{The lower bound}. We first derive a lower bound for $M_{2}\left(\delta,\,\mathcal{U}_{\gamma+1,1}\right)$.
Let $\left(\phi_{k}\right)_{k=0}^{\gamma}$ be the Legendre polynomials
on $\left[-1,\,1\right]$. For any function $f\in\mathcal{U}_{\gamma+1,1}$,
we can write 
\begin{equation}
f\left(x\right)=\sum_{k=0}^{\gamma}\tilde{\theta}_{k}\phi_{k}\left(x\right)\label{eq:orthogonal expansion}
\end{equation}
such that 
\begin{equation}
\tilde{\theta}_{k}=\frac{\left(2k+1\right)}{2}\int_{-1}^{1}f(x)\phi_{k}(x)dx.\label{eq:coeff1}
\end{equation}
In Lemma A.1 of Section \ref{subsec:Lemma-A.1-and}, we show that
\[
\tilde{\theta}_{k}=\left(k+\frac{1}{2}\right)\sum_{m=0}^{\left\lfloor \gamma/2\right\rfloor }\frac{f^{(k+2m)}(0)}{2^{k+2m}m!\left(\frac{1}{2}\right)_{k+m+1}}
\]
where $\left(a\right)_{k}=a(a+1)\cdots(a+k-1)$ is known as the Pochhammer
symbol. Recall $\left|f^{(k)}\left(0\right)\right|\leq R_{k}$ for
$k=0,...,\gamma$ and $f^{(k)}\left(0\right)=0$ for $k>\gamma$.
We can rewrite
\[
\mathcal{U}_{\gamma+1,1}=\left\{ f=\sum_{k=0}^{\gamma}\tilde{\theta}_{k}\phi_{k}\left(x\right):\textrm{\ensuremath{\left(\tilde{\theta}_{k}\right)_{k=0}^{\gamma}}\ensuremath{\ensuremath{\in}}\ensuremath{\mathcal{P}_{\gamma}^{L}}},\,x\in\left[-1,\,1\right]\right\} 
\]
with the $\left(\gamma+1\right)-$dimensional polyhedron 
\[
\ensuremath{\mathcal{P}_{\gamma}^{L}=\left\{ \ensuremath{\left(\tilde{\theta}_{k}\right)_{k=0}^{\gamma}}\ensuremath{\ensuremath{\in}}\mathbb{R}^{\gamma+1}:\tilde{\theta}_{k}\in\left[-\overline{R}_{k},\,\overline{R}_{k}\right]\right\} }
\]
where $\overline{R}_{k}:=\sum_{m=0}^{\left\lfloor \gamma/2\right\rfloor }b_{k,m}R_{k+2m}$
and $b_{k,m}=\frac{\left(k+\frac{1}{2}\right)}{2^{k+2m}m!\left(\frac{1}{2}\right)_{k+m+1}}$.
If we can bound $\overline{R}_{k}$ from below by $\underline{R}_{k}$,
then we have 
\begin{equation}
\mathcal{P}_{\gamma}^{L}\supseteq\mathcal{\underline{P}}_{\gamma}^{L}=\left\{ \ensuremath{\left(\tilde{\theta}_{k}\right)_{k=0}^{\gamma}}\ensuremath{\ensuremath{\in}}\mathbb{R}^{\gamma+1}:\tilde{\theta}_{k}\in\left[-\underline{R}_{k},\,\underline{R}_{k}\right]\right\} .\label{eq:12}
\end{equation}
Let us derive $\underline{R}_{k}$. Because $f^{(l)}\left(0\right)=0$
for $l>\gamma$,
\[
\frac{f^{(k+2m)}(0)}{2^{k+2m}m!\left(\frac{1}{2}\right)_{k+m+1}}=0\qquad\textrm{if }k+2m>\gamma.
\]
There are at most $\gamma+1$ terms that are multiplied in the product
$m!\left(\frac{1}{2}\right)_{k+m+1}$. Note that $m\leq\frac{\gamma}{2}\leq\frac{3\gamma}{2}+1$
and 

\begin{eqnarray*}
\left(\frac{1}{2}\right)_{k+m+1} & = & \frac{1}{2}\frac{1+2}{2}\cdots\frac{1+2(k+m)}{2}\\
 & \leq & \frac{2}{2}\frac{2+2}{2}\cdots\frac{2+2(k+m)}{2}\\
 & = & \left(k+m+1\right)!
\end{eqnarray*}
where $k+m+1\leq\frac{3\gamma}{2}+1$. Hence, we have

\[
m!\left(\frac{1}{2}\right)_{k+m+1}\leq m!\left(k+m+1\right)!\leq1\cdot\left(\frac{3\gamma}{2}+1\right)^{\gamma}\leq\left(3\gamma\right)^{\gamma}.
\]
As a result, we have 
\begin{eqnarray}
\overline{R}_{k}=\sum_{m=0}^{\left\lfloor \gamma/2\right\rfloor }b_{k,m}R_{k+2m} & = & \sum_{m=0}^{\left\lfloor \gamma/2\right\rfloor }\frac{\left(k+\frac{1}{2}\right)R_{k+2m}}{2^{k+2m}m!\left(\frac{1}{2}\right)_{k+m+1}}\nonumber \\
 & \geq & \left(k+\frac{1}{2}\right)2^{-\gamma}3^{-\gamma}\gamma^{-\gamma}\sum_{m=0}^{\left\lfloor \gamma/2\right\rfloor }R_{k+2m}\nonumber \\
 & \geq & \frac{9^{-\gamma}\gamma^{-\gamma}}{2}\sum_{m=0}^{\left\lfloor \gamma/2\right\rfloor }R_{k+2m}=:\underline{R}_{k}.\label{eq:fact1}
\end{eqnarray}
Note that for any $f,g\in\mathcal{U}_{\gamma+1,1}$ where $f\left(x\right)=\sum_{k=0}^{\gamma}\tilde{\theta}_{k}\phi_{k}\left(x\right)$
and $g\left(x\right)=\sum_{k=0}^{\gamma}\tilde{\theta}_{k}^{'}\phi_{k}\left(x\right)$,
we have 
\begin{equation}
\left|f-g\right|_{2}^{2}=\sum_{k=0}^{\gamma}\left[\sqrt{\frac{2}{2k+1}}\left(\tilde{\theta}_{k}-\tilde{\theta}_{k}^{'}\right)\right]^{2}.\label{eq:27}
\end{equation}
In view of (\ref{eq:12}) and (\ref{eq:27}), to construct a packing
set of $\mathcal{U}_{\gamma+1,1}$ with $\delta-$separation, we find
a largest $\sqrt{\frac{2k+1}{2}}\frac{\delta}{\sqrt{\gamma+1}}-$packing
of $\left[-\underline{R}_{k},\,\underline{R}_{k}\right]$ for each
$k=0,...,\gamma$, $\left\{ \tilde{\theta}_{k}^{1},...,\tilde{\theta}_{k}^{M_{k}}\right\} $,
such that for any distinct $\theta_{k}^{i_{k}}$ and $\theta_{k}^{j_{k}}$
in the packing sets,
\begin{equation}
\sum_{k=0}^{\gamma}\left[\sqrt{\frac{2}{2k+1}}\left(\tilde{\theta}_{k}^{i_{k}}-\tilde{\theta}_{k}^{j_{k}}\right)\right]^{2}>\delta^{2}.\label{eq:48}
\end{equation}
Therefore,\footnote{In fact, we can make the same statement about an exactly $\delta-$separated
set, i.e., $=\delta^{2}$ rather than $>\delta^{2}$ in (\ref{eq:48}).} 
\begin{equation}
\log M_{2}\left(\delta,\,\mathcal{U}_{\gamma+1,1}\right)\succsim\sum_{k=0}^{\gamma}\log\frac{\sqrt{2\left(\gamma+1\right)}\underline{R}_{k}}{\sqrt{2k+1}\delta}.\label{eq:upper_poly-1}
\end{equation}
Bounds (\ref{eq:upper_poly-1}) and (\ref{eq:fact1}) together give
\begin{equation}
\log M_{2}\left(\delta,\,\mathcal{U}_{\gamma+1,1}\right)\geq\sum_{k=0}^{\gamma}\log\left(9^{-\gamma}\gamma^{-\gamma}\right)+\sum_{k=0}^{\gamma}\log\frac{C\sum_{m=0}^{\left\lfloor \gamma/2\right\rfloor }R_{k+2m}}{\delta}=:\underline{B}_{1}\left(\delta\right)\label{eq:lower}
\end{equation}
for some positive universal constant $C$. 

The following argument gives another useful bound for $\log M_{2}\left(\delta,\,\mathcal{U}_{\gamma+1,1}\right)$. 

Let $\tilde{k}\in\arg\max_{k\in\left\{ 0,...,\gamma\right\} }\frac{R_{k}}{k!}$.
We consider a $3\delta\left[\left(\tilde{k}+1\right)\vee\sum_{k=0}^{\gamma}\frac{R_{k}}{k!}\right]-$grid
of points on $\left[-\frac{R_{\tilde{k}}}{\tilde{k}!},\,\frac{R_{\tilde{k}}}{\tilde{k}!}\right]$
and denote the collection of these points by $\left(\theta_{\tilde{k}}^{*i}\right)_{i=1}^{M_{0}}$
where $M_{0}=\left\lceil \frac{2R_{\tilde{k}}}{3\tilde{k}!\delta\left[\left(\tilde{k}+1\right)\vee\sum_{k=0}^{\gamma}\frac{R_{k}}{k!}\right]}\right\rceil +1$. 

We choose $\delta$ such that $M_{0}\geq2^{\gamma+1}$ and $3\delta\left[\left(\tilde{k}+1\right)\vee\sum_{k=0}^{\gamma}\frac{R_{k}}{k!}\right]\leq\frac{2R_{\tilde{k}}}{\tilde{k}!}$.
Let us fix $\theta_{k}^{*}\in\left[0,\,\frac{R_{k}\delta}{bk!(\tilde{k}+1)}\right]$
for\footnote{The parameter $b\in[1,\infty)$ is chosen according to need later
when we derive the miminax lower bounds.} $k\in\left\{ 0,...,\gamma\right\} \setminus\tilde{k}$ and define
\begin{equation}
f_{\lambda_{i}}^{*}\left(x\right)=\theta_{\tilde{k}}^{*i}x^{\tilde{k}}+\sum_{k\in\left\{ 0,...,\gamma\right\} \setminus\tilde{k}}\lambda_{i,k}\theta_{k}^{*}x^{k},\quad x\in\left[-1,\,1\right]\label{eq:M-packing}
\end{equation}
where $\left(\lambda_{i,k}\right)_{k\in\left\{ 0,...,\gamma\right\} \setminus\tilde{k}}=:\lambda_{i}\in\left\{ 0,\,1\right\} ^{\gamma}$
for all $i=1,...,M_{0}$. For any $\lambda_{i,k},\,\lambda_{j,k}\in\left\{ 0,\,1\right\} ^{\gamma}$
such that $i\neq j$, we have 
\begin{eqnarray*}
\left|f_{\lambda_{i}}^{*}-f_{\lambda_{j}}^{*}\right|_{2} & = & \left[\int_{-1}^{1}\left[\left(\theta_{\tilde{k}}^{*i}-\theta_{\tilde{k}}^{*j}\right)x^{\tilde{k}}+\sum_{k\in\left\{ 0,...,\gamma\right\} \setminus\tilde{k}}\left(1\left\{ \lambda_{i,k}\neq\lambda_{j,k}\right\} \theta_{k}^{*}x^{k}\right)\right]^{2}dx\right]^{\frac{1}{2}}\\
 & \geq & \left[\int_{0}^{1}\left[\left(\theta_{\tilde{k}}^{*i}-\theta_{\tilde{k}}^{*j}\right)x^{\tilde{k}}+\sum_{k\in\left\{ 0,...,\gamma\right\} \setminus\tilde{k}}\left(1\left\{ \lambda_{i,k}\neq\lambda_{j,k}\right\} \theta_{k}^{*}x^{k}\right)\right]^{2}dx\right]^{\frac{1}{2}}\\
 & \geq & \int_{0}^{1}\left|\left(\theta_{\tilde{k}}^{*i}-\theta_{\tilde{k}}^{*j}\right)x^{\tilde{k}}+\sum_{k\in\left\{ 0,...,\gamma\right\} \setminus\tilde{k}}\left(1\left\{ \lambda_{i,k}\neq\lambda_{j,k}\right\} \theta_{k}^{*}x^{k}\right)\right|dx\\
 & \geq & \left[\left(\frac{2\delta}{\tilde{k}+1}\sum_{k=0}^{\gamma}\frac{R_{k}}{k!}\right)\vee\left(2\delta\right)\right]>\delta
\end{eqnarray*}
where the third line follows from the Jensen's inequality and the
concavity of $\sqrt{\cdot}$ on $[0,\,1]$, and the fourth line follows
from the triangle inequality. Hence, we have constructed a $\delta-$packing
set. The cardinality of this packing set is at least $2^{\gamma+1}$.
That is, 
\[
\log M_{2}\left(\delta,\,\mathcal{U}_{\gamma+1,1}\right)\succsim\gamma+1=:\underline{B}_{2}.
\]
Note that the lower bound $\log M_{2}\left(\delta,\,\mathcal{U}_{\gamma+1,1}\right)\succsim\gamma+1$
holds for all $\delta$ such that $\left\lceil \frac{2R_{\tilde{k}}}{3\tilde{k}!\delta\left[\left(\tilde{k}+1\right)\vee\sum_{k=0}^{\gamma}\frac{R_{k}}{k!}\right]}\right\rceil +1\geq2^{\gamma+1}$
and $3\delta\left[\left(\tilde{k}+1\right)\vee\sum_{k=0}^{\gamma}\frac{R_{k}}{k!}\right]\leq\frac{2R_{\tilde{k}}}{\tilde{k}!}$. 

Because $\left|\cdot\right|_{\infty}\geq\frac{1}{2}\left|\cdot\right|_{2}$,
we clearly have 

\[
\log M_{\infty}\left(\delta,\,\mathcal{U}_{\gamma+1,1}\right)\succsim\gamma+1
\]
for all $\delta$ such that $\left\lceil \frac{2R_{\tilde{k}}}{3\tilde{k}!\delta\left[\left(\tilde{k}+1\right)\vee\sum_{k=0}^{\gamma}\frac{R_{k}}{k!}\right]}\right\rceil +1\geq2^{\gamma+1}$
and $3\delta\left[\left(\tilde{k}+1\right)\vee\sum_{k=0}^{\gamma}\frac{R_{k}}{k!}\right]\leq\frac{2R_{\tilde{k}}}{\tilde{k}!}$,\footnote{Note that these two conditions can be reduced to $\left\lceil \frac{2R_{\tilde{k}}}{3\tilde{k}!\delta\left[\left(\tilde{k}+1\right)\vee\sum_{k=0}^{\gamma}\frac{R_{k}}{k!}\right]}\right\rceil +1\geq2^{\gamma+1}$.}
and also
\[
\log M_{\infty}\left(\delta,\,\mathcal{U}_{\gamma+1,1}\right)\succsim\underline{B}_{1}\left(\delta\right).
\]

Finally, if the density function $p(x)$ on $\left[-1,\,1\right]$
is bounded away from zero, i.e., $p(x)\geq c>0$, then 
\[
\log M_{2,\mathbb{P}}\left(\delta,\,\mathcal{U}_{\gamma+1,1}\right)\asymp\log M_{2}\left(\delta,\,\mathcal{U}_{\gamma+1,1}\right)
\]
and therefore we have claim (iii).

\subsection{Lemma A.1 and its proof\label{subsec:Lemma-A.1-and}}

\textbf{Lemma A.1}. \textit{Let $\left\{ \phi_{k}\right\} _{k=1}^{\infty}$
be the Legendre polynomials on $\left[-1,\,1\right]$. For any $f\in\mathcal{U}_{\gamma+1,1}\left[-1,\,1\right]$,
we have $f(x)=\sum_{k=0}^{\gamma}\tilde{\theta}_{k}\phi_{k}(x)$ such
that 
\[
\tilde{\theta}_{k}=\left(k+\frac{1}{2}\right)\sum_{m=0}^{\left\lfloor \gamma/2\right\rfloor }\frac{f^{(k+2m)}(0)}{2^{k+2m}m!\left(\frac{1}{2}\right)_{k+m+1}}
\]
where $\left(a\right)_{k}=a(a+1)\cdots(a+k-1)$ is known as the Pochhammer
symbol.}\\
\textbf{}\\
\textbf{Proof}. To obtain the correct formula for finite sums, we
carefully modify the derivations in Cantero and Iserles (2012) which
concerns infinite sums. The Legendre expansion of $x^{k}$ yields
\begin{equation}
\frac{x^{k}}{k!}=\frac{1}{2^{k}}\sum_{m=0}^{\left\lfloor k/2\right\rfloor }\frac{k-2m+\frac{1}{2}}{m!\left(\frac{1}{2}\right)_{k-m+1}}\phi_{k-2m}(x).\label{eq:legendre_polynomials}
\end{equation}
First, let us consider the case where $\gamma$ is odd. Applying (\ref{eq:legendre_polynomials})
gives
\begin{eqnarray}
f\left(x\right) & = & \sum_{k=0}^{\gamma}\frac{f^{(k)}(0)}{2^{k}}\sum_{m=0}^{\left\lfloor k/2\right\rfloor }\frac{k-2m+\frac{1}{2}}{m!\left(\frac{1}{2}\right)_{k-m+1}}\phi_{k-2m}(x)\nonumber \\
 & = & \sum_{k=0}^{\left\lfloor \gamma/2\right\rfloor }\frac{f^{(2k)}(0)}{2^{2k}}\sum_{m=0}^{k}\frac{2k-2m+\frac{1}{2}}{m!\left(\frac{1}{2}\right)_{2k-m+1}}\phi_{2k-2m}(x)\,\,\,\,\,\,\,\,\,\,\,\,\,\,\,\,\,\,\,\,\,\,\,\,(\textrm{even }k)\nonumber \\
 &  & +\sum_{k=0}^{\left\lfloor \gamma/2\right\rfloor }\frac{f^{(2k+1)}(0)}{2^{2k+1}}\sum_{m=0}^{k}\frac{2k-2m+\frac{3}{2}}{m!\left(\frac{1}{2}\right)_{2k-m+2}}\phi_{2k-2m+1}(x)\,\,\,\,\,\,\,\,\,\,\,\,\,\,(\textrm{odd }k)\label{eq:odd}\\
 & = & \sum_{m=0}^{\left\lfloor \gamma/2\right\rfloor }\sum_{k=m}^{\left\lfloor \gamma/2\right\rfloor }\frac{2k-2m+\frac{1}{2}}{m!\left(\frac{1}{2}\right)_{2k-m+1}}\frac{f^{(2k)}(0)}{2^{2k}}\phi_{2k-2m}(x)\,\,\,\,\,\,\,\,\,\,\,\,(\textrm{interchanging sums})\nonumber \\
 &  & +\sum_{m=0}^{\left\lfloor \gamma/2\right\rfloor }\sum_{k=m}^{\left\lfloor \gamma/2\right\rfloor }\frac{2k-2m+\frac{3}{2}}{m!\left(\frac{1}{2}\right)_{2k-m+2}}\frac{f^{(2k+1)}(0)}{2^{2k+1}}\phi_{2k-2m+1}(x)\nonumber \\
 & = & \sum_{m=0}^{\left\lfloor \gamma/2\right\rfloor }\sum_{l=0}^{\left\lfloor \gamma/2\right\rfloor }\frac{2l+\frac{1}{2}}{m!\left(\frac{1}{2}\right)_{2l+m+1}}\frac{f^{(2l+2m)}(0)}{2^{2l+2m}}\phi_{2l}(x)\,\,\,\,\,\,\,\,\,\,\,\,\,(\textrm{letting }l=k-m)\nonumber \\
 &  & +\sum_{m=0}^{\left\lfloor \gamma/2\right\rfloor }\sum_{l=0}^{\left\lfloor \gamma/2\right\rfloor }\frac{2l+\frac{3}{2}}{m!\left(\frac{1}{2}\right)_{2l+m+2}}\frac{f^{(2l+2m+1)}(0)}{2^{2l+2m+1}}\phi_{2l+1}(x)\nonumber \\
 & = & \sum_{l=0}^{\left\lfloor \gamma/2\right\rfloor }\sum_{m=0}^{\left\lfloor \gamma/2\right\rfloor }\frac{2l+\frac{1}{2}}{m!\left(\frac{1}{2}\right)_{2l+m+1}}\frac{f^{(2l+2m)}(0)}{2^{2l+2m}}\phi_{2l}(x)\,\,\,\,\,\,\,\,\,\,\,\,(\textrm{interchanging sums})\nonumber \\
 &  & +\sum_{l=0}^{\left\lfloor \gamma/2\right\rfloor }\sum_{m=0}^{\left\lfloor \gamma/2\right\rfloor }\frac{2l+\frac{3}{2}}{m!\left(\frac{1}{2}\right)_{2l+m+2}}\frac{f^{(2l+2m+1)}(0)}{2^{2l+2m+1}}\phi_{2l+1}(x)\nonumber 
\end{eqnarray}
which gives the claim in Lemma A.1. 

For the case of even $\gamma$, note that the term in (\ref{eq:odd})
takes the form 
\begin{eqnarray*}
 &  & \sum_{k=0}^{\left\lfloor \gamma/2\right\rfloor }\frac{f^{(2k+1)}(0)}{2^{2k+1}}\sum_{m=0}^{k}\frac{2k-2m+\frac{3}{2}}{m!\left(\frac{1}{2}\right)_{2k-m+2}}\phi_{2k-2m+1}(x)\\
 & = & \sum_{k=0}^{\left\lfloor \gamma/2\right\rfloor -1}\frac{f^{(2k+1)}(0)}{2^{2k+1}}\sum_{m=0}^{k}\frac{2k-2m+\frac{3}{2}}{m!\left(\frac{1}{2}\right)_{2k-m+2}}\phi_{2k-2m+1}(x)\\
 &  & +\underset{0}{\underbrace{\frac{f^{(\gamma+1)}(0)}{2^{\gamma+1}}}}\sum_{m=0}^{\gamma/2}\frac{\gamma-2m+\frac{3}{2}}{m!\left(\frac{1}{2}\right)_{\gamma-m+2}}\phi_{\gamma-2m+1}(x)
\end{eqnarray*}
and hence the previous derivations go through. 

\subsection{Proof for Lemma \ref{lm:entropy-generalized-H=0000F6lder}\label{subsec:Proof-for-entropy-holder-sub}}

\textbf{The upper bound}. The following derivations generalize Kolmogorov
and Tikhomirov (1959). Any function $f\in\mathcal{U}_{\gamma+1,2}$
can be written as 
\[
f(x+\Delta)=\underset{:=F_{\gamma-1}(x)}{\underbrace{f(x)+\Delta f^{'}(x)+\frac{\Delta^{2}}{2!}f^{''}(x)+\cdots+\frac{\Delta^{\gamma-1}}{\left(\gamma-1\right)!}f^{(\gamma-1)}(x)}}+\frac{\Delta^{\gamma}}{\gamma!}f^{(\gamma)}(z)
\]
where $x,\,x+\Delta\in\left(-1,\,1\right)$ and $z$ is some intermediate
value. Let $REM_{0}(x+\Delta):=f(x+\Delta)-F_{\gamma-1}(x)-\frac{\Delta^{\gamma}}{\gamma!}f^{(\gamma)}(x)$
and note that
\begin{eqnarray}
\left|REM_{0}(x+\Delta)\right| & = & \frac{\left|\Delta\right|^{\gamma}}{\gamma!}\left|f^{(\gamma)}(z)-f^{(\gamma)}(x)\right|\nonumber \\
 & \leq & \frac{\left|\Delta\right|^{\gamma+1}}{\gamma!}R_{\gamma+1}.\label{eq:rem}
\end{eqnarray}
In other words, 
\[
f(x+\Delta)=\sum_{k=0}^{\gamma}\frac{\Delta^{k}}{k!}f^{(k)}(x)+REM_{0}(x+\Delta)
\]
where $\left|REM_{0}(x+\Delta)\right|\leq\frac{\left|\Delta\right|^{\gamma+1}}{\gamma!}R_{\gamma+1}$.
Similarly, any $f^{(i)}\in\mathcal{U}_{\gamma+1-i,2}$ for $1\leq i\leq\gamma$
can be written as
\begin{equation}
f^{(i)}(x+\Delta)=\sum_{k=0}^{\gamma-i}\frac{\Delta^{k}}{k!}f^{(i+k)}(x)+REM_{i}(x+\Delta)\label{eq:a30}
\end{equation}
where $\left|REM_{i}(x+\Delta)\right|\leq\frac{\left|\Delta\right|^{\gamma+1-i}}{\left(\gamma-i\right)!}R_{\gamma+1-i}$. 

For some $\delta_{0},\ldots,\delta_{\gamma}>0$, suppose that $\left|f^{(k)}(x)-g^{(k)}(x)\right|\leq\delta_{k}$
for $k=0,\dots,\gamma$, where $f,g\in\mathcal{U}_{\gamma+1,2}$.
Then we have 
\[
\left|f(x+\Delta)-g(x+\Delta)\right|\leq\sum_{k=0}^{\gamma}\frac{|\Delta|^{k}\delta_{k}}{k!}+2\frac{\left|\Delta\right|^{\gamma+1}}{\gamma!}R_{\gamma+1}.
\]
Let $\left(\max_{k\in\left\{ 1,...,\gamma+1\right\} }\frac{R_{k}}{\left(k-1\right)!}\right)\vee1=:R^{*}$.
Consider $|\Delta|\leq\left(R^{*-1}\delta\right)^{\frac{1}{\gamma+1}}$
and $\delta_{k}=R^{*\frac{k}{\gamma+1}}\delta^{1-\frac{k}{\gamma+1}}$
for $k=0,\dots,\gamma$. Then, 
\begin{eqnarray}
\left|f(x+\Delta)-g(x+\Delta)\right| & \leq & \delta\sum_{k=0}^{\gamma}\left(R^{*}{}^{\frac{-k+k}{\gamma+1}}\frac{1}{k!}\right)+2R^{*}\left|\Delta\right|^{\gamma+1}\nonumber \\
 & \leq & \delta\sum_{k=0}^{\gamma}\frac{1}{k!}+2\delta\leq5\delta.\label{eq:a90}
\end{eqnarray}
Let us consider the following $\left(R^{*-1}\delta\right)^{\frac{1}{\gamma+1}}-$grid
of points in $\left[-1,\,1\right]$:
\[
x_{-s}<x_{-s+1}\cdots<x_{-1}<x_{0}<x_{1}<\cdots<x_{s-1}<x_{s},
\]
with 
\[
x_{0}=0\textrm{ and }s\precsim\left(R^{*-1}\delta\right)^{\frac{-1}{\gamma+1}}.
\]
It suffices to cover the $k$th derivatives of functions in $\mathcal{U}_{\gamma+1,2}$
within $\delta_{k}-$precision at each grid point. Then by (\ref{eq:a90}),
we obtain a $5\delta-$cover of $\mathcal{U}_{\gamma+1,2}$. Following
the arguments in Kolmogorov and Tikhomirov (1959), bounding $N_{\infty}\left(\delta,\,\mathcal{U}_{\gamma+1,2}\right)$
can be reduced to bounding the cardinality of 
\[
\Lambda=\left\{ \left(\left\lfloor \frac{f^{(k)}\left(x_{i}\right)}{\delta_{k}}\right\rfloor ,-s\leq i\leq s,\,0\leq k\leq\gamma\right):f\in\mathcal{U}_{\gamma+1,2}\right\} 
\]
with $\left\lfloor x\right\rfloor $ denoting the largest integer
smaller than or equal to $x$. Starting with $x_{0}=0$, the number
of possible values of the vector $\left(\left\lfloor \frac{f^{(k)}\left(x_{0}\right)}{\delta_{k}}\right\rfloor \right)_{k=0}^{\gamma}$
when $f$ ranges over $\mathcal{U}_{\gamma+1,2}$ is $1$. For $i=1,...,s$,
given the value of $\left(\left\lfloor \frac{f^{(k)}\left(x_{i-1}\right)}{\delta_{k}}\right\rfloor \right)_{k=0}^{\gamma}$,
let us count the number of possible values of $\left(\left\lfloor \frac{f^{(k)}\left(x_{i}\right)}{\delta_{k}}\right\rfloor \right)_{k=0}^{\gamma}$.
The counting for $\left(\left\lfloor \frac{f^{(k)}\left(x_{-i}\right)}{\delta_{k}}\right\rfloor \right)_{k=0}^{\gamma}$
is similar. For each $0\leq k\leq\gamma$, let $B_{k,i-1}:=\left\lfloor \frac{f^{(k)}\left(x_{i-1}\right)}{\delta_{k}}\right\rfloor $.
Observe that $B_{k,i-1}\delta_{k}\leq f^{(k)}\left(x_{0}\right)<\left(B_{k,i-1}+1\right)\delta_{k}$. 

Taking (\ref{eq:a30}) with $x=x_{i-1}$ and $\Delta=x_{i}-x_{i-1}$
gives 
\[
\left|f^{(i)}\left(x_{i}\right)-\sum_{k=0}^{\gamma-i}\frac{\Delta^{k}}{k!}f^{(i+k)}\left(x_{i-1}\right)\right|\leq\frac{\left|\Delta\right|^{\gamma+1-i}}{\left(\gamma-i\right)!}R_{\gamma+1-i}.
\]
As a result, 
\begin{align*}
 & \left|f^{(i)}\left(x_{i}\right)-\sum_{k=0}^{\gamma-i}\frac{\Delta^{k}}{k!}B_{i+k,i-1}\right|\\
 & \leq\left|f^{(i)}\left(x_{i}\right)-\sum_{k=0}^{\gamma-i}\frac{\Delta^{k}}{k!}f^{(i+k)}\left(x_{i-1}\right)\right|+\left|\sum_{k=0}^{\gamma-i}\frac{\Delta^{k}}{k!}\left(f^{(i+k)}\left(x_{i-1}\right)-B_{i+k,i-1}\right)\right|\\
 & \leq\frac{\left|\Delta\right|^{\gamma+1-i}}{\left(\gamma-i\right)!}R_{\gamma+1-i}+\sum_{k=0}^{\gamma-i}\frac{|\Delta|^{k}}{k!}\delta_{i+k}\\
 & \leq\left(R^{*-1}\delta\right)^{1-\frac{i}{\gamma+1}}R_{\gamma+1-i}^{*}+\sum_{k=0}^{\gamma-i}\left[\frac{1}{k!}\left(R^{*-1}\delta\right)^{\frac{k}{\gamma+1}}R^{*\frac{i+k}{\gamma+1}}\delta^{1-\frac{i+k}{\gamma+1}}\right]\\
 & \leq R^{*\frac{i}{\gamma+1}}\delta^{1-\frac{i}{\gamma+1}}+R^{*\frac{i}{\gamma+1}}\delta^{1-\frac{i}{\gamma+1}}\sum_{k=0}^{\gamma-i}\frac{1}{k!}\leq4\delta_{i}.
\end{align*}
Hence, the number of possible values of $\left(\left\lfloor \frac{f^{(k)}\left(x_{i}\right)}{\delta_{k}}\right\rfloor \right)_{k=0}^{\gamma}$
is at most $4$ given the value of $\left(\left\lfloor \frac{f^{(k)}\left(x_{i-1}\right)}{\delta_{k}}\right\rfloor \right)_{k=0}^{\gamma}$.
Consequently, we have 
\[
\textrm{card}\left(\Lambda\right)\precsim4^{2s}\precsim16^{\left(R^{*-1}\delta\right)^{\frac{-1}{\gamma+1}}}
\]
which implies 
\begin{equation}
\log N_{2}\left(\delta,\,\mathcal{U}_{\gamma+1,2}\right)\leq\log N_{\infty}\left(\delta,\,\mathcal{U}_{\gamma+1,2}\right)\precsim R^{*\frac{1}{\gamma+1}}\delta^{\frac{-1}{\gamma+1}}.\label{eq:upper2}
\end{equation}
\textbf{The lower bound}. In the derivation of the lower bound, Kolmogorov
and Tikhomirov (1959) considers a $\delta^{\frac{1}{\gamma+1}}-$grid
of points
\[
\cdots<\underline{a}_{1}<\overline{a}_{1}<\underline{a}_{2}<\overline{a}_{2}<\cdots<\underline{a}_{2s}<\overline{a}_{2s}
\]
where $\overline{a}_{i}-\underline{a}_{i}=\delta^{\frac{1}{\gamma+1}}$
and $s\succsim\delta^{\frac{-1}{\gamma+1}}$. Recall that we have
previously considered a $\left(R^{*-1}\delta\right)^{\frac{1}{\gamma+1}}-$grid
of points in $\left[-1,\,1\right]$ in the derivation of the upper
bound for $\log N_{\infty}\left(\delta,\,\mathcal{U}_{\gamma+1,2}\right)$.
To obtain a lower bound for $\log M_{\infty}\left(\delta,\,\mathcal{U}_{\gamma+1,2}\right)$
with the same scaling as our upper bound, the key modification we
need is to replace the $\overline{a}_{i}-\underline{a}_{i}=\delta^{\frac{1}{\gamma+1}}$
with $\overline{a}_{i}-\underline{a}_{i}=\left(R^{*-1}\delta\right)^{\frac{1}{\gamma+1}}$
and $s\succsim\delta^{\frac{-1}{\gamma+1}}$ with $s\succsim R^{*\frac{1}{\gamma+1}}\delta^{\frac{-1}{\gamma+1}}$.
The rest of the arguments are similar to those in Kolmogorov and Tikhomirov
(1959). In particular, let us consider 
\[
f_{\lambda}(x)=R^{*}\sum_{i=1}^{2s}\lambda_{i}\left(\overline{a}_{i}-\underline{a}_{i}\right)^{\gamma+1}h_{0}\left(\frac{x-\underline{a}_{i}}{\overline{a}_{i}-\underline{a}_{i}}\right)
\]
where $\lambda_{i}\in\left\{ 0,1\right\} $ and $\lambda\in\left\{ 0,1\right\} ^{2s}$,
and $h_{0}$ is a function on $\mathbb{R}$ satisfying: (1) $h_{0}$
restricted to $\left[-1,\,1\right]$ belongs to $\mathcal{U}_{\gamma+1,2}$;
(2) $h_{0}(x)=0$ for $x\notin[0,\,1]$ and $h_{0}(x)>0$ for $x\in[0,\,1]$;
(3) $h_{0}\left(\frac{1}{2}\right)=\max_{x\in\left[0,\,1\right]}h_{0}(x)=R_{0}$.
As an example, we can take $h_{0}(x)=\begin{cases}
0 & x\notin[0,\,1]\\
bx^{\gamma+1}(1-x)^{\gamma+1} & x\in[0,\,1]
\end{cases}$ for some properly chosen $b$ that ensures (1). Note that the functions
$h(x):=R^{*}\left(\overline{a}_{i}-\underline{a}_{i}\right)^{\gamma+1}h_{0}\left(\frac{x-\underline{a}_{i}}{\overline{a}_{i}-\underline{a}_{i}}\right)$
and also $f_{\lambda}(x)$ belong to $\mathcal{U}_{\gamma+1,2}$ if
$\delta\in\left(0,\,1\right)$. For any distinct $\lambda,\,\lambda^{'}\in\left\{ 0,1\right\} ^{2s}$,
we have 
\[
\left|f_{\lambda}-f_{\lambda^{'}}\right|_{\infty}=R^{*}\left(\overline{a}_{i}-\underline{a}_{i}\right)^{\gamma+1}h_{0}\left(\frac{1}{2}\right)=R_{0}\delta.
\]
If $R_{0}\succsim1$, then $R_{0}\delta\succsim\delta$ and 
\[
\log M_{\infty}\left(\delta,\,\mathcal{U}_{\gamma+1,2}\right)\succsim R^{*\frac{1}{\gamma+1}}\delta^{\frac{-1}{\gamma+1}}.
\]
If $R_{0}\precsim1$, then we obtain
\[
\log M_{\infty}\left(R_{0}\delta,\,\mathcal{U}_{\gamma+1,2}\right)\succsim R^{*\frac{1}{\gamma+1}}\delta^{\frac{-1}{\gamma+1}}
\]
which implies that 
\[
\log M_{\infty}\left(\delta,\,\mathcal{U}_{\gamma+1,2}\right)\succsim R^{*\frac{1}{\gamma+1}}\left(\frac{\delta}{R_{0}}\right)^{\frac{-1}{\gamma+1}}.
\]
Standard argument in the literature based on the Vasharmov-Gilbert
Lemma further gives 
\begin{equation}
\log M_{2}\left(\delta,\,\mathcal{U}_{\gamma+1,2}\right)\succsim\begin{cases}
R^{*\frac{1}{\gamma+1}}\delta^{\frac{-1}{\gamma+1}} & \textrm{if }R_{0}\succsim1\\
\left(R^{*}R_{0}\right)^{\frac{1}{\gamma+1}}\delta^{\frac{-1}{\gamma+1}} & \textrm{if }R_{0}\precsim1
\end{cases}.\label{eq:lower2}
\end{equation}

To show the last two bounds in Lemma \ref{lm:entropy-generalized-H=0000F6lder},
we apply the same arguments for showing claim (iii) in Lemma \ref{lm:entropy-generalized-polynomial}.
\\
\\
\textbf{Remark A}. Sections \ref{subsec:Proof-for-entropy-poly} and
\ref{subsec:Proof-for-entropy-holder-sub} derive bounds for $\mathcal{U}_{\gamma+1,1}$
and $\mathcal{U}_{\gamma+1,2}$ on $\left[-1,\,1\right]$. These derivations
can be easily extended for $\mathcal{U}_{\gamma+1,1}$ and $\mathcal{U}_{\gamma+1,2}$
on a general bounded interval $\left[c_{1},\,c_{2}\right]$, where
$c_{1}$ and $c_{2}$ are universal constants that are independent
of $\gamma$ and $\left\{ R_{k}\right\} _{k=0}^{\gamma+1}$. In particular,
the resulting bounds have the same scaling (in terms of $\delta$,
$\gamma$ and $\left\{ R_{k}\right\} _{k=0}^{\gamma+1}$) as those
in Sections \ref{subsec:Proof-for-entropy-poly} and \ref{subsec:Proof-for-entropy-holder-sub}.

\subsection{Proof for Lemma \ref{lm:entropy-ellipsoid-subclass}\label{subsec:Proof-for-entropy_ellipsoid_sub}}

In the special case of $R_{\gamma+1}=1$, the argument below sharpens
the upper bound for $\log N_{2}\left(\delta,\,\mathcal{H}_{\gamma+1}\right)$
in Wainwright (2019) from $\left(\gamma\vee1\right)\delta^{-\frac{1}{\gamma+1}}$
to $\delta^{-\frac{1}{\gamma+1}}$. We find the cause of the gap lies
in that the ``pivotal'' eigenvalue (that balances the ``estimation
error'' and the ``approximation error'' from truncating for a given
resolution $\delta$) in Wainwright (2019) is not optimal. We close
the gap by finding the optimal ``pivotal'' eigenvalue.

More generally, for the case of $R_{\gamma+1}\precsim\gamma+1$, we
consider two different truncations, one giving the upper bound $\delta^{\frac{-1}{\gamma+1}}$
and the other giving the lower bound $\left(R_{\gamma+1}\delta^{-1}\right)^{\frac{1}{\gamma+1}}$.
Note that $\left(R_{\gamma+1}\delta^{-1}\right)^{\frac{1}{\gamma+1}}\asymp\delta^{\frac{-1}{\gamma+1}}$
when $R_{\gamma+1}\asymp1$. For the case of $R_{\gamma+1}\succsim\gamma+1$,
we use only one truncation to show that both the upper bound and the
lower bound scale as $\left(R_{\gamma+1}\delta^{-1}\right)^{\frac{1}{\gamma+1}}$.

In view of (\ref{eq:Ellipsoid-1}), given $\left(\phi_{m}\right)_{m=1}^{\infty}$
and $\left(\mu_{m}\right)_{m=1}^{\infty}$, to compute $N_{2}\left(\delta,\,\mathcal{H}_{\gamma+1}\right)$,
it suffices to compute $N_{2}\left(\delta,\,\mathcal{E}_{\gamma+1}\right)$
where 
\[
\mathcal{E}_{\gamma+1}=\left\{ \left(\theta_{m}\right)_{m=1}^{\infty}:\sum_{m=1}^{\infty}\frac{\theta_{m}^{2}}{\mu_{m}}\leq R_{\gamma+1}^{2},\,\mu_{m}=\left(cm\right)^{-2\left(\gamma+1\right)}\right\} .
\]
Let us introduce the $M-$dimensional ellipsoid 
\[
\overline{\mathcal{E}}_{\gamma+1}=\left\{ \left(\theta_{m}\right)_{m=1}^{M}\textrm{ coincide with the first \ensuremath{M}}\textrm{elements of \ensuremath{\left(\theta_{m}\right)_{m=1}^{\infty}} in \ensuremath{\mathcal{E}_{\gamma+1}} }\right\} 
\]
where $M\left(=M\left(\gamma+1,\,\delta\right)\right)$ is the smallest
integer such that, for a given resolution $\delta>0$ and weight $w_{\gamma+1}$,
$w_{\gamma+1}^{2}\delta^{2}\geq\mu_{M}$. In other words, $\mu_{m}\geq w_{\gamma+1}^{2}\delta^{2}$
for all indices $m\leq M$. Consequently, we have: (1) 
\begin{equation}
\mathbb{B}_{2}^{M}\left(w_{\gamma+1}R_{\gamma+1}\delta\right)\subseteq\overline{\mathcal{E}}_{\gamma+1};\label{eq:contain-1}
\end{equation}
(2) $\mu_{M-1}=\left(c\left(M-1\right)\right)^{-2\left(\gamma+1\right)}>w_{\gamma+1}^{2}\delta^{2}$
and $\mu_{M-1}=\left(c\left(M+1\right)\right)^{-2\left(\gamma+1\right)}<w_{\gamma+1}^{2}\delta^{2}$,
which yield 
\begin{equation}
M\asymp\left(w_{\gamma+1}\delta\right)^{-\frac{1}{\gamma+1}}.\label{eq:M-1-1}
\end{equation}
Note that (\ref{eq:contain-1}), (\ref{eq:M-1-1}), and the fact $\mathcal{E}_{\gamma+1}\supseteq\overline{\mathcal{E}}_{\gamma+1}$
give 
\begin{align}
\log N_{2}\left(\delta,\,\mathcal{E}_{\gamma+1}\right) & \geq\log N_{2}\left(\delta,\,\overline{\mathcal{E}}_{\gamma+1}\right)\nonumber \\
 & \succsim M\log\left(w_{\gamma+1}R_{\gamma+1}\right)\nonumber \\
 & \asymp\left(w_{\gamma+1}\delta\right)^{-\frac{1}{\gamma+1}}\log\left(w_{\gamma+1}R_{\gamma+1}\right).\label{eq:50}
\end{align}

In the following, let $A_{1}+A_{2}:=\left\{ a_{1}+a_{2}:a_{1}\in A_{1},\,a_{2}\in A_{2}\right\} $
for sets $A_{1}$ and $A_{2}$. For the upper bound, we have 
\begin{eqnarray}
N_{2}\left(\delta,\,\overline{\mathcal{E}}_{\gamma+1}\right) & \leq & \frac{\textrm{vol}\left(\frac{2}{\delta}\overline{\mathcal{E}}_{\gamma+1}+\mathbb{B}_{2}^{M}\left(1\right)\right)}{\textrm{vol}\left(\mathbb{B}_{2}^{M}\left(1\right)\right)}\nonumber \\
 & \leq & \left(\frac{2}{\delta}\right)^{M}\frac{\textrm{vol}\left(\overline{\mathcal{E}}_{\gamma+1}+\mathbb{B}_{2}^{M}\left(\frac{\delta}{2}\right)\right)}{\textrm{vol}\left(\mathbb{B}_{2}^{M}\left(1\right)\right)}\nonumber \\
 & \leq & \left(\frac{2}{\delta}\right)^{M}\max\left\{ \frac{\textrm{vol}\left(2\overline{\mathcal{E}}_{\gamma+1}\right)}{\textrm{vol}\left(\mathbb{B}_{2}^{M}\left(1\right)\right)},\,\frac{\textrm{vol}\left(2\mathbb{B}_{2}^{M}\left(\frac{\delta}{2}\right)\right)}{\textrm{vol}\left(\mathbb{B}_{2}^{M}\left(1\right)\right)}\right\} \nonumber \\
 & \leq & \max\left\{ \left(\frac{4R_{\gamma+1}}{\delta}\right)^{M}\prod_{m=1}^{M}\sqrt{\mu_{m}},\,2^{M}\right\} \label{eq:54-1}
\end{eqnarray}
where the first inequality follows from the standard volumetric argument,
and the last inequality follows from the standard result for the volume
of ellipsoids. The fact $\mu_{m}=\left(cm\right)^{-2\left(\gamma+1\right)}$
and the elementary inequality $\sum_{m=1}^{M}\log m\geq M\log M-M$
give
\begin{eqnarray}
\log\left[\left(\frac{4R_{\gamma+1}}{\delta}\right)^{M}\prod_{m=1}^{M}\sqrt{\mu_{m}}\right] & \leq & M\left(\log\left(4R_{\gamma+1}\right)+\gamma+1\right)+\nonumber \\
 &  & M\left(\log\frac{1}{\delta}-\left(\gamma+1\right)\log\left(cM\right)\right)\nonumber \\
 & = & M\left(\log\left(4R_{\gamma+1}\right)+\gamma+1\right)+\nonumber \\
 &  & M\left(\log\frac{1}{\delta}-\left(\gamma+1\right)\log\left(cM\right)+\log\frac{1}{w_{\gamma+1}}-\log\frac{1}{w_{\gamma+1}}\right)\nonumber \\
 & \leq & M\left(\log4R_{\gamma+1}+\gamma+1\right)+M\log w_{\gamma+1}\nonumber \\
 & \precsim & M\log\left(w_{\gamma+1}\left(\left(\gamma+1\right)\vee R_{\gamma+1}\right)\right)\label{eq:33}
\end{eqnarray}
where we have used the fact $\mu_{M}=\left(cM\right)^{-2\left(\gamma+1\right)}\leq w_{\gamma+1}^{2}\delta^{2}$
in the second inequality. Inequalities (\ref{eq:M-1-1}), (\ref{eq:54-1})
and (\ref{eq:33}) together yield 
\[
\log N_{2}\left(\delta,\,\overline{\mathcal{E}}_{\gamma+1}\right)\precsim\left(w_{\gamma+1}\delta\right)^{-\frac{1}{\gamma+1}}\max\left\{ \log\left(w_{\gamma+1}\left(\left(\gamma+1\right)\vee R_{\gamma+1}\right)\right),\,\log2\right\} .
\]

For any $\theta\in\mathcal{E}_{\gamma+1}$, note that for a given
$\delta$, we have
\begin{equation}
\sum_{m=M+1}^{\infty}\theta_{m}^{2}\leq\mu_{M}\sum_{m=M+1}^{\infty}\frac{\theta_{m}^{2}}{\mu_{m}}\leq w_{\gamma+1}^{2}R_{\gamma+1}^{2}\delta^{2}.\label{eq:approx error-1}
\end{equation}
To cover $\mathcal{E}_{\gamma+1}$ within $\left(1+w_{\gamma+1}^{2}R_{\gamma+1}^{2}\right)^{\frac{1}{2}}\delta-$precision,
we find a smallest $\delta-$cover of $\overline{\mathcal{E}}_{\gamma+1}$,
$\left\{ \theta^{1},...,\theta^{N}\right\} $, such that for any $\theta\in\mathcal{E}_{\gamma+1}$,
there exists some $i$ from the covering set with
\[
\left|\theta-\theta^{i}\right|_{2}^{2}\leq\sum_{m=1}^{M}\left(\theta_{m}-\theta_{m}^{i}\right)^{2}+w_{\gamma+1}^{2}R_{\gamma+1}^{2}\delta^{2}\leq\left(1+w_{\gamma+1}^{2}R_{\gamma+1}^{2}\right)\delta^{2}
\]
where we have used (\ref{eq:approx error-1}). Consequently, we have
\begin{eqnarray}
 &  & \log N_{2}\left(\delta,\,\mathcal{E}_{\gamma+1}\right)\nonumber \\
 & \precsim & \log N_{2}\left(\delta\left(1+w_{\gamma+1}^{2}R_{\gamma+1}^{2}\right)^{\frac{-1}{2}},\,\overline{\mathcal{E}}_{\gamma+1}\right)\nonumber \\
 & \precsim & \left(w_{\gamma+1}\delta\left(1+w_{\gamma+1}^{2}R_{\gamma+1}^{2}\right)^{\frac{-1}{2}}\right)^{-\frac{1}{\gamma+1}}\max\left\{ \log\left(w_{\gamma+1}\left(\left(\gamma+1\right)\vee R_{\gamma+1}\right)\right),\,\log2\right\} .\label{eq:54}
\end{eqnarray}
\textbf{Case 1}: $R_{\gamma+1}\succsim\gamma+1$. Setting $w_{\gamma+1}\asymp R_{\gamma+1}^{-1}$
in (\ref{eq:50}) and (\ref{eq:54}) solves 
\begin{eqnarray}
 &  & \left(w_{\gamma+1}\delta\left(1+w_{\gamma+1}^{2}R_{\gamma+1}^{2}\right)^{\frac{-1}{2}}\right)^{-\frac{1}{\gamma+1}}\max\left\{ \log\left(w_{\gamma+1}\left(\left(\gamma+1\right)\vee R_{\gamma+1}\right)\right),\,\log2\right\} \nonumber \\
 & \asymp & \left(w_{\gamma+1}\delta\right)^{-\frac{1}{\gamma+1}}\log\left(w_{\gamma+1}R_{\gamma+1}\right)\label{eq:equation}
\end{eqnarray}
and gives 
\[
\log N_{2}\left(\delta,\,\mathcal{E}_{\gamma+1}\right)\asymp\left(R_{\gamma+1}\delta^{-1}\right)^{\frac{1}{\gamma+1}}.
\]
\textbf{Case 2}: $R_{\gamma+1}\precsim\gamma+1$. Setting $w_{\gamma+1}\asymp\left(\gamma+1\right)^{-1}$
in (\ref{eq:54}) gives 
\[
\log N_{2}\left(\delta,\,\mathcal{E}_{\gamma+1}\right)\precsim\delta^{\frac{-1}{\gamma+1}}.
\]
Note that the lower bound obtained by setting $w_{\gamma+1}\asymp\left(\gamma+1\right)^{-1}$
in (\ref{eq:50}) is not particularly useful. Instead, we consider
a different truncation with $w_{\gamma+1}\asymp R_{\gamma+1}^{-1}$.
Then (\ref{eq:50}) with $w_{\gamma+1}\asymp R_{\gamma+1}^{-1}$ gives
\[
\log N_{2}\left(\delta,\,\mathcal{E}_{\gamma+1}\right)\succsim R_{\gamma+1}^{\frac{1}{\gamma+1}}\delta^{\frac{-1}{\gamma+1}}.
\]

To show the last claim in Lemma \ref{lm:entropy-ellipsoid-subclass},
we apply the same arguments for showing claim (iii) in Lemma \ref{lm:entropy-generalized-polynomial}. 

\section{Proofs for Section \ref{sec:Minimax-standard} and Appendix \ref{sec:MISE-nonstandard}
\label{sec:Proofs-for-MISE}}

\subsection{Proof for Theorem \ref{thm:MISE-lower-standard}\label{subsec:Proof-for-MISE-lower-standard}}

The lower bound $\left(\frac{\sigma^{2}}{n}\right)^{\frac{2\left(\gamma+1\right)}{2\left(\gamma+1\right)+1}}$
in part (i) of Theorem \ref{thm:MISE-lower-standard} is the standard
one in the literature. One simply applies Lemma C.1(ii) in Appendix
\ref{sec:Supporting-lemmas} and the standard metric entropy bound
$\left(\frac{1}{\delta}\right)^{\frac{1}{\gamma+1}}$. Therefore,
we only show part (ii) of Theorem \ref{thm:MISE-lower-standard} below.
\\
\textbf{}\\
\textbf{Standard Sobolev $\mathcal{S}_{\gamma+1}$}. Let us first
define a subset of $\mathcal{S}_{\gamma+1}$: 
\begin{align*}
\overline{\mathcal{S}}_{\gamma+1} & :=\{f:\,\left[0,\,1\right]\rightarrow\mathbb{R}\vert\,f\textrm{ is \ensuremath{\gamma+1} times differentiable a.e.,}\\
 & f^{(k)}\textrm{ is absolutely continuous and},\\
 & \left|f\right|_{\mathcal{H},k+1}\leq\overline{C},\:\textrm{for all }k=0,...,\gamma\}
\end{align*}
where 
\begin{equation}
\left|f\right|_{\mathcal{H},k+1}=\sqrt{\sum_{j=0}^{k}\left(f^{(j)}(0)\right)^{2}+\int_{0}^{1}\left[f^{(k+1)}\left(t\right)\right]^{2}dt}.\label{eq:sobolev norm-1}
\end{equation}
Note that $\overline{\mathcal{S}}_{k+1}\subseteq\overline{\mathcal{S}}_{k^{'}+1}$
for any $k\geq k^{'}$. 

We apply Lemma C.1(i) in Appendix \ref{sec:Supporting-lemmas} and
the construction in Appendix \ref{subsec:Proof-for-Lemma B.1}. To
construct the packing subset in Lemma C.1(i), we construct $\mathcal{M}$,
a set in the $\left(\gamma^{*}+1\right)$th order polynomial subclass.
By (\ref{eq:KL}) in Appendix \ref{sec:Supporting-lemmas}, for a
set 
\[
\left\{ f^{1},f^{2},...,f^{M}\right\} =\mathcal{M}\subseteq\mathcal{U}_{\gamma^{*}+1,1}^{*}
\]
(where $\mathcal{U}_{\gamma^{*}+1,1}^{*}$ is defined the same way
as (\ref{eq:polynomial subspace in Sobolev})), we have 
\begin{eqnarray*}
D_{KL}\left(\mathbb{P}_{f^{j}}\times\mathbb{P}_{X}\parallel\mathbb{P}_{f^{l}}\times\mathbb{P}_{X}\right) & = & \frac{n}{2\sigma^{2}}\left|f^{j}-f^{l}\right|_{2,\mathbb{P}}^{2}.
\end{eqnarray*}
Let us consider the packing set consisting of $M_{0}$ elements in
the form (\ref{eq:M-packing}) for $\mathcal{U}_{\gamma^{*}+1,1}^{*}$. 

Following the notation in Appendix \ref{subsec:Proof-for-entropy-poly}
(the derivations for $\underline{B}_{2}$), we have 
\begin{eqnarray}
\left|f^{j}-f^{l}\right|_{2,\mathbb{P}}^{2} & = & \left|f_{\lambda_{j}}^{*}-f_{\lambda_{l}}^{*}\right|_{2,\mathbb{P}}^{2}\nonumber \\
 & = & \int_{0}^{1}\left[\left(\theta_{\tilde{k}}^{*j}-\theta_{\tilde{k}}^{*l}\right)x^{\tilde{k}}+\sum_{k\in\left\{ 0,...,\gamma^{*}\right\} \setminus\tilde{k}}\left(1\left\{ \lambda_{j,k}\neq\lambda_{l,k}\right\} \theta_{k}^{*}x^{k}\right)\right]^{2}p(x)dx\nonumber \\
 & \precsim & \left(\theta_{\tilde{k}}^{*j}-\theta_{\tilde{k}}^{*l}\right)^{2}+\left(\sum_{k\in\left\{ 0,...,\gamma^{*}\right\} \setminus\tilde{k}}\theta_{k}^{*}\right)^{2}\nonumber \\
 & \precsim & \delta^{2}\left[\left(\tilde{k}+1\right)\vee\sum_{k=0}^{\gamma^{*}}\frac{R_{k}}{k!}\right]^{2}\nonumber \\
 & \precsim & \delta^{2}\textrm{ }\quad\textrm{for any }j,l\in\mathcal{M},j\neq l\label{eq:f-packing_standard}
\end{eqnarray}
where the last line follows from that $R_{0}=\frac{\overline{C}}{6}$
and $R_{k}=\frac{\overline{C}}{6\gamma^{*}}$ for $k=1,...,\gamma^{*}$,
in which case we have $\tilde{k}=0$ and $\sum_{k=0}^{\gamma^{*}}\frac{R_{k}}{k!}\asymp1$.
Recall that the lower bound $\log M_{0}\succsim\gamma^{*}+1$ in Appendix
\ref{subsec:Proof-for-entropy-poly} holds for all $\delta$ such
that $\frac{R_{\tilde{k}}}{\tilde{k}!\delta\left[\left(\tilde{k}+1\right)\vee\sum_{k=0}^{\gamma^{*}}\frac{R_{k}}{k!}\right]}\geq c2^{\gamma^{*}+1}$
and $\delta\left[\left(\tilde{k}+1\right)\vee\sum_{k=0}^{\gamma^{*}}\frac{R_{k}}{k!}\right]\leq\frac{2}{3}\frac{R_{\tilde{k}}}{\tilde{k}!}$,
where $\tilde{k}=0$, $\sum_{k=0}^{\gamma^{*}}\frac{R_{k}}{k!}\asymp1$
and $c\in\left(1,\,\infty\right)$ is a universal constant. Let us
take $\delta^{2}=c^{'2}\frac{\sigma^{2}\left(\gamma^{*}+1\right)}{n}$
for a sufficiently small universal constant $c^{'}\in(0,\,1]$ such
that $\left[c^{'}cR_{0}^{-1}\left(1\vee\sum_{k=0}^{\gamma^{*}}\frac{R_{k}}{k!}\right)\right]^{2}=:c_{0}\in(0,\,1]$,
$c^{'}\left(1\vee\sum_{k=0}^{\gamma^{*}}\frac{R_{k}}{k!}\right)\sqrt{\frac{\sigma^{2}\left(\gamma^{*}+1\right)}{n}}\leq\frac{2}{3}R_{0}$,
and
\begin{align*}
\frac{1}{M^{2}}\sum_{j,l\in\left\{ 1,...,M\right\} }D_{KL}\left(\mathbb{P}_{f^{j}}\times\mathbb{P}_{X}\parallel\mathbb{P}_{f^{l}}\times\mathbb{P}_{X}\right) & \leq\underline{C}_{0}\left(\gamma^{*}+1\right),\\
\log M=\log\left|\mathcal{M}\right| & \geq\overline{C}_{0}\left(\gamma^{*}+1\right),
\end{align*}
where the positive universal constants satisfy $\underline{C}_{0}\leq\frac{1}{2}\overline{C}_{0}$.
Therefore,
\[
\delta^{2}\left(1-\frac{\log2+\frac{1}{M^{2}}\sum_{j,l\in\left\{ 1,...,M\right\} }D_{KL}\left(\mathbb{P}_{f^{j}}\times\mathbb{P}_{X}\parallel\mathbb{P}_{f^{l}}\times\mathbb{P}_{X}\right)}{\log M}\right)\geq\underline{c}\frac{\sigma^{2}\left(\gamma^{*}+1\right)}{n}
\]
for some universal constant $\underline{c}\in(0,\,1]$. Note that
we can choose $c^{'}$ in $\delta$ above to be small enough such
that $\underline{c}\leq\frac{1}{2}\underline{c}_{0}$ as stated in
Theorem \ref{thm:MISE-lower-standard}(ii). \\
\textbf{}\\
\textbf{Standard Hölder $\mathcal{U}_{\gamma+1}$}. The proof for
the standard $\mathcal{U}_{\gamma+1}$ is identical to the proof for
$\mathcal{S}_{\gamma+1}$ shown previously. We consider the $\delta-$packing
set consisting of $M_{0}$ elements in the form (\ref{eq:M-packing})
for $\mathcal{U}_{\gamma^{*}+1,1}^{*}$. 

\subsubsection{Allowing for heteroscedasticity and non-Gaussian noise\label{subsec:proof-allowing-for-heter-nonGaussian}}

Parallel to the \textit{standard} Hölder class in definition (\ref{Standard smoothness classes}),
we let $\sigma(\cdot)$ in (\ref{eq:heter model}) range over $\mathcal{U}_{\gamma+1}$
with $R_{k}=\overline{C}\asymp1$ for all $k=0,...,\gamma+1$ such
that the function members take positive values in $\left[\underline{\sigma},\,\overline{C}\right]$
and $\underline{\sigma}\asymp1$; we denote this class as the standard
$\mathcal{U}_{\gamma+1}^{+}$. Parallel to the standard Sobolev class
in definition (\ref{Standard smoothness classes}), we let $\sigma(\cdot)$
in (\ref{eq:heter model}) range over
\begin{align}
\mathcal{S}_{\gamma+1}^{+}:= & \{g:\,\left[0,\,1\right]\rightarrow\mathbb{R}^{+}\vert\,g\textrm{ is \ensuremath{\gamma+1} times differentiable a.e.,}\nonumber \\
 & g(x)\in\left[\underline{\sigma},\,\overline{C}\right]\textrm{ for all }x\in\left[0,\,1\right]\textrm{ such that }\underline{\sigma}\asymp1,\nonumber \\
 & g^{(\gamma)}\textrm{ is absolutely continuous and}\left|g\right|_{\mathcal{H},\gamma+1}\leq\overline{C}\asymp1\}\label{eq:sobolev_std-sigma}
\end{align}
with $\left|g\right|_{\mathcal{H},\gamma+1}$ defined in (\ref{eq:sobolev norm}).
Note that $\mathcal{U}_{\gamma+1}^{+}$ and $\mathcal{S}_{\gamma+1}^{+}$
defined above consist of scale functions that are bounded away from
zero and bounded from above, similar to that $\sigma\asymp1$ in the
homoscedasstic case under Assumption \ref{Assumption 1}. Given (\ref{eq:cover_KL-hetero})
and (\ref{eq:minimax_hetero}), the lower bound $\left(\frac{\sigma^{2}}{n}\right)^{\frac{2\left(\gamma+1\right)}{2\left(\gamma+1\right)+1}}$
is again the standard one in the literature. Given (\ref{eq:KL_hetero})
and (\ref{eq:minimax_hetero_0}), the arguments are very similar to
the previous for the lower bound $\frac{\sigma^{2}\left(\gamma^{*}+1\right)}{n}$.
We spell out the slight differences below. 

Let us consider
\begin{equation}
\mathcal{U}_{\gamma+1,1}^{*+}=\left\{ f=\sum_{k=0}^{\gamma_{1}}\theta_{2k}x^{2k}:\textrm{\ensuremath{\left(\theta_{2k}\right)_{k=0}^{\gamma_{1}}}\ensuremath{\ensuremath{\in}}\ensuremath{\mathcal{P}_{\gamma_{1}}^{*}}},\,x\in\left[0,\,1\right]\right\} \label{eq:poly+}
\end{equation}
where $2\gamma_{1}=\gamma$ (if $\gamma$ is even) and $2\gamma_{1}+1=\gamma$
(if $\gamma$ is odd). Note that $\mathcal{U}_{\gamma+1,1}^{*+}\subseteq\mathcal{S}_{\gamma+1}^{+}$
and $\mathcal{U}_{\gamma+1,1}^{*+}\subseteq\mathcal{U}_{\gamma+1}^{+}$
if we define the $\left(\gamma_{1}+1\right)-$dimensional polyhedron
\[
\ensuremath{\mathcal{P}_{\gamma_{1}}^{*+}=\left\{ \ensuremath{\left(\theta_{2k}\right)_{k=0}^{\gamma_{1}}}:\theta_{0}\in\left[\underline{\sigma},\,\frac{\overline{C}}{6}\right],\,\theta_{2k}\in\left[0,\,\frac{\overline{C}}{6\gamma_{1}}\frac{1}{\left(2k\right)!}\right],\,k=1,...,\gamma_{1}\right\} .}
\]

For $\mathcal{U}_{\gamma^{*}+1,1}^{*+}$ (by letting $\gamma=\gamma^{*}$
in (\ref{eq:poly+})), we consider the packing set $\mathcal{M}^{+}$
which consists of $M_{0}^{+}$ elements in the form (\ref{eq:M-packing}).
Following the previous derivation that leads to (\ref{eq:f-packing_standard}),
we have 
\begin{align*}
\left|\sigma^{j^{'}}-\sigma^{l^{'}}\right|_{2,\mathbb{P}}^{2} & \precsim\delta^{2}\left[1\vee\sum_{k=0}^{\gamma_{1}^{*}}\frac{R_{2k}}{\left(2k\right)!}\right]^{2}\precsim\delta^{2}\textrm{ }\quad\textrm{for any }j^{'},l^{'}\in\mathcal{M}^{+},j^{'}\neq l^{'}
\end{align*}
where $2\gamma_{1}^{*}=\gamma^{*}$ (if $\gamma^{*}$ is even) and
$2\gamma_{1}^{*}+1=\gamma^{*}$ (if $\gamma^{*}$ is odd), and the
last inequality follows from that $R_{0}=\frac{\overline{C}}{6}$
and $R_{2k}=\frac{\overline{C}}{6\gamma_{1}^{*}}$ for $k=1,...,\gamma_{1}^{*}$.
Likely previously, by choosing $\delta^{2}\asymp\frac{\sigma^{2}\left(\gamma^{*}+1\right)}{n}$,
(\ref{eq:KL_hetero}) implies 
\begin{align*}
\frac{1}{\left(MM^{+}\right)^{2}}\sum_{\begin{array}{c}
j,l\in\left\{ 1,...,M\right\} \\
j^{'},l^{'}\in\left\{ 1,...,M^{+}\right\} 
\end{array}}D_{KL}\left(\mathbb{P}_{f^{j},\sigma^{j^{'}}}\times\mathbb{P}_{X}\parallel\mathbb{P}_{f^{l},\sigma^{l^{'}}}\times\mathbb{P}_{X}\right) & \leq\underline{C}_{0}^{'}\left(\gamma^{*}+1\right),\\
\log(MM^{+})=\log\left|\mathcal{M}\right|+\log\left|\mathcal{M}^{+}\right| & \geq\overline{C}_{0}^{'}\left(\gamma^{*}+1\right),
\end{align*}
where the positive universal constants satisfy $\underline{C}_{0}^{'}\leq\frac{1}{2}\overline{C}_{0}^{'}$.
Therefore, 
\[
\delta^{2}\left(1-\frac{\log2+\frac{1}{M^{2}M^{+2}}\sum_{j,l\in\left\{ 1,...,M\right\} ,j^{'},l^{'}\in\left\{ 1,...,M^{+}\right\} }D_{KL}\left(\mathbb{P}_{f^{j},\sigma^{j^{'}}}\times\mathbb{P}_{X}\parallel\mathbb{P}_{f^{l},\sigma^{l^{'}}}\times\mathbb{P}_{X}\right)}{\log\left(MM^{+}\right)}\right)\succsim\frac{\sigma^{2}\left(\gamma^{*}+1\right)}{n}.
\]

\subsection{Lemma B.1 and its proof\label{subsec:Proof-for-Lemma B.1}}

\textbf{Lemma B.1}. \textit{We have 
\begin{eqnarray*}
\log M_{2}\left(\delta,\,\mathcal{S}_{\gamma+1}\right)\geq\log M_{2}\left(\delta,\,\overline{\mathcal{S}}_{\gamma+1}\right) & \succsim & \gamma+1,\:\forall\frac{1}{\delta}\succsim2^{\gamma+1}.
\end{eqnarray*}
If the density function $p(x)$ on $\left[0,\,1\right]$ is bounded
away from zero, i.e., $p(x)\geq c>0$, then
\begin{eqnarray}
\log M_{2,\mathbb{P}}\left(\delta,\,\mathcal{S}_{\gamma+1}\right)\geq\log M_{2,\mathbb{P}}\left(\delta,\,\overline{\mathcal{S}}_{\gamma+1}\right) & \succsim & \gamma+1,\:\forall\frac{1}{\delta}\succsim2^{\gamma+1}.\label{eq:sobolev_lower_packing}
\end{eqnarray}
}\textbf{Proof}. To prove Lemma B.1, we use the bound $\underline{B}_{2}$
in Lemma \ref{lm:entropy-generalized-polynomial} (and Remark A at
the end of Appendix \ref{subsec:Proof-for-entropy-holder-sub}). Let
us consider
\begin{equation}
\mathcal{U}_{\gamma+1,1}^{*}=\left\{ f=\sum_{k=0}^{\gamma}\theta_{k}x^{k}:\textrm{\ensuremath{\left(\theta_{k}\right)_{k=0}^{\gamma}}\ensuremath{\ensuremath{\in}}\ensuremath{\mathcal{P}_{\gamma}^{*}}},\,x\in\left[0,\,1\right]\right\} \label{eq:polynomial subspace in Sobolev}
\end{equation}
with the $\left(\gamma+1\right)-$dimensional polyhedron 
\[
\ensuremath{\mathcal{P}_{\gamma}^{*}=\left\{ \ensuremath{\left(\theta_{k}\right)_{k=0}^{\gamma}}:\theta_{0}\in\left[-\frac{\overline{C}}{6},\,\frac{\overline{C}}{6}\right],\,\theta_{k}\in\left[-\frac{\overline{C}}{6\gamma}\frac{1}{k!},\,\frac{\overline{C}}{6\gamma}\frac{1}{k!}\right],\,k=1,...,\gamma\right\} }.
\]
Note that 
\begin{equation}
\mathcal{U}_{\gamma+1,1}^{*}\subset\overline{\mathcal{S}}_{\gamma+1}\subset\mathcal{S}_{\gamma+1}.\label{eq:contain}
\end{equation}
As a result, we can apply (\ref{eq:10}) in Lemma \ref{lm:entropy-generalized-polynomial}
to show \textit{
\begin{eqnarray*}
\log M_{2}\left(\delta,\,\mathcal{U}_{\gamma+1,1}^{*}\right) & \succsim & \gamma+1,\:\forall\frac{1}{\delta}\succsim2^{\gamma+1}.
\end{eqnarray*}
}If the density function $p(x)$ on $\left[0,\,1\right]$ is bounded
away from zero, i.e., $p(x)\geq c>0$, then 
\begin{eqnarray*}
\log M_{2,\mathbb{P}}\left(\delta,\,\overline{\mathcal{S}}_{\gamma+1}\right) & \asymp & \log M_{2}\left(\delta,\,\overline{\mathcal{S}}_{\gamma+1}\right).
\end{eqnarray*}
Therefore we have (\ref{eq:sobolev_lower_packing}).

\subsection{Proof for Theorem \ref{thm:MISE-upper-sobolev-std}\label{subsec:Proof-for-MISE-upper-sobolev-std} }

\textbf{Step 1}. We apply Lemma C.8 (in Appendix \ref{sec:Supporting-lemmas})
where $\mathcal{W}$ corresponds to the Sobolev space containing $\mathcal{S}_{k+1}$
($k=\gamma$ in part (i) and $k=\gamma^{*}$ in part (ii)) and the
kernel functions correspond to (\ref{eq:kernel_1}) and (\ref{eq:kernel_2}).
Then solving (\ref{eq:99}) is reduced to solving
\[
\left(r\sqrt{\frac{\left(k+1\right)\wedge n}{n}}\right)\vee\left(\frac{1}{\sqrt{n}}r^{\frac{2k+1}{2k+2}}\right)\asymp\frac{\overline{C}r^{2}}{\sigma}.
\]
Under the condition $n\geq k+1$, $r=\bar{\delta}_{1}=c_{1}\left[\sqrt{\frac{\sigma^{2}(k+1)}{n}}\vee\left(\frac{\sigma^{2}}{n}\right)^{\frac{k+1}{2k+3}}\right]$
solves the above. In a similar fashion, $r=\bar{\delta}_{2}=c_{2}\left[\sqrt{\frac{k+1}{n}}\vee\left(\frac{1}{n}\right)^{\frac{k+1}{2k+3}}\right]$
solves (\ref{eq:100}). Note that both $\bar{\delta}_{1}$ and $\bar{\delta}_{2}$
are non-random and do not depend on the values of $\left\{ x_{i}\right\} _{i=1}^{n}$.\\
\\
\textbf{Step 2}. Given $\bar{\delta}_{1}$, we apply Lemma C.6 to
show that 
\begin{equation}
\left|\hat{f}-f\right|_{n}^{2}\precsim t_{1}^{2}\quad\textrm{for any }t_{1}\geq\bar{\delta}_{1}\label{eq:70}
\end{equation}
with probability at least $1-c^{'}\exp\left(-c^{''}\frac{nt_{1}^{2}}{\sigma^{2}}\right)$,
whenever $\lambda\asymp t_{1}^{2}$. By choosing $\lambda\asymp t_{1}^{2}=\bar{\delta}_{1}^{2}$,
we have $\lambda\asymp\left[\frac{k+1}{n}\vee\left(\frac{1}{n}\right)^{\frac{2k+2}{2k+3}}\right]$
since $\sigma\asymp1$.\\
\\
\textbf{Step 3}. Given $\bar{\delta}_{2}$, we now connect $\left|\hat{f}-f\right|_{n}^{2}$
with $\left|\hat{f}-f\right|_{2,\mathbb{P}}^{2}$. We divide the argument
into two cases depending on whether $\bar{\delta}_{2}\geq r^{*}$,
the smallest positive solution to (\ref{eq:crit_radius_pop}) with
$c=\overline{C}$.

Case 1 (when $\bar{\delta}_{2}\leq r^{*}$). Note that $\bar{\delta}_{2}$
is an upper bound for $\hat{r}^{*}$, the smallest positive solution
to (\ref{eq:crit_radius_emp}) with $c_{0}=\overline{C}^{-1}$. For
case 1, we can apply Lemma C.4 together with (\ref{eq:70}) to show
that 
\[
\left|\hat{f}-f\right|_{2,\mathbb{P}}^{2}\precsim t_{1}^{2}+t_{2}^{2}\quad\textrm{for any }t_{1}\geq\bar{\delta}_{1},t_{2}\in[\bar{\delta}_{2},\,r^{*}]
\]
with probability at least \textit{
\[
1-c^{'}\exp\left(-c^{''}\frac{nt_{1}^{2}}{\sigma^{2}}\right)-c_{0}^{'}\exp\left(-c_{0}^{''}nr^{*2}\right),
\]
}which is greater than 
\begin{equation}
1-c^{'}\exp\left(-c^{''}\frac{nt_{1}^{2}}{\sigma^{2}}\right)-c_{0}^{'}\exp\left(-c_{0}^{''}nt_{2}^{2}\right).\label{eq:tail}
\end{equation}

Case 2 (when $\bar{\delta}_{2}>r^{*}$). In this case, we can apply
Lemma C.2 (where we take $\bar{r}=r^{*}$) together with (\ref{eq:70})
to show that 
\[
\left|\hat{f}-f\right|_{2,\mathbb{P}}^{2}\precsim t_{1}^{2}+t_{2}^{2}\quad\textrm{for any }t_{1}\geq\bar{\delta}_{1},t_{2}\geq\bar{\delta}_{2}
\]
with probability at least (\ref{eq:tail}). 

Applying Lemmas C.2 and C.4 requires the shifted class $\bar{\mathcal{F}}$
associated with $\mathcal{F}=\mathcal{S}_{k+1}$ to be a class such
that for all $g\in\bar{\mathcal{F}}$, $\left|g\right|_{\infty}\precsim1$.
This condition holds easily given the kernel functions (\ref{eq:kernel_1})
and (\ref{eq:kernel_2}) and Lemma C.9.\\
\textbf{}\\
\textbf{Step 4}. Integrating the tail probability in the form of (\ref{eq:tail})
gives 
\begin{eqnarray*}
\mathbb{E}\left(\left|\hat{f}-f\right|_{2,\mathbb{P}}^{2}\right) & \precsim & \check{r}^{2}+\exp\left\{ -cn\check{r}^{2}\right\} ,
\end{eqnarray*}
where $\check{r}^{2}=\frac{\sigma^{2}(k+1)}{n}\vee\left(\frac{\sigma^{2}}{n}\right)^{\frac{2\left(k+1\right)}{2\left(k+1\right)+1}}\asymp\frac{k+1}{n}\vee\left(\frac{1}{n}\right)^{\frac{2\left(k+1\right)}{2\left(k+1\right)+1}}$
since $\sigma\asymp1$. Finally, we take sup and obtain 
\begin{eqnarray*}
\sup_{f\in\overline{\mathcal{S}}_{\gamma+1}}\mathbb{E}\left(\left|\hat{f}-f\right|_{2,\mathbb{P}}^{2}\right) & \precsim & \check{r}^{2}+\exp\left\{ -cn\check{r}^{2}\right\} .
\end{eqnarray*}

For part (i), we have $\lambda\asymp\left(\frac{1}{n}\right)^{\frac{2k+2}{2k+3}}$
and
\begin{eqnarray*}
\sup_{f\in\overline{\mathcal{S}}_{\gamma+1}}\mathbb{E}\left(\left|\hat{f}-f\right|_{2,\mathbb{P}}^{2}\right) & \leq & \overline{c}\left[r_{1}^{2}+\exp\left\{ -cnr_{1}^{2}\right\} \right]
\end{eqnarray*}
for some universal constant $\overline{c}\in\left(1,\,\infty\right)$,
where $r_{1}^{2}=\left(\frac{\sigma^{2}}{n}\right)^{\frac{2\left(\gamma+1\right)}{2\left(\gamma+1\right)+1}}$. 

For part (ii), we have $\lambda\asymp\frac{k+1}{n}$ and
\begin{equation}
\sup_{f\in\overline{\mathcal{S}}_{\gamma+1}}\mathbb{E}\left(\left|\hat{f}-f\right|_{2,\mathbb{P}}^{2}\right)\leq\overline{c}\left[r_{2}^{2}+\exp\left\{ -cnr_{2}^{2}\right\} \right]\label{eq:upper_r2}
\end{equation}
for some universal constant $\overline{c}\in\left(1,\,\infty\right)$,
where $r_{2}^{2}=\frac{\sigma^{2}\left(\gamma^{*}+1\right)}{n}$. 

\subsection{Proof for Theorem \ref{thm:MISE-upper-holder-std}\label{subsec:Proof-for-MISE-upper-holder-std} }

\textbf{Step 1}. We apply Lemma C.7 (in Appendix \ref{sec:Supporting-lemmas})
where $\mathcal{F}$ corresponds to the standard $\mathcal{U}_{k+1}$
($k=\gamma$ in part (i) and $k=\gamma^{*}$ in part (ii)). Taking
$R_{j}=\overline{C}$ for $j=0,...,k+1$ in (\ref{eq:upper_together})
(the second bound) yields 
\begin{equation}
\log N_{\infty}\left(\delta,\,\mathcal{U}_{k+1}\right)\precsim\left(k+1\right)\log\frac{1}{\delta}+\left(\frac{1}{\delta}\right)^{\frac{1}{k+1}},\label{eq:21}
\end{equation}
where we have used the fact that $R^{*}\asymp1$. Note that 
\begin{eqnarray*}
 &  & \frac{1}{\sqrt{n}}\int_{0}^{r}\sqrt{\log N_{n}(\delta,\,\Omega(r;\,\bar{\mathcal{F}}))}d\delta\\
 & \leq & \frac{1}{\sqrt{n}}\int_{0}^{r}\sqrt{\log N_{\infty}\left(\delta,\,\bar{\mathcal{F}}\right)}d\delta\\
 & \precsim & \underset{T(r)}{\underbrace{r\sqrt{\frac{k+1}{n}}+\frac{1}{\sqrt{n}}r^{\frac{2k+1}{2k+2}}}}
\end{eqnarray*}
where the last line follows from (\ref{eq:21}). Setting $\sigma T(r)\asymp r^{2}$
yields $r=\bar{\delta}_{1}=c_{1}\left[\sqrt{\frac{\sigma^{2}(k+1)}{n}}\vee\left(\frac{\sigma^{2}}{n}\right)^{\frac{k+1}{2k+3}}\right]$,
which solves (\ref{eq:97}). In a similar fashion, $r=\bar{\delta}_{2}=c_{2}\left[\sqrt{\frac{k+1}{n}}\vee\left(\frac{1}{n}\right)^{\frac{k+1}{2k+3}}\right]$
solves (\ref{eq:98}). Note that both $\bar{\delta}_{1}$ and $\bar{\delta}_{2}$
are non-random and do not depend on the values of $\left\{ x_{i}\right\} _{i=1}^{n}$.\textbf{}\\
\textbf{}\\
\textbf{Step 2}. Given $\bar{\delta}_{1}$, we apply Lemma C.5 to
show that 
\begin{equation}
\left|\hat{f}-f\right|_{n}^{2}\precsim t_{1}^{2}\quad\textrm{for any }t_{1}\geq\bar{\delta}_{1}\label{eq:70-1}
\end{equation}
with probability at least $1-c^{'}\exp\left(-c^{''}\frac{nt_{1}^{2}}{\sigma^{2}}\right)$.
\\
\\
\textbf{Steps 3 and 4}. The arguments are identical to those in Step
3 and Step 4 of the proof for Theorem \ref{thm:MISE-upper-sobolev-std}
in Appendix \ref{subsec:Proof-for-MISE-upper-sobolev-std}. The verification
that $\left|g\right|_{\infty}\precsim1$ for all $g\in\bar{\mathcal{F}}$
associated with the standard $\mathcal{U}_{k+1}$ is obvious given
its definition.

\subsection{Proof for Corollary \ref{corr:analytic functions}\label{subsec:Proof-for-analytic_MISE}}

The proof for the lower bound in Corollary \ref{corr:analytic functions}
is identical to Appendix \ref{subsec:Proof-for-MISE-lower-standard}.
For the upper bound associated with $\mathcal{S}_{\infty}$ in Corollary
\ref{corr:analytic functions}, we replace $\bar{\delta}_{1}=c_{1}\left[\sqrt{\frac{\sigma^{2}(k+1)}{n}}\vee\left(\frac{\sigma^{2}}{n}\right)^{\frac{k+1}{2k+3}}\right]$
and $\bar{\delta}_{2}=c_{2}\left[\sqrt{\frac{k+1}{n}}\vee\left(\frac{1}{n}\right)^{\frac{k+1}{2k+3}}\right]$
in Appendix \ref{subsec:Proof-for-MISE-upper-sobolev-std} with $\bar{\delta}_{1}=c_{1}\sqrt{\frac{\sigma^{2}(k+1)}{n}}$
and $\bar{\delta}_{2}=c_{2}\sqrt{\frac{k+1}{n}}$. Then Step 4 of
Appendix \ref{subsec:Proof-for-MISE-upper-sobolev-std} is reduced
to (\ref{eq:upper_r2}) only. Similar modifications are made in Appendix
\ref{subsec:Proof-for-MISE-upper-holder-std} for the derivation of
the upper bound associated with $\mathcal{U}_{\infty}$.

\subsection{Proof for Theorem \ref{thm:MISE-gen-sub-ellipsoid}\label{subsec:Proof-for-MISE-gen-sub-ellipsoid}}

To show the lower bounds in Theorem \ref{thm:MISE-gen-sub-ellipsoid},
we can apply either Lemma C.1(i) or Lemma C.1(ii) in Appendix \ref{sec:Supporting-lemmas},
but the latter gives more insight about where the rates in Theorem
\ref{thm:MISE-gen-sub-ellipsoid} are coming from. The arguments for
the upper bound in Claim (i) of Theorem \ref{thm:MISE-gen-sub-ellipsoid}
are similar to those in Appendix \ref{subsec:Proof-for-MISE-upper-sobolev-std}.
The arguments for the upper bound in Claim (ii) are similar to those
in Appendix \ref{subsec:Proof-for-MISE-upper-holder-std}.\textbf{}\\
\textbf{}\\
\textbf{The lower bound} \textbf{(Sobolev)}. We apply Lemma C.1(ii)
in Appendix \ref{sec:Supporting-lemmas} with $\mathcal{F}=\mathcal{H}_{\gamma+1}$,
the results in Lemma \ref{lm:entropy-ellipsoid-subclass} and (\ref{eq:sandwich}).
By (\ref{eq:cover_KL}) and Lemma \ref{lm:entropy-ellipsoid-subclass},
we have
\begin{align*}
\log N_{KL}\left(\epsilon,\,\mathcal{Q}\right) & \precsim\left(\frac{R_{\gamma+1}\sqrt{n}}{\epsilon}\right)^{\frac{1}{\gamma+1}}.
\end{align*}
Setting $\left(\frac{R_{\gamma+1}\sqrt{n}}{\epsilon}\right)^{\frac{1}{\gamma+1}}\asymp\epsilon^{2}$
yields $\epsilon^{2}\asymp\left(nR_{\gamma+1}^{2}\right)^{\frac{1}{2\left(\gamma+1\right)+1}}=:\epsilon^{*2}$.
Observe that setting 
\[
\delta\asymp R_{\gamma+1}^{\frac{1}{2\left(\gamma+1\right)+1}}\left(\frac{\sigma^{2}}{n}\right)^{\frac{\gamma+1}{2\left(\gamma+1\right)+1}}
\]
and the assumption $\sigma\asymp1$ ensure 
\[
\left(R_{\gamma+1}\delta^{-1}\right)^{\frac{1}{\gamma+1}}\asymp R_{\gamma+1}^{\frac{2}{2\gamma+3}}\left(\frac{n}{\sigma^{2}}\right)^{\frac{1}{2\left(\gamma+1\right)+1}}\succsim\epsilon^{*2}.
\]
Consequently, we have 
\[
1-\frac{\log2+\log N_{KL}\left(\epsilon^{*},\,\mathcal{Q}\right)+\epsilon^{*2}}{\log M_{2,\mathbb{P}}\left(\delta,\,\mathcal{H}_{\gamma+1}\right)}\geq\frac{1}{2}
\]
and 
\[
\inf_{\tilde{f}}\sup_{f\in\mathcal{H}_{\gamma+1}}\mathbb{E}\left(\left|\tilde{f}-f\right|_{2,\mathbb{P}}^{2}\right)\succsim R_{\gamma+1}^{\frac{2}{2\left(\gamma+1\right)+1}}\left(\frac{\sigma^{2}}{n}\right)^{\frac{2\left(\gamma+1\right)}{2\left(\gamma+1\right)+1}}.
\]
\textbf{The upper bound} \textbf{(Sobolev)}. \textbf{Step 1}. We apply
Lemma C.8 where $\mathcal{W}$ corresponds to the RKHS associated
with the kernel function $\mathcal{K}$, which is continuous, positive
semidefinite, and satisfies $\mathcal{K}\left(x,x^{'}\right)\precsim1$
for all $x,\,x^{'}\in\left[0,\,1\right]$. Moreover, $\mathcal{W}$
contains $\mathcal{H}_{\gamma+1}$. Then solving (\ref{eq:99}) is
reduced to 
\[
\frac{1}{\sqrt{n}}r^{\frac{2\gamma+1}{2\gamma+2}}\asymp\frac{R_{\gamma+1}r^{2}}{\sigma}.
\]
Note that $r=\bar{\delta}_{1}=c_{1}R_{\gamma+1}^{\frac{-2\left(\gamma+1\right)}{2\left(\gamma+1\right)+1}}\left(\frac{\sigma^{2}}{n}\right)^{\frac{\gamma+1}{2\gamma+3}}$
solves the above. In a similar fashion, $r=\bar{\delta}_{2}=c_{2}R_{\gamma+1}^{\frac{-2\left(\gamma+1\right)}{2\left(\gamma+1\right)+1}}\left(\frac{1}{n}\right)^{\frac{\gamma+1}{2\gamma+3}}$
solves (\ref{eq:100}). Note that both $\bar{\delta}_{1}$ and $\bar{\delta}_{2}$
are non-random and do not depend on the values of $\left\{ x_{i}\right\} _{i=1}^{n}$.\\
\\
\textbf{Steps 2-4}. The arguments are identical to those in Steps
2-4 in Appendix \ref{subsec:Proof-for-MISE-upper-sobolev-std}. \\
\textbf{}\\
\textbf{The lower bound} \textbf{(Hölder)}. We apply Lemma C.1(ii)
with $\mathcal{F}=\mathcal{U}_{\gamma+1,2}$ and the results in Lemma
\ref{lm:entropy-generalized-H=0000F6lder}. By (\ref{eq:cover_KL})
and Lemma \ref{lm:entropy-generalized-H=0000F6lder}, we have
\begin{align*}
\log N_{KL}\left(\epsilon,\,\mathcal{Q}\right) & \precsim\left(\frac{R^{*}\sqrt{n}}{\epsilon}\right)^{\frac{1}{\gamma+1}}.
\end{align*}
The rest of the arguments are identical to those for the lower bound
concerning $\mathcal{H}_{\gamma+1}$ by simply replacing $R_{\gamma+1}$
with $R^{*}$. \\
\\
\textbf{The upper bound (Hölder)}. \textbf{Step 1}. We apply Lemma
C.7 where $\mathcal{F}$ corresponds to $\mathcal{U}_{\gamma+1,2}$.
Note that 
\begin{eqnarray*}
 &  & \frac{1}{\sqrt{n}}\int_{0}^{r}\sqrt{\log N_{n}(\delta,\,\Omega(r;\,\bar{\mathcal{F}}))}d\delta\\
 & \leq & \frac{1}{\sqrt{n}}\int_{0}^{r}\sqrt{\log N_{\infty}\left(\delta,\,\bar{\mathcal{F}}\right)}d\delta\\
 & \precsim & \underset{T(r)}{\underbrace{\frac{1}{\sqrt{n}}\left(R^{*}\right)^{\frac{1}{2\gamma+2}}r^{\frac{2\gamma+1}{2\gamma+2}}}}
\end{eqnarray*}
where the last line follows from Lemma \ref{lm:entropy-generalized-H=0000F6lder}. 

Setting $\sigma T(r)\asymp r^{2}$ yields $r=\bar{\delta}_{1}=c_{1}\left(R^{*}\right)^{\frac{1}{2\left(\gamma+1\right)+1}}\left(\frac{\sigma^{2}}{n}\right)^{\frac{\gamma+1}{2\left(\gamma+1\right)+1}}$,
which solves (\ref{eq:97}). In a similar fashion, $r=\bar{\delta}_{2}=c_{2}\left(R^{*}\right)^{\frac{1}{2\left(\gamma+1\right)+1}}\left(\frac{1}{n}\right)^{\frac{\gamma+1}{2\left(\gamma+1\right)+1}}$
solves (\ref{eq:98}). Note that both $\bar{\delta}_{1}$ and $\bar{\delta}_{2}$
are non-random and do not depend on the values of $\left\{ x_{i}\right\} _{i=1}^{n}$.
\\
\\
\textbf{Steps 2-4}. The arguments are identical to those in Steps
2-4 in Appendix \ref{subsec:Proof-for-MISE-upper-holder-std}. 

\subsection{Proof for Theorem \ref{thm:MISE_(k-1)!}\label{subsec:Proof-for-MISE_(k-1)!}}

The arguments for the lower bounds are almost identical to those in
Appendix \ref{subsec:Proof-for-MISE-lower-standard}. The arguments
for the upper bounds are almost identical to those in Appendix \ref{subsec:Proof-for-MISE-upper-holder-std}.\textbf{
}In proving Theorem \ref{thm:MISE_(k-1)!}, we use the bounds $\overline{B}_{1}\left(\delta\right)$
and $\underline{B}_{2}$ in Lemma \ref{lm:entropy-generalized-polynomial},
as well as the bounds in Lemma \ref{lm:entropy-generalized-H=0000F6lder}.\textbf{}\\
\textbf{}\\
\textbf{The lower bound}. We apply Lemma C.1(i) in Appendix \ref{sec:Supporting-lemmas}.
Let us consider the packing set consisting of $M_{0}$ elements in
the form (\ref{eq:M-packing}) for $\mathcal{U}_{\gamma^{*}+1,1}$
with $R_{0}=\overline{C}$ and $R_{k}$ taking a value in $\left[\overline{C},\,\overline{C}\left(k-1\right)!\right]$
for all $k=1,...,\gamma^{*}$. Let us choose $b=\gamma\vee1$ in $\theta_{k}^{*}$s
of (\ref{eq:M-packing}) for this construction. Recall that the lower
bound $\log M_{0}\succsim\gamma^{*}+1$ in Appendix \ref{subsec:Proof-for-entropy-poly}
holds for all $\delta$ such that $\frac{R_{\tilde{k}}}{\tilde{k}!\delta\left[\left(\tilde{k}+1\right)\vee\sum_{k=0}^{\gamma^{*}}\frac{R_{k}}{k!}\right]}\geq c2^{\gamma^{*}+1}$
and $\delta\left[\left(\tilde{k}+1\right)\vee\sum_{k=0}^{\gamma^{*}}\frac{R_{k}}{k!}\right]\leq\frac{2}{3}\frac{R_{\tilde{k}}}{\tilde{k}!}$,
where $\tilde{k}=0$ and $c\in\left(1,\,\infty\right)$ is a universal
constant. These conditions are reduced to 
\begin{eqnarray}
\frac{R_{0}}{\delta\left(1\vee\sum_{k=0}^{\gamma^{*}}\frac{R_{k}}{k!}\right)} & \geq & c2^{\gamma^{*}+1},\label{eq:38-1}\\
\delta\left(1\vee\sum_{k=0}^{\gamma^{*}}\frac{R_{k}}{k!}\right) & \leq & \frac{2}{3}R_{0}.\nonumber 
\end{eqnarray}
The rest of the arguments are identical to those in Appendix \ref{subsec:Proof-for-MISE-lower-standard}.
\\
\textbf{}\\
\textbf{The upper bound}. Taking $R_{0}=\overline{C}$ and $R_{j}=\overline{C}\left(j-1\right)!$
for $j=1,...,k+1$ in (\ref{eq:upper_together}) (the first bound)
yields 
\begin{equation}
\log N_{\infty}\left(\delta,\,\mathcal{U}_{k+1}\right)\precsim\left(k+1\right)\log\frac{1}{\delta}+\left(\frac{1}{\delta}\right)^{\frac{1}{k+1}},\label{eq:21-2}
\end{equation}
where we have used the fact that $R^{*}\asymp1$. The rest of the
arguments are identical to those in Appendix \ref{subsec:Proof-for-MISE-upper-holder-std}. 

\subsection{Proof for Theorem \ref{thm:MISE_k!}\label{subsec:Proof-for-MISE_k!}}

In proving Theorem \ref{thm:MISE_k!}, we use the bounds $\overline{B}_{1}\left(\delta\right)$
and $\underline{B}_{1}\left(\delta\right)$ in Lemma \ref{lm:entropy-generalized-polynomial}
for $\mathcal{U}_{k+1,1}$, as well as the bounds in Lemma \ref{lm:entropy-generalized-H=0000F6lder}
for $\mathcal{U}_{k+1,2}$. When $R_{j}=\overline{C}j!$ for all $j=0,...,k+1$,
$R^{*\frac{1}{k+1}}\asymp1$ and $\log M_{2,\mathbb{P}}\left(\delta,\,\mathcal{U}_{k+1,2}\right)\succsim\left(\frac{1}{\delta}\right)^{\frac{1}{k+1}}$.
Moreover,
\begin{equation}
\log N_{\infty}\left(\delta,\,\mathcal{U}_{k+1}\right)\precsim\left(k+1\right)\log\frac{k+1}{\delta}+\left(\frac{1}{\delta}\right)^{\frac{1}{k+1}};\label{eq:32}
\end{equation}
for $k>1$, we have

\begin{eqnarray}
\log M_{2,\mathbb{P}}\left(\delta,\,\mathcal{U}_{k+1,1}\right) & \geq & c^{'}\left[\left(k+1\right)\log\frac{1}{\delta}-k^{2}\right]\nonumber \\
 &  & -\left(k+1\right)\log k^{k}+\left(k+1\right)\log\left(k-1\right)!\nonumber \\
 & \asymp & k\log\frac{1}{\delta}-k^{2};\label{eq:79}
\end{eqnarray}
for $k\in\left\{ 0,\,1\right\} $, we simply have 
\begin{equation}
\log M_{2,\mathbb{P}}\left(\delta,\,\mathcal{U}_{k+1,1}\right)\succsim\left(k+1\right)\log\frac{1}{\delta}-1.\label{eq:80}
\end{equation}
\textbf{The lower bound}. The arguments for the lower bounds are similar
to those in Appendix \ref{subsec:Proof-for-MISE-lower-standard}.
Let us spell out the differences below. 

To show part (ii), we can construct a packing set $\mathcal{M}$ of
$\mathcal{U}_{\gamma^{*}+1,1}$ such that each element in this subset
is $\delta-$apart and by setting $\delta\leq2^{-\left(\gamma^{*}+1\right)}$
in (\ref{eq:79}) and (\ref{eq:80}), we have 
\[
\log\left|\mathcal{M}\right|\succsim\left(\gamma^{*}+1\right)^{2}.
\]
Let us take $\delta^{2}=c\frac{\sigma^{2}\left(\gamma^{*}+1\right)\log\left(\gamma^{*}\vee2\right)}{n}$
for a sufficiently small universal constant $c\in(0,\,1]$. We need
the condition $c4^{\left(\gamma^{*}+1\right)}\left(\gamma^{*}+1\right)\log\left(\gamma^{*}\vee2\right)\leq\frac{n}{\sigma^{2}}$
to ensure $\delta=\sqrt{\frac{c\sigma^{2}\left(\gamma^{*}+1\right)\log\left(\gamma^{*}\vee2\right)}{n}}\leq2^{-\left(\gamma^{*}+1\right)}$
so we can apply $\log\left|\mathcal{M}_{1}\right|\succsim\left(\gamma^{*}+1\right)^{2}\succsim\left(\gamma^{*}+1\right)\log\left(\gamma^{*}\vee2\right)$
below. When $c\in(0,\,1]$ is chosen to be sufficiently small, we
are guaranteed to have
\begin{align*}
\frac{1}{M^{2}}\sum_{j,l\in\left\{ 1,...,M\right\} }D_{KL}\left(\mathbb{P}_{j}\times\mathbb{P}_{X}\parallel\mathbb{P}_{l}\times\mathbb{P}_{X}\right) & \leq\underline{C}_{0}\left(\gamma^{*}+1\right)\log\left(\gamma^{*}\vee2\right),\\
\log M=\log\left|\mathcal{M}_{1}\right|+\log\left|\mathcal{M}_{2}\right| & \geq\overline{C}_{0}\left(\gamma^{*}+1\right)\log\left(\gamma^{*}\vee2\right),
\end{align*}
where the positive universal constants satisfy $\underline{C}_{0}\leq\frac{1}{2}\overline{C}_{0}$.
Therefore,
\[
\delta^{2}\left(1-\frac{\log2+\frac{1}{M^{2}}\sum_{j,l\in\left\{ 1,...,M\right\} }D_{KL}\left(\mathbb{P}_{j}\times\mathbb{P}_{X}\parallel\mathbb{P}_{l}\times\mathbb{P}_{X}\right)}{\log M}\right)\geq\underline{c}\frac{\sigma^{2}\left(\gamma^{*}+1\right)\log\left(\gamma^{*}\vee2\right)}{n}
\]
for some universal constant $\underline{c}\in(0,\,1]$. Note that
we can choose $c$ in $\delta^{2}$ above to be small enough such
that $\underline{c}\leq\frac{1}{3}\underline{c}_{0}$ as stated in
Theorem \ref{thm:MISE_k!}(ii). \\
\textbf{}\\
\textbf{The upper bound}. In this case, solving (\ref{eq:97}) yields
$r=\bar{\delta}_{1}=c_{1}\left[\sqrt{\frac{\sigma^{2}\left(k+1\right)\log\left(k\vee2\right)}{n}}\vee\left(\frac{\sigma^{2}}{n}\right)^{\frac{k+1}{2k+3}}\right]$.
In a similar fashion, solving (\ref{eq:98}) yields $r=\bar{\delta}_{2}=c_{2}\left[\sqrt{\frac{\left(k+1\right)\log\left(k\vee2\right)}{n}}\vee\left(\frac{1}{n}\right)^{\frac{k+1}{2k+3}}\right]$.
Note that both $\bar{\delta}_{1}$ and $\bar{\delta}_{2}$ are non-random
and do not depend on the values of $\left\{ x_{i}\right\} _{i=1}^{n}$.
The rest of the arguments are similar to those in Appendix \ref{subsec:Proof-for-MISE-upper-holder-std}.

\section{Additional supporting lemmas for Appendix \ref{sec:Proofs-for-MISE}
\label{sec:Supporting-lemmas}}

The set of lemmas in this section support our proofs in Appendix \ref{sec:Proofs-for-MISE}
and are based on Dudley (1967), Ledoux and Talagrand (1991), Yang
and Barron (1999), van de Geer (2000), Bartlett and Mendelson (2002),
Mendelson (2002), and Wainwright (2019). Suppose $f$ belongs to the
class $\mathcal{F}$. \\
\textbf{}\\
Lemma C.1 below provides two versions of Fano's inequality. In some
of our derivations of the minimax lower bounds, we apply the first
version as it is more useful for showing the minimax optimal MISE
rates in the small sample regime. In other cases, we apply the second
version as it gives more insight about where the rates are from.\\
\textbf{}\\
\textbf{Lemma C.1}. \textit{(i) Let $\left\{ f^{1},f^{2},...,f^{M}\right\} $
be a $c\delta-$separated set in the $L^{2}\left(\mathbb{P}\right)$
norm of $\mathcal{F}$. Then 
\begin{equation}
\inf_{\tilde{f}_{data}}\sup_{f\in\mathcal{F}}\mathbb{E}\left(\left|\tilde{f}_{data}-f\right|_{2,\mathbb{P}}^{2}\right)\succsim\delta^{2}\left(1-\frac{\log2+\frac{1}{M^{2}}\sum_{j,l\in\left\{ 1,...,M\right\} }D_{KL}\left(\mathbb{P}_{f^{j}}\times\mathbb{P}_{X}\parallel\mathbb{P}_{f^{l}}\times\mathbb{P}_{X}\right)}{\log M}\right)\label{eq:minimax0}
\end{equation}
where $D_{KL}\left(\mathbb{P}_{f^{j}}\times\mathbb{P}_{X}\parallel\mathbb{P}_{f^{l}}\times\mathbb{P}_{X}\right)$
is the $KL-$divergence between the distributions of $\left(Y,\,\left\{ X_{i}\right\} _{i=1}^{n}\right)$
under $f^{j}$ and $f^{l}$, $\mathbb{P}_{X}$ denotes the product
distribution of $\left\{ X_{i}\right\} _{i=1}^{n}$, $\mathbb{P}_{f^{j}}$
denotes the the distribution of $Y$ given $\left\{ X_{i}\right\} _{i=1}^{n}$
when the truth is $f^{j}$, and $data=\left\{ Y_{i},X_{i}\right\} _{i=1}^{n}$.}

\textit{(ii) Let $N_{KL}\left(\epsilon,\,\mathcal{Q}\right)$ denote
the $\epsilon-$covering number of $\mathcal{Q}:=\left\{ \mathbb{P}_{f}\times\mathbb{P}_{X}:f\in\mathcal{F}\right\} $
with respect to the square root of the $KL-$divergence, and $M_{2,\mathbb{P}}\left(\delta,\,\mathcal{F}\right)$
denote the $\delta-$packing number of $\mathcal{F}$ with respect
to $\left|\cdot\right|_{2,\mathbb{P}}$. Then the Yang and Barron
version of Fano's inequality gives}

\textit{
\begin{equation}
\inf_{\tilde{f}_{data}}\sup_{f\in\mathcal{F}}\mathbb{E}\left(\left|\tilde{f}_{data}-f\right|_{2,\mathbb{P}}^{2}\right)\succsim\sup_{\delta,\epsilon}\delta^{2}\left(1-\frac{\log2+\log N_{KL}\left(\epsilon,\,\mathcal{Q}\right)+\epsilon^{2}}{\log M_{2,\mathbb{P}}\left(\delta,\,\mathcal{F}\right)}\right).\label{eq:minimax}
\end{equation}
}\textbf{Remark}. Under Assumption \ref{Assumption 1}, we have 
\begin{eqnarray}
D_{KL}\left(\mathbb{P}_{f^{j}}\times\mathbb{P}_{X}\parallel\mathbb{P}_{f^{l}}\times\mathbb{P}_{X}\right) & = & \mathbb{E}_{X}\left[D_{KL}\left(\mathbb{P}_{f^{j}}\parallel\mathbb{P}_{f^{l}}\right)\right]\nonumber \\
 & = & \frac{n}{2\sigma^{2}}\left|f^{j}-f^{l}\right|_{2,\mathbb{P}}^{2}\label{eq:KL}
\end{eqnarray}
and 
\begin{align}
\log N_{KL}\left(\epsilon,\,\mathcal{Q}\right) & =\log N_{2,\mathbb{P}}\left(\sqrt{\frac{2}{n}}\sigma\epsilon,\,\mathcal{F}\right).\label{eq:cover_KL}
\end{align}

To allow for heteroscedasticity and non-Gaussian noise, let us consider
(\ref{eq:heter model}) where, conditional on $X_{i}=x_{i}$, $Y_{i}$
has a density $\frac{1}{\sigma(x_{i})}p_{w}\left(\frac{y-f(x_{i})}{\sigma(x_{i})}\right)$
subjective to conditions detailed in Zhao and Yang (2022), where $p_{w}(\cdot)$
is the probability density function of the error $w_{i}$. We further
assume that $\mathbb{E}\left(w_{i}\vert X_{i}\right)=0$ for all $i=1,...,n$.
Let $\mathcal{V}$ be the class of functions that $\sigma(\cdot)$
in (\ref{eq:heter model}) ranges over. The class $\mathcal{V}$ consists
of $\sigma(\cdot)$ such that $\sigma(x)\in\left[\underline{\sigma},\,\overline{\sigma}\right]$
for all $x$ in the domain of $\mathcal{V}$ and $\underline{\sigma},\,\overline{\sigma}\asymp1$.
In other words, the scale functions are bounded away from zero and
bounded from above. Then we can apply Lemma 1 in Zhao and Yang (2022)
to arrive at
\begin{eqnarray}
D_{KL}\left(\mathbb{P}_{f^{j},\sigma^{j^{'}}}\times\mathbb{P}_{X}\parallel\mathbb{P}_{f^{l},\sigma^{l^{'}}}\times\mathbb{P}_{X}\right) & \precsim & n\left[\left|f^{j}-f^{l}\right|_{2,\mathbb{P}}^{2}+\left|\sigma^{j^{'}}-\sigma^{l^{'}}\right|_{2,\mathbb{P}}^{2}\right]\label{eq:KL_hetero}
\end{eqnarray}
where $\sigma^{j^{'}},\sigma^{l^{'}}\in\left\{ \sigma^{1},\sigma^{2},...,\sigma^{M^{+}}\right\} $,
a $c\delta-$separated set in the $L^{2}\left(\mathbb{P}\right)$
norm of $\mathcal{V}$. 

In (\ref{eq:KL_hetero}), $D_{KL}\left(\mathbb{P}_{f^{j},\sigma^{j^{'}}}\times\mathbb{P}_{X}\parallel\mathbb{P}_{f^{l},\sigma^{l^{'}}}\times\mathbb{P}_{X}\right)$
is the $KL-$divergence between the distributions of $\left(Y,\,\left\{ X_{i}\right\} _{i=1}^{n}\right)$
under $\left(f^{j},\,\sigma^{j^{'}}\right)$ and $\left(f^{l},\,\sigma^{l^{'}}\right)$,
$\mathbb{P}_{X}$ denotes the product distribution of $\left\{ X_{i}\right\} _{i=1}^{n}$,
and $\mathbb{P}_{f^{j},\sigma^{j^{'}}}$ denotes the the distribution
of $Y$ given $\left\{ X_{i}\right\} _{i=1}^{n}$ when the truth is
$\left(f^{j},\,\sigma^{j^{'}}\right)$. Letting $\mathcal{P}:=\left\{ \mathbb{P}_{h,\sigma}:h\in\mathcal{F},\sigma\in\mathcal{V}\right\} $,
we have\footnote{Note that ``$data$'' (the part of randomness) depends on $\mathbb{P}_{h,\sigma}$,
where $h\in\mathcal{F}$ and $\sigma\in\mathcal{V}$. In stating (\ref{eq:minimax0})
and (\ref{eq:minimax}), the ``sup'' part takes an easier form under
a constant $\sigma$. In other places of the paper, we suppress ``$data$''
and write $\tilde{f}$ for simplicity. } \textit{
\begin{align}
 & \inf_{\tilde{f}_{data}}\sup_{\mathbb{P}_{h,\sigma}\in\mathcal{\mathcal{P}}}\mathbb{E}\left(\left|\tilde{f}_{data}-f(\mathbb{P}_{h,\sigma})\right|_{2,\mathbb{P}}^{2}\right)\nonumber \\
\succsim & \delta^{2}\left(1-\frac{\log2+\frac{1}{M^{2}M^{+2}}\sum_{j,l\in\left\{ 1,...,M\right\} ,j^{'},l^{'}\in\left\{ 1,...,M^{+}\right\} }D_{KL}\left(\mathbb{P}_{f^{j},\sigma^{j^{'}}}\times\mathbb{P}_{X}\parallel\mathbb{P}_{f^{l},\sigma^{l^{'}}}\times\mathbb{P}_{X}\right)}{\log\left(MM^{+}\right)}\right).\label{eq:minimax_hetero_0}
\end{align}
}

Moreover, we have 
\begin{align}
\log N_{KL}\left(\epsilon,\,\mathcal{Q}\right) & \precsim\log N_{2,\mathbb{P}}\left(\sqrt{\frac{1}{n}}\epsilon,\,\mathcal{F}\right)+\log N_{2,\mathbb{P}}\left(\sqrt{\frac{1}{n}}\epsilon,\,\mathcal{V}\right)\label{eq:cover_KL-hetero}
\end{align}
where $\mathcal{Q}:=\left\{ \mathbb{P}_{h,\sigma}\times\mathbb{P}_{X}\right\} $
and \textit{
\begin{equation}
\inf_{\tilde{f}_{data}}\sup_{\mathbb{P}_{h,\sigma}\in\mathcal{\mathcal{P}}}\mathbb{E}\left(\left|\tilde{f}_{data}-f(\mathbb{P}_{h,\sigma})\right|_{2,\mathbb{P}}^{2}\right)\succsim\sup_{\delta,\epsilon}\delta^{2}\left(1-\frac{\log2+\log N_{KL}\left(\epsilon,\,\mathcal{Q}\right)+\epsilon^{2}}{\log M_{2,\mathbb{P}}\left(\delta,\,\mathcal{F}\right)+\log M_{2,\mathbb{P}}\left(\delta,\,\mathcal{V}\right)}\right).\label{eq:minimax_hetero}
\end{equation}
}\textbf{Definition} (local complexity). Let $\Omega_{n}(r;\,\bar{\mathcal{F}})=\left\{ g\in\bar{\mathcal{F}}:\,\left|g\right|_{n}\leq r\right\} $
with 
\begin{equation}
\bar{\mathcal{F}}:=\left\{ g=g_{1}-g_{2}:\,g_{1},\,g_{2}\in\mathcal{F}\right\} .\label{eq:shifted F}
\end{equation}
Conditional on $\left\{ X_{i}\right\} _{i=1}^{n}=\left\{ x_{i}\right\} _{i=1}^{n}$,
the \textit{empirical (conditional) local sub-Gaussian complexity
}is defined as\textit{
\begin{equation}
\mathcal{G}_{n}(r;\,\bar{\mathcal{F}}):=\mathbb{E}_{\tilde{\varepsilon}\vert X}\left[\sup_{g\in\Omega_{n}(r;\,\bar{\mathcal{F}})}\left|\frac{1}{n}\sum_{i=1}^{n}\tilde{\varepsilon}_{i}g(x_{i})\right|\right],\label{eq:13-1}
\end{equation}
}where $\tilde{\varepsilon}=\left\{ \tilde{\varepsilon}_{i}\right\} _{i=1}^{n}$
are independent sub-Gaussian random variables with sub-Gaussian parameters
at most $1$, conditional on $\left\{ X_{i}\right\} _{i=1}^{n}$.
The \textit{empirical local Rademacher complexity $\mathcal{R}_{n}(r;\,\bar{\mathcal{F}})$
}is defined in a similar fashion where $\tilde{\varepsilon}=\left\{ \tilde{\varepsilon}_{i}\right\} _{i=1}^{n}$
are i.i.d. Rademacher variables taking the values of either $-1$
or $1$ equiprobably, and independent of $\left\{ X_{i}\right\} _{i=1}^{n}$. 

Let $\Omega(r;\,\bar{\mathcal{F}})=\left\{ g\in\bar{\mathcal{F}}:\,\left|g\right|_{2,\mathbb{P}}\leq r\right\} $.
The \textit{population local Rademacher complexity }is defined as\textit{
\begin{equation}
\mathcal{R}(r;\,\bar{\mathcal{F}}):=\mathbb{E}_{\tilde{\varepsilon},X}\left[\sup_{g\in\Omega(r;\,\bar{\mathcal{F}})}\left|\frac{1}{n}\sum_{i=1}^{n}\tilde{\varepsilon}_{i}g(X_{i})\right|\right],\label{eq:13-1-1}
\end{equation}
}where $\tilde{\varepsilon}=\left\{ \tilde{\varepsilon}_{i}\right\} _{i=1}^{n}$
are i.i.d. Rademacher variables taking the values of either $-1$
or $1$ equiprobably, and independent of $\left\{ X_{i}\right\} _{i=1}^{n}$.\\
\textbf{}\\
\textbf{Definition} (star-shaped function class).\textbf{ }The class
$\bar{\mathcal{F}}$ is a star-shaped function class if for any $g\in\bar{\mathcal{F}}$
and $\alpha\in[0,\,1]$, $\alpha g\in\bar{\mathcal{F}}$. \textbf{}\\
\\
\textbf{Remark}. The smoothness classes considered in this paper are
star-shaped.\\
\\
\textbf{Lemma C.2}. \textit{Suppose the class $\bar{\mathcal{F}}$
is star-shaped, and for all $g\in\bar{\mathcal{F}}$, $\left|g\right|_{\infty}\leq c$
for some universal constant $c$. Let $\bar{r}$ be any positive solution
to the critical inequality 
\begin{equation}
\mathcal{R}(r;\,\bar{\mathcal{F}})\leq\frac{r^{2}}{c}.\label{eq:crit_radius_pop}
\end{equation}
Then for any $\delta\geq\bar{r}$ and all $g\in\bar{\mathcal{F}}$,
we have 
\begin{equation}
\frac{1}{2}\left|g\right|_{2,\mathbb{P}}^{2}\leq\left|g\right|_{n}^{2}+\frac{\delta^{2}}{2}\label{eq:delta}
\end{equation}
with probability at least $1-c_{1}\exp\left(-c_{2}n\delta^{2}\right)$.
}\\
\\
\textbf{Lemma C.3}. \textit{For any star-shaped $\bar{\mathcal{F}}$,
the function $r\mapsto\frac{\mathcal{R}_{n}(r;\,\mathcal{\bar{\mathcal{F}}})}{r}$
is non-increasing on $(0,\,\infty)$. As a result, the critical inequality
\begin{equation}
\mathcal{R}_{n}(r;\,\bar{\mathcal{F}})\leq c_{0}r^{2}\label{eq:crit_radius_emp}
\end{equation}
has a smallest positive solution for any constant $c_{0}>0$. Similarly,
the function $r\mapsto\frac{\mathcal{R}(r;\,\mathcal{\bar{\mathcal{F}}})}{r}$
is non-increasing on $(0,\,\infty)$. As a result, the critical inequality
\begin{equation}
\mathcal{R}(r;\,\bar{\mathcal{F}})\leq c_{0}r^{2}\label{eq:crit_radius_pop-2}
\end{equation}
has a smallest positive solution for any constant $c_{0}>0$.}\\
\\
\textbf{Lemma C.4}.\textit{ Suppose the class $\bar{\mathcal{F}}$
is star-shaped, and for all $g\in\bar{\mathcal{F}}$, $\left|g\right|_{\infty}\leq c$
for some universal constant $c$. }Let $\hat{r}^{*}$ be the smallest
positive solution to (\ref{eq:crit_radius_emp}) with $c_{0}=c^{-1}$
and $r^{*}$ be the smallest positive solution to (\ref{eq:crit_radius_pop-2})
with $c_{0}=c^{-1}$.\textit{ We have 
\begin{equation}
\left|g\right|_{2,\mathbb{P}}^{2}\precsim\left|g\right|_{n}^{2}+\hat{r}^{*2}\label{eq:84-1}
\end{equation}
with probability at least 
\begin{equation}
1-c_{1}\exp\left(-c_{2}nr^{*2}\right).\label{eq:84}
\end{equation}
}

In the least squares problem (\ref{eq:least squares-1}), if we can
bound $\left|\hat{f}-f\right|_{n}^{2}$ with high probability, then
we can apply Lemmas C.2 or C.4 to bound $\left|\hat{f}-f\right|_{2,\mathbb{P}}^{2}$
with high probability. The following lemmas provide bounds for $\left|\hat{f}-f\right|_{n}^{2}$.
\textbf{}\\
\textbf{}\\
\textbf{Lemma C.5}. \textit{Suppose the class $\bar{\mathcal{F}}$
is star-shaped. Let $\bar{\delta}$ be any positive solution to the
critical inequality 
\begin{equation}
\mathcal{G}_{n}(r;\,\bar{\mathcal{F}})\leq\frac{r^{2}}{2\sigma}.\label{eq:crit_radius_emp_Gaussian}
\end{equation}
Then for any $\delta\geq\bar{\delta}$, we have 
\begin{equation}
\left|\hat{f}-f\right|_{n}^{2}\precsim\delta\bar{\delta}\label{eq:delta-1}
\end{equation}
with probability at least $1-c_{1}\exp\left(-c_{2}\frac{n\delta\bar{\delta}}{\sigma^{2}}\right)$.
}\\

The following lemma concerns the regularized least squares in the
form 
\begin{equation}
\hat{f}\in\arg\min_{\check{f}\in\mathcal{W}}\frac{1}{2n}\sum_{i=1}^{n}\left(y_{i}-\check{f}\left(x_{i}\right)\right)^{2}+\lambda\left|\check{f}\right|_{\mathcal{H}}^{2}\label{eq:regularized least squares}
\end{equation}
where $\mathcal{W}$ is a space of real-valued functions with an associated
semi-norm and contains $\mathcal{F}$. When $\mathcal{W}$ is an RKHS
with its RKHS norm $\left|\cdot\right|_{\mathcal{H}}$, (\ref{eq:regularized least squares})
is referred to as the Kernel Ridge Regression (KRR) estimators. In
particular, as we discuss in Section \ref{sec:Minimax-standard},
when $\mathcal{F}=\mathcal{S}_{k+1}$ in (\ref{eq:least squares-1}),
we can transform (\ref{eq:least squares-1}) into the form (\ref{eq:regularized least squares}),
which is equivalent to solving (\ref{eq:constraints2}) by exploiting
the (reproducing) kernel function associated with the Sobolev space.
To state the following lemma, let us introduce the \textit{empirical
(conditional) local sub-Gaussian complexity }specifically for RKHS: 

\[
\mathcal{G}_{n}(r;\,\bar{\mathcal{W}}):=\mathbb{E}_{\tilde{\varepsilon}\vert X}\left[\sup_{g\in\Omega_{n}(r;\,\bar{\mathcal{W}})}\left|\frac{1}{n}\sum_{i=1}^{n}\tilde{\varepsilon}_{i}g(x_{i})\right|\right],
\]
where 
\[
\Omega_{n}(r;\,\bar{\mathcal{W}})=\left\{ g\in\bar{\mathcal{W}}:\,\left|g\right|_{n}\leq r,\,\left|g\right|_{\mathcal{H}}\leq3\right\} 
\]
and 
\[
\bar{\mathcal{W}}:=\left\{ g=g_{1}-g_{2}:\,g_{1},\,g_{2}\in\mathcal{W}\right\} .
\]
\textbf{Lemma C.6}. \textit{Suppose the class $\mathcal{W}(\supseteq\mathcal{F})$
is convex. Let $\bar{\delta}$ be any positive solution to the critical
inequality 
\begin{equation}
\mathcal{G}_{n}(r;\,\bar{\mathcal{W}})\leq\frac{Rr^{2}}{2\sigma}\label{eq:crit_radius_Gaussian_RKHS}
\end{equation}
where $R$ is a user defined radius. Then for any $\delta\geq\bar{\delta}$,
if (\ref{eq:regularized least squares}) is solved with $\lambda\geq2\delta^{2}$,
we have 
\begin{equation}
\left|\hat{f}-f\right|_{n}^{2}\precsim R^{2}\delta^{2}+R^{2}\lambda\label{eq:delta-1-1}
\end{equation}
with probability at least $1-c_{1}\exp\left(-c_{2}\frac{nR^{2}\delta^{2}}{\sigma^{2}}\right)$.
}\\
\textbf{}\\
\textbf{Remark}.\textbf{ }Concerning the problem in Section \ref{sec:Minimax-standard},
$\mathcal{W}$ corresponds to the Sobolev space, which is convex and
contains $\mathcal{F}$. Moreover, when $\mathcal{F}=\mathcal{S}_{k+1}$
in (\ref{eq:least squares-1}), we can take $R=\overline{C}$.\\

In order to make use of Lemmas C.5 and C.6 to establish sharp bounds
on $\left|\hat{f}-f\right|_{n}^{2}$, we need good candidates for
$\bar{\delta}$ that solves (\ref{eq:crit_radius_emp_Gaussian}) and
(\ref{eq:crit_radius_Gaussian_RKHS}), respectively. To make use of
Lemmas C.2 and C.4 to connect $\left|\hat{f}-f\right|_{n}^{2}$ with
$\left|\hat{f}-f\right|_{2,\mathbb{P}}^{2}$, we need a good candidate
that solves (\ref{eq:crit_radius_emp}). The following lemmas serve
this purpose. \\
\textbf{}\\
\textbf{Lemma C.7}. \textit{Suppose the class $\bar{\mathcal{F}}$
is star-shaped. Let $N_{n}(\delta,\,\Omega_{n}(r;\,\bar{\mathcal{F}}))$
be the $\delta-$covering number of the set $\Omega_{n}(r;\,\bar{\mathcal{F}})$
in the $\left|\cdot\right|_{n}$ norm.}

\textit{(i) Any $\delta\in(0,\,\sigma]$ that solves 
\begin{equation}
\frac{1}{\sqrt{n}}\int_{0}^{r}\sqrt{\log N_{n}(\delta,\,\Omega_{n}(r;\,\bar{\mathcal{F}}))}d\delta\asymp\frac{r^{2}}{\sigma}\label{eq:97}
\end{equation}
solves (\ref{eq:crit_radius_emp_Gaussian}).}

\textit{(ii) Suppose $\left|g\right|_{\infty}\leq c$ for all $g\in\bar{\mathcal{F}}$.
Then any $\delta>0$ that solves 
\begin{equation}
\frac{1}{\sqrt{n}}\int_{0}^{r}\sqrt{\log N_{n}(\delta,\,\Omega_{n}(r;\,\bar{\mathcal{F}}))}d\delta\asymp r^{2}\label{eq:98}
\end{equation}
solves (\ref{eq:crit_radius_emp}).}\\

The following lemma concerns the KRR estimator (\ref{eq:regularized least squares})
when $\mathcal{W}$ in Lemma C.6 is an RKHS.\\
\textbf{}\\
\textbf{Lemma C.8}. \textit{Suppose $\mathcal{W}$ is a convex RKHS
and the KRR estimator (\ref{eq:regularized least squares}) is of
interest. Let $\tilde{\mu}_{1}\geq\tilde{\mu}_{2}\geq...\geq\tilde{\mu}_{n}\geq0$
be the eigenvalues of the kernel matrix $\mathbb{K}\in\mathbb{R}^{n\times n}$
consisting of entries $\frac{1}{n}\mathcal{K}\left(x_{i},\,x_{j}\right)$,
where $\mathcal{K}$ is the kernel function associated with $\mathcal{W}$.
Suppose $\mathcal{K}:\,\mathcal{X}\times\mathcal{X}\rightarrow\mathbb{R}$
is a positive semidefinite kernel function such that $\mathcal{K}\left(x,x^{'}\right)\precsim1$
for all $x,\,x^{'}\in\mathcal{X}$.}

\textit{(i) Any $\delta>0$ that solves 
\begin{equation}
\sqrt{\frac{1}{n}\sum_{i=1}^{n}\min\left\{ r^{2},\,\tilde{\mu}_{i}\right\} }\asymp\frac{Rr^{2}}{\sigma}\label{eq:99}
\end{equation}
solves (\ref{eq:crit_radius_Gaussian_RKHS}). }

\textit{(ii) Any $\delta>0$ that solves 
\begin{equation}
\sqrt{\frac{1}{n}\sum_{i=1}^{n}\min\left\{ r^{2},\,\tilde{\mu}_{i}\right\} }\asymp r^{2}\label{eq:100}
\end{equation}
solves (\ref{eq:crit_radius_emp}) where $c_{0}\asymp1$.}\\

The following lemma is useful when applying Lemmas C.2 anc C.4 in
the case of RKHS and KRR. \\
\\
\textbf{Lemma C.9}. \textit{Let $\mathcal{K}:\,\mathcal{X}\times\mathcal{X}\rightarrow\mathbb{R}$
be a positive semidefinite kernel function such that $\mathcal{K}\left(x,x^{'}\right)\precsim1$
for all $x,\,x^{'}\in\mathcal{X}$. Then, $\left|f\right|_{\infty}\precsim1$
for any function f in the ball of the associated RKHS where the ball
has a constant radius (with respect to the RKHS norm). }

\section{Proofs for Appendix \ref{sec:multi-dim-extensions}\label{sec:Proofs-for-multi-dim-extension}}

\subsection{Proof for Lemma \ref{lm:entropy_multi_dim_Holder_sub} \label{subsec:Proof-for-entropy-multi-dim-holder-sub}}

Like in Appendix \ref{subsec:Proof-for-entropy-holder-sub}, the proper
choice of the grid of points on each dimension of $\left[-1,\,1\right]^{d}$
is the key in this case. Any function $f\in\mathcal{U}_{\gamma+1,2}^{d}$
can be written as 
\begin{align*}
f(x+\Delta) & =\sum_{k=0}^{\gamma}\sum_{p:P=k}\frac{\Delta^{p}D^{p}f\left(x\right)}{k!}+\\
 & \underset{:=REM_{0}(x+\Delta)}{\underbrace{\sum_{p:P=\gamma}\left[\frac{\Delta^{p}D^{p}f\left(z\right)}{\gamma!}-\frac{\Delta^{p}D^{p}f\left(x\right)}{\gamma!}\right]}}
\end{align*}
where $x,\,x+\Delta\in\left(-1,\,1\right)^{d}$ and $z$ is some intermediate
value. For a given $k\in\left\{ 0,...,\gamma\right\} $, recall $\textrm{card}\left(\left\{ p:P=k\right\} \right)=\left(\begin{array}{c}
d+k-1\\
d-1
\end{array}\right)=D_{k}^{*}$. Therefore, we have 
\begin{equation}
\left|REM_{0}(x+\Delta)\right|\leq\frac{D_{\gamma}^{*}R_{\gamma+1}\left|\Delta\right|_{\infty}^{\gamma+1}}{\gamma!}.\label{eq:45}
\end{equation}
In a similar way, writing
\begin{align*}
 & D^{\tilde{p}}f(x+\Delta)=\sum_{k=0}^{\gamma-\tilde{P}}\sum_{p:P=k}\frac{\Delta^{p}D^{p+\tilde{p}}f\left(x\right)}{k!}+\\
 & \underset{:=REM_{\tilde{P}}(x+\Delta)}{\underbrace{\sum_{p:P=\gamma-\tilde{P}}\left[\frac{\Delta^{p}D^{p+\tilde{p}}f\left(\tilde{z}\right)}{\left(\gamma-\tilde{P}\right)!}-\frac{\Delta^{p}D^{p+\tilde{p}}f\left(x\right)}{\left(\gamma-\tilde{P}\right)!}\right]}}
\end{align*}
for $1\leq\tilde{P}:=\sum_{j=1}^{d}\tilde{p}_{j}\leq\gamma$ and $\tilde{p}=\left(\tilde{p}_{j}\right)_{j=1}^{d}$,
we have 
\begin{equation}
\left|REM_{\tilde{P}}(x+\Delta)\right|\leq\frac{D_{\gamma-\tilde{P}}^{*}R_{\gamma-\tilde{P}+1}\left|\Delta\right|_{\infty}^{\gamma+1-\tilde{P}}}{\left(\gamma-\tilde{P}\right)!}.\label{eq:46}
\end{equation}
For some $\delta_{0},\ldots,\delta_{\gamma}>0$, suppose that $\left|D^{p}f\left(w\right)-D^{p}g\left(w\right)\right|\leq\delta_{k}$
for all $p$ with $P=k\in\left\{ 0,\dots,\gamma\right\} $, where
$f,g\in\mathcal{U}_{\gamma+1,2}^{d}$. Then we have 

\begin{eqnarray*}
 &  & \left|f(x+\Delta)-g(x+\Delta)\right|\\
 & \leq & \left|\sum_{k=0}^{\gamma}\sum_{p:P=k}\frac{\Delta^{p}}{k!}\left(D^{p}f\left(x\right)-D^{p}g\left(x\right)\right)\right|+2\frac{D_{\gamma}^{*}R_{\gamma+1}\left|\Delta\right|_{\infty}^{\gamma+1}}{\gamma!}\\
 & \leq & \sum_{k=0}^{\gamma}\frac{D_{k}^{*}\left|\Delta\right|_{\infty}^{k}\delta_{k}}{k!}+2\frac{D_{\gamma}^{*}R_{\gamma+1}\left|\Delta\right|_{\infty}^{\gamma+1}}{\gamma!}.
\end{eqnarray*}
Let $\left(\max_{k\in\left\{ 1,...,\gamma+1\right\} }\frac{D_{k-1}^{*}R_{k}}{\left(k-1\right)!}\vee1\right)=:R^{*}$.
Consider $\left|\Delta\right|_{\infty}\leq d^{-1}\left(R^{*-1}\delta\right)^{\frac{1}{\gamma+1}}$
and $\delta_{k}=R^{*\frac{k}{\gamma+1}}\delta^{1-\frac{k}{\gamma+1}}$
for $k=0,\dots,\gamma$. Then, 

\begin{eqnarray}
\left|f(x+\Delta)-g(x+\Delta)\right| & \leq & \delta\sum_{k=0}^{\gamma}\left(R^{*}{}^{\frac{-k+k}{\gamma+1}}\frac{1}{k!}\right)+2R^{*}\left|\Delta\right|_{\infty}^{\gamma+1}\nonumber \\
 & \leq & \delta\sum_{k=0}^{\gamma}\frac{1}{k!}+2\delta\leq5\delta\label{eq:a90-1}
\end{eqnarray}
where we have used the fact that $D_{k}^{*}\leq d^{k}$. On each dimension
of $\left[-1,\,1\right]^{d}$, we consider a $d^{-1}\left(R^{*-1}\delta\right)^{\frac{1}{\gamma+1}}-$grid
of points. The rest of the arguments follow closely those in Kolmogorov
and Tikhomirov (1959). 

\subsection{Proof for Lemma \ref{lm:entropy_multi_dim_poly} \label{subsec:Proof-for-entropy_multi_poly}}

For a given $k\in\left\{ 0,...,\gamma\right\} $, let $\textrm{card}\left(\left\{ p:P=k\right\} \right)=\left(\begin{array}{c}
d+k-1\\
d-1
\end{array}\right)=\left(\begin{array}{c}
d+k-1\\
k
\end{array}\right)=D_{k}^{*}$. Recall the definition of $\mathcal{U}_{\gamma+1,1}^{d}$:
\[
\mathcal{U}_{\gamma+1,1}^{d}=\left\{ f=\sum_{k=0}^{\gamma}\sum_{p:P=k}x^{p}\theta_{(p,k)}:\{\theta_{(p,k)}\}_{(p,k)}\ensuremath{\in}\ensuremath{\mathcal{P}_{\Gamma}},\,x\in\left[-1,\,1\right]^{d}\right\} 
\]
with the $\Gamma:=\sum_{k=0}^{\gamma}D_{k}^{*}-$dimensional polyhedron
\[
\ensuremath{\mathcal{P}_{\Gamma}=\left\{ \{\theta_{(p,k)}\}_{(p,k)}\ensuremath{\ensuremath{\in}}\mathbb{R}^{\Gamma}:\textrm{for any given }k,\,\{\theta_{(p,k)}\}_{p}\in\left[\frac{-R_{k}}{k!},\,\frac{R_{k}}{k!}\right]\right\} }
\]
where $\theta=\{\theta_{(p,k)}\}_{(p,k)}$ denotes the collection
of $\theta_{(p,k)}$ over all $(p,k)$ configurations and $\{\theta_{(p,k)}\}_{p}$
denotes the collection of $\theta_{(p,k)}$ over all $p$ configurations
for a given $k\in\left\{ 0,...,\gamma\right\} $. 

To bound $\log N_{\infty}\left(\delta,\,\mathcal{U}_{\gamma+1,1}^{d}\right)$
from above, note that for $f,\,f^{'}\in\mathcal{U}_{\gamma+1,1}^{d}$,
we have
\[
\left|f-f^{'}\right|_{\infty}\leq\sum_{k=0}^{\gamma}\sum_{p:P=k}\left|\theta_{(p,k)}-\theta_{(p,k)}^{'}\right|
\]
where $f^{'}=\sum_{k=0}^{\gamma}\sum_{p:P=k}x^{p}\theta_{(p,k)}^{'}$
such that $\theta^{'}=\{\theta_{(p,k)}^{'}\}_{(p,k)}\in\mathcal{P}_{\Gamma}$.
Therefore, the problem is reduced to finding $N_{1}\left(\delta,\,\mathcal{P}_{\Gamma}\right)$. 

To cover $\mathcal{P}_{\Gamma}$ within $\delta-$precision, using
arguments similar to those in Section \ref{subsec:Proof-for-entropy-poly},
we find a smallest $\frac{\delta}{\left(\gamma+1\right)D_{k}^{*}}-$cover
of $\left[\frac{-R_{k}}{k!},\,\frac{R_{k}}{k!}\right]$ for each $k=0,...,\gamma$,
$\left\{ \theta_{k}^{1},...,\theta_{k}^{N_{k}}\right\} $, such that
for any $\theta\in\mathcal{P}_{\Gamma}$, there exists some $i_{(p,k)}\in\left\{ 1,...,N_{k}\right\} $
with 
\[
\sum_{k=0}^{\gamma}\sum_{p:P=k}\left|\theta_{(p,k)}-\theta_{k}^{i_{(p,k)}}\right|\leq\delta.
\]
As a consequence, we have 
\[
\log N_{1}\left(\delta,\,\mathcal{P}_{\Gamma}\right)\leq\sum_{k=0}^{\gamma}D_{k}^{*}\log\frac{4\left(\gamma+1\right)D_{k}^{*}R_{k}}{\delta k!}
\]
and 
\[
\log N_{2,\mathbb{P}}\left(\delta,\,\mathcal{U}_{\gamma+1,1}^{d}\right)\leq\log N_{\infty}\left(\delta,\,\mathcal{U}_{\gamma+1,1}^{d}\right)\leq\sum_{k=0}^{\gamma}D_{k}^{*}\log\frac{4\left(\gamma+1\right)D_{k}^{*}R_{k}}{\delta k!}.
\]
If $\delta$ is large enough such that $\min_{k\in\left\{ 0,...,\gamma\right\} }\log\frac{4\left(\gamma+1\right)D_{k}^{*}R_{k}}{\delta k!}<0$,
we use the counting argument in Kolmogorov and Tikhomirov (1959) to
obtain 
\[
\log N_{2,\mathbb{P}}\left(\delta,\,\mathcal{U}_{\gamma+1,1}^{d}\right)\leq\log N_{\infty}\left(\delta,\,\mathcal{U}_{\gamma+1,1}^{d}\right)\precsim\left(\sum_{k=0}^{\gamma}D_{k}^{*}\right)\log\frac{1}{\delta}+\sum_{k=0}^{\gamma}D_{k}^{*}\log R_{k}.
\]

\end{document}